\RequirePackage{snapshot}

\newif\ifsubsections
\subsectionstrue 

\documentclass[]{pcmi}

\makeatletter
\def\swappedhead#1#2#3{%
  \thmnumber{\@upn{\the\thm@headfont#2\@ifnotempty{#1}{.~}}}%
  \thmname{#1}%
  \thmnote{ {\the\thm@notefont(#3)}}}
\makeatother
\swapnumbers

\let\oldsubsubsection\subsubsection
\renewcommand{\subsubsection}[1]{%
\setcounter{subsubsection}{\value{equation}}%
\oldsubsubsection{#1}
\refstepcounter{equation}
}

\usepackage[sc]{mathpazo}          
\usepackage{eulervm}               
\usepackage[scaled=0.86]{berasans} 
\usepackage[scaled=1]{inconsolata} 
\usepackage[T1]{fontenc}

\usepackage[%
	protrusion=true,
	expansion=false,
	auto=false
	]{microtype}

\usepackage{xcolor}
\usepackage{graphicx}
\graphicspath{{./figures/}}

\ifdraft
	\definecolor{linkred}{rgb}{0.7,0.2,0.2}
	\definecolor{linkblue}{rgb}{0,0.2,0.6}
\else
	\definecolor{linkred}{rgb}{0.0,0.0,0.0}
	\definecolor{linkblue}{rgb}{0,0.0,0.0}
\fi

\usepackage{amssymb, stmaryrd}

\usepackage{mathrsfs}
\input xy
\xyoption{all}
\RequirePackage{xspace}

\PassOptionsToPackage{hyphens}{url} 
\usepackage[
    setpagesize=false,
    pagebackref,
	pdfpagelabels=false,
    pdfstartview={FitH 1000},
    bookmarksnumbered=false,
    linktoc=all,
    colorlinks=true,
    anchorcolor=black,
    menucolor=black,
    runcolor=black,
    filecolor=black,
    linkcolor=linkblue,
	citecolor=linkblue,
	urlcolor=linkred,
]{hyperref}
\usepackage[backrefs,msc-links,nobysame]{amsrefs}

\customizeamsrefs 

\theoremstyle{plain}
\newtheorem{theorem}[equation]{Theorem}
\newtheorem{lemma}[equation]{Lemma}
\newtheorem{prop}[equation]{Proposition}

\theoremstyle{definition}

\newtheorem{remark}[equation]{Remark}

\def\AA{\mathbb{A}}

\def\CC{\mathbb{C}}
\def\DD{\mathbb{D}}

\def\FF{\mathbb{F}}
\def\GG{\mathbb{G}}

\def\NN{\mathbb{N}}

\def\PP{\mathbb{P}}
\def\QQ{\mathbb{Q}}
\def\RR{\mathbb{R}}
\def\SS{\mathbb{S}}

\def\ZZ{\mathbb{Z}}

\def\calA{\mathcal{A}}
\def\calB{\mathcal{B}}

\def\calE{\mathcal{E}}
\def\calF{\mathcal{F}}
\def\calG{\mathcal{G}}
\def\calH{\mathcal{H}}
\def\calI{\mathcal{I}}

\def\calL{\mathcal{L}}
\def\calM{\mathcal{M}}
\def\calN{\mathcal{N}}
\def\calO{\mathcal{O}}
\def\calP{\mathcal{P}}

\def\calS{\mathcal{S}}

\def\calU{\mathcal{U}}
\def\calV{\mathcal{V}}

\long\def\commentout#1{}
\commentout{

\newcommand{\sB}{\ensuremath{\mathscr{B}}\xspace}

\newcommand{\sN}{\ensuremath{\mathscr{N}}\xspace}

\newcommand{\sP}{\ensuremath{\mathscr{P}}\xspace}

\newcommand{\sS}{\ensuremath{\mathscr{S}}\xspace}

\newcommand{\sU}{\ensuremath{\mathscr{U}}\xspace}

\newcommand{\sX}{\ensuremath{\mathscr{X}}\xspace}
\newcommand{\sY}{\ensuremath{\mathscr{Y}}\xspace}

}

\newcommand{\sB}{\ensuremath{\mathcal{B}}\xspace}

\newcommand{\sN}{\ensuremath{\mathcal{N}}\xspace}

\newcommand{\sP}{\ensuremath{\mathcal{P}}\xspace}

\newcommand{\sS}{\ensuremath{\mathcal{S}}\xspace}

\newcommand{\sU}{\ensuremath{\mathcal{U}}\xspace}

\newcommand{\sX}{\ensuremath{\mathcal{X}}\xspace}
\newcommand{\sY}{\ensuremath{\mathcal{Y}}\xspace}

\def\bI{\mathbf{I}}

\def\bP{\mathbf{P}}
\def\bQ{\mathbf{Q}}
\def\bR{\mathbf{R}}

\newcommand\frL{\mathfrak{L}}

\newcommand\frb{\mathfrak{b}}
\newcommand\frc{\mathfrak{c}}

\newcommand\frg{\mathfrak{g}}
\newcommand\frh{\mathfrak{h}}

\newcommand\frl{\mathfrak{l}}

\newcommand\frn{\mathfrak{n}}

\newcommand\frt{\mathfrak{t}}

\newcommand\tilA{\widetilde{A}}

\newcommand\tilC{\widetilde{C}}

\newcommand\tilG{\widetilde{G}}

\newcommand\tilJ{\widetilde{J}}
\newcommand\tilK{\widetilde{K}}

\newcommand\tilR{\widetilde{R}}
\newcommand\tilS{\widetilde{S}}

\newcommand\tilU{\widetilde{U}}
\newcommand\tilV{\widetilde{V}}
\newcommand\tilW{\widetilde{W}}

\newcommand{\Bun}{\textup{Bun}}

\newcommand{\diag}{\textup{diag}}
\DeclareMathOperator\Dyn{Dyn}
\newcommand\ev{\textup{ev}}

\newcommand{\Fl}{\textup{Fl}}

\newcommand\Frac{\textup{Frac}}
\newcommand\Frob{\textup{Frob}}
\newcommand\Gal{\textup{Gal}}

\newcommand{\Gr}{\textup{Gr}}

\newcommand{\Hilb}{\textup{Hilb}}
\newcommand\IC{\textup{IC}}

\renewcommand{\Im}{\textup{Im}}
\newcommand{\Ind}{\textup{Ind}}

\newcommand\inv{\textup{inv}}
\newcommand\Irr{\textup{Irr}}

\newcommand\Lie{\textup{Lie}\ }

\newcommand\nil{\textup{nil}}
\newcommand{\Nm}{\textup{Nm}}

\newcommand{\Pic}{\textup{Pic}}
\newcommand{\cPic}{\overline{\textup{Pic}}}

\newcommand{\red}{\textup{red}}
\newcommand{\reg}{\textup{reg}}

\newcommand{\Res}{\textup{Res}}

\newcommand\rk{\textup{rk}}
\newcommand\rs{\textup{rs}}

\newcommand\Span{\textup{Span}}
\newcommand\Spec{\textup{Spec}\ }
\newcommand\St{\textup{St}}
\newcommand\st{\textup{st}}
\newcommand\Stab{\textup{Stab}}

\newcommand\Sym{\textup{Sym}}

\newcommand\Tot{\textup{Tot}}
\newcommand{\Tr}{\textup{Tr}}

\newcommand{\val}{\textup{val}}

\newcommand{\vol}{\textup{vol}}

\newcommand\Aut{\textup{Aut}}
\newcommand\Hom{\textup{Hom}}
\newcommand\End{\textup{End}}
\newcommand\Isom{\textup{Isom}}

\newcommand\GL{\textup{GL}}
\newcommand\gl{\mathfrak{gl}}
\newcommand\PGL{\textup{PGL}}

\newcommand\SL{\textup{SL}}
\renewcommand\sl{\mathfrak{sl}}

\newcommand\SO{\textup{SO}}

\newcommand\Sp{\textup{Sp}}
\renewcommand\sp{\mathfrak{sp}}

\newcommand{\Gm}{\GG_m}
\def\Ga{\GG_a}

\newcommand{\ad}{\textup{ad}}
\newcommand{\Ad}{\textup{Ad}}

\newcommand\xch{\mathbb{X}^*}
\newcommand\xcoch{\mathbb{X}_*}

\newcommand{\incl}{\hookrightarrow}
\newcommand{\isom}{\stackrel{\sim}{\to}}
\newcommand{\bij}{\leftrightarrow}
\newcommand{\surj}{\twoheadrightarrow}
\newcommand{\ep}{\epsilon}
\renewcommand{\l}{\lambda}
\renewcommand{\L}{\Lambda}
\newcommand{\om}{\omega}
\newcommand{\Om}{\Omega}
\newcommand{\one}{\mathbf{1}}

\newcommand{\Ql}{\QQ_{\ell}}

\newcommand{\twtimes}[1]{\stackrel{#1}{\times}}
\newcommand{\jiao}[1]{\langle{#1}\rangle}
\newcommand{\wt}{\widetilde}
\newcommand{\wh}{\widehat}
\newcommand{\bs}{\backslash}
\newcommand\mat[4]{\left(\begin{array}{cc} #1 & #2 \\ #3 & #4 \end{array}\right)}  
\newcommand\un{\underline}
\newcommand\ov{\overline}
\newcommand\bu{\bullet}
\newcommand\chark{\textup{char}(k)}
\newcommand{\tl}[1]{[\![#1]\!]}
\newcommand{\lr}[1]{(\!(#1)\!)}

\newcommand\Wa{W_{\textup{aff}}}

\newcommand{\homog}[2]{\textup{H}_{#1}({#2})}  
\newcommand{\cohog}[2]{\textup{H}^{#1}({#2})}     

\newcommand\upH{\textup{H}}

\newcommand{\kbar}{\overline{k}}

\newcommand{\Gk}{\Gal(\kbar/k)}

\newcommand{\tN}{\wt{\sN}}
\newcommand{\tX}{\wt{\sX}}
\newcommand{\tY}{\wt{\sY}}
\newcommand{\tB}{\wt{\sB}}
\newcommand{\tg}{\wt{\frg}}
\newcommand{\Ah}{\calA^{\heartsuit}}
\DeclareMathOperator\subreg{subreg}

\newcommand{\Prym}{\textup{Prym}}
\newcommand{\cPrym}{\overline{\textup{Prym}}}
\newcommand{\ur}{\textup{ur}}
\newcommand{\Mpar}{\calM^{\textup{par}}}
\newcommand{\fpar}{f^{\textup{par}}}

\begin{document}

%
%
%
%
%
%

\title{Lectures on Springer theories and orbital integrals}

%
%
\author{Zhiwei Yun}
\thanks{Research supported by the NSF grant DMS-1302071, the Packard Foundation and the PCMI.}
\address{Department of Mathematics, Stanford University, 450 Serra Mall, Building 380, Stanford, CA 94305}
\email{zwyun@stanford.edu}
%
%

\begin{abstract}
These are the expanded lecture notes from the author's mini-course during the graduate summer school of the Park City Math Institute in 2015. The main topics covered are: geometry of Springer fibers, affine Springer fibers and Hitchin fibers; representations of (affine) Weyl groups arising from these objects; relation between affine Springer fibers and orbital integrals.
\end{abstract}

%
%
\maketitle

\tableofcontents

%
%

\section{Introduction}

\subsection{Topics of these notes} These are the lectures notes from a mini-course that I gave at the PCMI graduate summer school in 2015. The goal is twofold. First I would like to introduce to the audience some interesting geometric objects that have representation-theoretic applications, and Springer fibers and their generalizations are nice examples of such. Second is to introduce orbital integrals with emphasis on its relationship with affine Springer fibers, and thereby supplying background material for B-C. Ng\^o's mini-course on the Fundamental Lemmas.

The geometric part of these lectures (everything except \S\ref{s:orb}) consists of the study of three types of ``fibers'': Springer fibers, affine Springer fibers and Hitchin fibers, with increasing complexity. We will study their geometric properties such as connectivity and irreducible components. We will construct certain group actions on these varieties, and use the action to study several nontrivial examples. Most importantly we will study certain Weyl group actions on the cohomology of these fibers which do not come from actions on the varieties themselves. The representation-theoretic significance of these three types of fibers and the analogy between them can be summarized in the following table.

\begin{table*}
	\centering
		\begin{tabular}{|c|c|c|c|}
		\hline
		& Springer fibers & affine  Springer fibers & Hitchin fibers\\
		\hline
		field & $k$ & local field $F=k\lr{t}$ & global field $k(X)$\\
		\hline
		symmetry & $W$ & $\tilW$ &  $\tilW$\\
		\hline
		extended sym& graded AHA & graded DAHA & graded DAHA\\
		\hline
		rep theory & characters of  & orbital integrals &  trace formula \\
		when $k=\FF_{q}$ & $G(k)$ & for $G(F)$ & for $G$ over $k(X)$\\
		\hline
		\end{tabular}
\end{table*}
Here $AHA$ stands for the Affine Hecke Algebra, while DAHA stands for the Double Affine Hecke Algebra; $X$ denotes an algebraic curve over $k$; $W$ and $\tilW$ are the Weyl group and extended affine Weyl group. 

In these lecture notes we do not try to give complete proofs to all statements but instead to point out interesting phenomena and examples. We do, however, give more or less complete proofs of several key results, such as 
\begin{itemize}
\item Theorem \ref{th:Spr corr} (the Springer correspondence);
\item Theorem \ref{th:ft quot} (finiteness properties of affine Springer fibers);
\item Theorem \ref{th:orb coho} (cohomological interpretation of stable orbital integrals).
\end{itemize}

\subsection{What we assume from the readers} The target readers for these lectures are beginning graduate students interested in geometric representation theory. We assume some basic algebraic geometry (scheme theory, moduli problems, point counting over a finite field, etc) though occasionally we will use the language of algebraic stacks. We also assume some Lie theory (reductive groups over an algebraic closed field and over a local field), but knowing $\GL_{n}$ and $\SL_{n}$ should be enough to understand most of these notes. 

The next remark is about the cohomology theory we use in these notes. Since we work with algebraic varieties over a general field $k$ instead of $\CC$, we will be using the \'etale cohomology with coefficients in $\ell$-adic sheaves (usually the constant sheaf $\Ql$) on these varieties. We denote the \'etale cohomology of a scheme $X$ over $k$ with $\Ql$-coefficients simply by $\cohog{*}{X}$. Readers not familiar with \'etale cohomology are encouraged to specialize to the case $k=\CC$ and understand $\cohog{*}{X}$ as the singular cohomology of $X(\CC)$ with $\QQ$-coefficients. Perverse sheaves will be used only in \S\ref{ss:Spr corr}.

\subsection*{Acknowledgement} I would like to thank the co-organizers, lecturers and the staff of the PCMI summer program in 2015. I would also like to thank the audience of my lectures for their feedback. I am especially grateful to Jingren Chi who carefully read through the first draft of these notes and provided helpful suggestions.

\section{Lecture I: Springer fibers}

Springer fibers are classical and fundamental objects in geometric representation theory. Springer \cite{Springer} first discovered that their cohomology groups realized representations of Weyl groups, a phenomenon known as the Springer correspondence. As singular algebraic varieties, Springer fibers are interesting geometric objects by themselves. They are also connected to the representation theory of finite groups of Lie type via character sheaves.

\subsection{The setup}\label{ss:notation}  In this section, let $k$ be an algebraically closed field. Let $G$ be a connected reductive group over $k$ whose adjoint group is simple (so the adjoint group $G^{\ad}$ is determined by one of the seven series of Dynkin diagrams). Assume that $\chark$ is large compared to $G$.  Let $r$ be the rank of $G$.

Sometimes it will be convenient to fix a Cartan subalgebra $\frt$ of $\frg$, or equivalently a maximal torus $T\subset G$. Once we have done this, we may talk about the roots of the $T$-action on $\frg$.
 
Let $\sB$ be the flag variety of $G$: this is the $G$-homogeneous projective variety parametrizing Borel subgroups of $G$. Choosing a Borel subgroup $B\subset G$, we may identify $\sB$ with $G/B$.

Let $\frg$ be the Lie algebra of $G$.  For $X\in\frg$, let $C_{G}(X)$ denote the centralizer of $X$ in $G$, i.e., the stabilizer at $X$ of the adjoint action of $G$ on $\frg$. 

Let $\sN\subset\frg$ be the subvariety of nilpotent elements. This is a cone: it is stable under the action of $\Gm$ on $\frg$ by scaling.  It is known that there are finitely many $G$-orbits on $\sN$ under the adjoint action.

\subsection{Springer fibers}

\subsubsection{The Springer resolution} The cotangent bundle $\tN:=T^{*}\sB$ classifies pairs $(e, B)$ where $e\in\sN$ and $B$ is a Borel subgroup of $G$ such that $e\in \frn$, where $\frn$ is the nilpotent radical of $\Lie B$. The {\em Springer resolution} is the forgetful map
\begin{equation*}
\pi:\tN\to\sN
\end{equation*}
sending $(e, B)$ to $e$. This map is projective.

For $e\in\sN$, the fiber $\sB_{e}:=\pi^{-1}(e)$ is called the {\em Springer fiber} of $e$. By definition, $\sB_{e}$ is the closed subscheme of $\sB$ consisting of those Borel subgroups $B\subset G$ such that $e$ is contained in the nilpotent radical of $\Lie B$. 

\subsubsection{The Grothendieck alteration}\label{sss:Gr alt}
Consider the variety $\tg$ of pairs $(X,B)$ where $X\in\frg$ and $B\in\sB$ such that $X\in\Lie B$. The forgetful map $(X,B)\mapsto X$
\begin{equation*}
\pi_{\frg}: \tg\to \frg
\end{equation*}  
is called the {\em Grothendieck alteration} \footnote{The term ``alteration'' refers to a proper, generically finite map whose source is smooth over $k$.}, also known as the {\em Grothendieck simultaneous resolution}. 
 
Let $\frc=\frg\sslash G\cong \frt\sslash W$ be the categorical quotient of $\frg$ by the adjoint action of $G$. Then we have a commutative diagram
\begin{equation}\label{tgg}
\xymatrix{\tg\ar[d]^{\pi_{\frg}}\ar[r]^{\wt\chi} & \frt\ar[d]\\
\frg \ar[r]^{\chi} & \frc}
\end{equation}
Here $\chi:\frg\to \frc$ is the natural quotient map and $\wt\chi:\tg\to \frt$ sends $(X,B)\in\tg$ to the image of $X$ in $\frb/\frn$ (where $\frb=\Lie B$ with nilpotent radical $\frn$), which can be canonically identified with $\frt$ (upon choosing a Borel containing $\frt$). The diagram \eqref{tgg} is Cartesian when restricting the left column to regular elements  \footnote{An element $X\in\frg$ is called {\em regular} if its centralizer in $G$ has dimension $r$, the rank of $G$.} in $\frg$. In particular, if we restrict the diagram \eqref{tgg} to the regular semisimple locus $\frc^{\rs}\subset \frc$   \footnote{An element $X\in\frt$ is {\em regular semisimple} if $\alpha(X)\neq0$ for any root $\alpha$. Denote by $\frt^{\rs}\subset\frt$ the open subscheme of regular semisimple elements in $\frt$. Since $\frt^{\rs}$ is stable under $W$, it is the full preimage of an open subset $\frc^{\rs}\subset\frc$. The open subset $\frg^{\rs}=\chi^{-1}(\frc^{\rs})$ is by definition the locus of {\em regular semisimple} elements of $\frg$.}, we see that $\pi_{\frg}|_{\tg^{\rs}}:\tg^{\rs}\to \frg^{\rs}$ is a $W$-torsor; i.e., there is an action of $W$ on $\pi_{\frg}^{-1}(\frg^{\rs})$ preserving the projection to $\frg^{\rs}$ making the fibers of $\pi_{\frg}$ into principal homogeneous spaces for $W$. The map $\pi_{\frg}$ becomes branched but still finite over the regular locus $\frg^{\reg}\subset \frg$.

For $X\in \frg$, we denote the fiber $\pi_{\frg}^{-1}(X)$ simply by $\tB_{X}$.  Restricting $\pi_{\frg}$ to $\sN$, the fibers $\tB_{e}$ for a nilpotent element $e\in\sN$ is the closed subvariety of $\sB$ consisting of those $B$ such that $e\in\Lie B$.  Clearly $\sB_{e}$ is a subscheme of $\tB_{e}$, and the two schemes have the same reduced structure. However, as schemes, $\tB_{e}$ and $\sB_{e}$ are different in general. See \S\ref{sss:regular Spr} below and Exercise \ref{exI:nonred Spr}.

\subsection{Examples of Springer fibers} 

\subsubsection{} When $e=0$, $\sB_{e}=\sB$.  

\subsubsection{Regular nilpotent elements}\label{sss:regular Spr} The unique dense open $G$-orbit consists of {\em regular nilpotent elements}, i.e., those $e$ such that $\dim C_{G}(e)=r$. When $e$ is a regular nilpotent element, $\sB_{e}$ is a single point: there is a unique Borel subgroup of $G$ whose Lie algebra contains $e$. What is this Borel subgroup? 

By the Jacobson-Morosov theorem, we may extend $e$ to an $\sl_{2}$-triple $(e,h,f)$ in $\frg$. The adjoint action of $h$ on $\frg$ has integer weights, and it decomposes $\frg$ into weight spaces $\frg(n)$, $n\in\ZZ$. Let $\frb=\oplus_{n\geq 0}\frg(n)$. This is a Borel subalgebra of $\frg$, and the corresponding Borel subgroup is the unique point in $\sB_{e}$.

When $e$ is regular,  the fiber $\tB_{e}$ of the Grothendieck alteration is a non-reduced scheme whose underlying reduced scheme is a point. The coordinate ring of $\tB_{e}$ is isomorphic to the coinvariant algebra $\Sym(\frt^{*})/(\Sym(\frt^{*})^{W}_{+})$, which is, interestingly, also isomorphic to the cohomology ring with $k$-coefficients of the complex flag variety $\calB_{\CC}$.

\subsubsection{} When $G=\SL(V)$ for some vector space $V$ of dimension $n$, $\sB$ is the moduli space of full flags $0=V_{0}\subset V_{1}\subset V_{2}\subset\cdots\subset V_{n-1}\subset V_{n}=V$. The Springer fiber $\sB_{e}$ consists of those flags such that $eV_{i}\subset V_{i-1}$.

\subsubsection{}\label{sss:SL3 subreg} Consider the case $G=\SL_{3}$ and $e=\left(\begin{array}{ccc} 0 & 0 & 1\\  0 & 0& 0 \\  0 & 0& 0 \end{array}\right)$ under the standard basis $\{v_{1}, v_{2}, v_{3}\}$ of $V$. Then $\sB_{e}$ is the union of two $\PP^{1}$s: the first $\PP^{1}$ consisting of flags $0\subset V_{1}\subset \jiao{v_{1}, v_{2}}\subset V$ with varying $V_{1}$ inside the fixed plane $\ker(e)=\jiao{v_{1},v_{2}}$; the second $\PP^{1}$ consisting of flags $0\subset \jiao{v_{1}}\subset V_{2}\subset V$ with a varying $V_{2}$ containing $\Im(e)=\jiao{v_{1}}$.

\subsubsection{} Consider the case $G=\SL_{4}$ and the nilpotent element $e: v_{3}\mapsto v_{1}\mapsto 0, v_{4}\mapsto v_{2}\mapsto 0$ under the standard basis $\{v_{1},v_{2},v_{3},v_{4}\}$. If a flag $0\subset V_{1}\subset V_{2}\subset V_{3}\subset V$ is in $\sB_{e}$, then $V_{1}\subset \ker(e)=\jiao{v_{1},v_{2}}$, and $V_{3}\supset\Im(e)=\jiao{v_{1},v_{2}}$. We denote $H:=\jiao{v_{1},v_{2}}$. There are two cases:
\begin{enumerate}
\item $V_{2}=H$. Then we may choose $V_{1}$ to be any line in $H$ and $V_{3}$ to be any hyperplane containing $V_{2}$. We get a closed subvariety of $\sB_{e}$ isomorphic to $\PP(H)\times \PP(H)\cong \PP^{1}\times \PP^{1}$. We denote this closed subvariety of $\sB_{e}$ by $C_{1}$. 
\item $V_{2}\neq H$. This defines an open subscheme $U$ of $\sB_{e}$. Suppose the line $V_{1}\subset H$ is spanned by  $av_{1}+bv_{2}$ for some $[a:b]\in\PP(H)$, then the image of $V_{2}$ in $V/H=\jiao{v_{3},v_{4}}$ is spanned by $av_{3}+bv_{4}$. Fixing $V_{1}$, the choices of $V_{2}$ are given by $\Hom(\jiao{av_{3}+bv_{4}}, H/V_{1})\cong\Hom(V_{1},H/V_{1})$. Once $V_{2}$ is fixed, $V_{3}=V_{2}+H$ is also fixed. Therefore $U$ is isomorphic to the total space of the line bundle $\calO(2)$ over $\PP(H)\cong \PP^{1}$.
\end{enumerate}
From the above discussion we see that $\sB_{e}$ has dimension $2$, $C_{1}$ is an irreducible component of $\sB_{e}$ and so is the closure of $U$, which we denote by $C_{2}$. We have $C_{1}\cup C_{2}=\sB_{e}$ and $C_{1}\cap C_{2}$ is the diagonal inside $C_{1}\cong\PP^{1}\times\PP^{1}$. 
 
\subsubsection{Components of type $A$ Springer fibers}\label{sss:Spr type A} When $G=\SL_{n}=\SL(V)$, Spaltenstein \cite{Spa76} and Steinberg \cite[]{St} gave a description of the irreducible components of  $\sB_{e}$ using standard Young tableaux of size $n$. This will be relevant to the Springer correspondence that we will discuss later, see \S\ref{sss:A Spr corr}. Below we follow the presentation of \cite[Ch II, \S5]{Spaltenstein}.

Fix a nilpotent element $e\in\sN$ whose Jordan type is a partition $\l$ of $n$. This means, if the partition $\l$ is $n=\l_{1}+\l_{2}+\cdots$, $e$ has Jordan blocks of sizes $\l_{1},\l_{2},\cdots$. We shall construct a (non-algebraic) map $\sB_{e}\to ST(\l)$, where $ST(\l)$ is the discrete set of standard Young tableau for the partition $\l$. For each full flag $0=V_{0}\subset V_{1}\subset\cdots\subset V_{n-1}\subset V_{n}=V$ such that $eV_{i}\subset V_{i-1}$, $e$ induces a nilpotent endomorphism of $V/V_{n-i}$. Let $\mu_{i}$ be the Jordan type of the $e$ on $V/V_{n-i}$, then $\mu_{i}$ is a partition of $i$. The increasing sequence of partitions $\mu_{1},\mu_{2},\cdots, \mu_{n}=\l$ satisfies that $\mu_{i}$ is obtained from $\mu_{i-1}$ by increasing one part of $\mu_{i-1}$ by 1 (including creating a part of size $1$). This gives an increasing sequence $Y_{1},\cdots, Y_{n}=Y(\l)$ of subdiagrams of the Young diagram $Y(\l)$ of $\l$. We label the unique box in $Y_{i}-Y_{i-1}$ by $i$ to get a standard Young tableau.
 
Spaltenstein \cite[Ch II, Prop 5.5]{Spaltenstein} showed that the closure of the preimage of each standard Young tableaux in $\sB_{e}$ is an irreducible component. Moreover, all irreducible components of $\sB_{e}$ arise in this way and they all have the same dimension $d_{e}=\frac{1}{2}\sum_{i}\l^{*}_{i}(\l^{*}_{i}-1)$, where $\l^{*}$ is the conjugate partition of $\l$. In particular, the top dimensional cohomology $\cohog{2d_{e}}{\sB_{e}}$ has dimension equal to $\#ST(\l)$, which is also the dimension of an irreducible representation of the symmetric group $S_{n}$. This statement is a numerical shadow of the Springer correspondence, which says that $\cohog{2d_{e}}{\sB_{e}}$ is naturally an irreducible representation of $S_{n}$. 

Spaltenstein \cite[Ch II, Prop 5.9]{Spaltenstein} also showed that there exists a stratification of $\sB$ into affine spaces such that $\sB_{e}$ is a union of strata. This implies that the restriction map on cohomology $\cohog{*}{\sB}\to \cohog{*}{\sB_{e}}$ is surjective.

\subsubsection{}\label{sss:Sp4} Consider the case $G=\Sp(V)$ for some symplectic vector space $V$ of dimension $2n$, then $\sB$ is the moduli space of full flags
\begin{equation*}
0=V_{0}\subset V_{1}\subset V_{2}\subset\cdots\subset V_{n}\subset\cdots\subset V_{2n-1}\subset V_{2n}=V
\end{equation*}
such that $V^{\bot}_{i}=V_{2n-i}$ for $i=1,\cdots, n$. The Springer fiber $\sB_{e}$ consists of those flags such that $eV_{i}\subset V_{i-1}$ for all $i$.

Consider the case where $\dim V=4$. We choose a basis $\{v_{1},v_{2},v_{3},v_{4}\}$ for $V$ such that the symplectic form $\om$ on $V$ satisfies $\om(v_{i},v_{5-i})=1$ if $i=1$ and $2$, and $\om(v_{i},v_{j})=0$ for $i+j\neq 5$. Let $e$ be the nilpotent element in $\frg=\sp(V)$ given by $e:v_{4}\mapsto v_{1}\mapsto 0, v_{3}\mapsto 0,v_{2}\mapsto 0$. Then a flag $0\subset V_{1}\subset V_{2}\subset V_{1}^{\bot}\subset V$ in $\sB_{e}$ must satisfy $\jiao{v_{1}}\subset V_{2}\subset\jiao{v_{1},v_{2},v_{3}}$, and this is the only condition for it to lie in $\sB_{e}$ (Exercise \ref{exI:211Sp4}). Such a $V_{2}$ corresponds to a line $V_{2}/\jiao{v_{1}}\subset \jiao{v_{2},v_{3}}$, hence a point in $\PP^{1}=\PP(\jiao{v_{2},v_{3}})$. Over this $\PP^{1}$ we have a tautological rank two bundle $\calV_{2}$ whose fiber at $V_{2}/\jiao{v_{1}}$ is the two-dimensional vector space $V_{2}$. The further choice of $V_{1}$ gives a point in the projectivization of $\calV_{2}$. The exact sequence $0\to\jiao{v_{1}}\to V_{2}\to V_{2}/\jiao{v_{1}}\to0$ gives an exact sequence of vector bundles $0\to\calO\to\calV_{2}\to\calO(-1)\to0$ over $\PP^{1}$. Therefore $\calV_{2}$ is isomorphic to $\calO(-1)\oplus\calO$, and $\sB_{e}\cong \PP(\calO(-1)\oplus\calO)$ is a Hirzebruch surface.

\subsubsection{Subregular Springer fibers}\label{sss:subreg} The example considered in \S\ref{sss:SL3 subreg} is a simplest case of a {\em subregular Springer fiber}. There is a unique nilpotent orbit $\calO_{\subreg}$ of codimension 2 in $\sN$, which is called the {\em subregular nilpotent orbit}. For $e\in\calO_{\subreg}$, it is known that $\sB_{e}$ is a union of $\PP^{1}$'s whose configuration we now describe. We may form the dual graph to $\sB_{e}$ whose vertices are the irreducible components of $\sB_{e}$ and two vertices are joined by an edge if the two corresponding components intersect (it turns out that they intersection at a single point).

For simplicity assume $G$ is of adjoint type. Let $G'$ be another adjoint simple group whose type is defined as follows. When $G$ is simply-laced, take $G'=G$. When $G$ is of type $B_{n}$, $C_{n}$, $F_{4}$ and $G_{2}$, take $G'$ to be of type $A_{2n-1}$, $D_{n+1}$, $E_{6}$ and $D_{4}$ respectively. One can show that $\sB_{e}$ is always a union of $\PP^{1}$ whose dual graph is the Dynkin diagram of $G'$. The rule in the non-simply-laced case is that each long simple root corresponds to 2 or 3 $\PP^{1}$'s while each short simple root corresponds to a unique $\PP^{1}$. Such a configuration of $\PP^{1}$'s is called a Dynkin curve, see \cite[\S3.10, Definition and Prop 2]{StConj} and \cite[\S6.3]{Slodowy}.

For example, when $G$ is of type $A_{n}$, then $\sB_{e}$ is a chain consisting of $n$ $\PP^{1}$'s:  $\sB_{e}=C_{1}\cup C_{2}\cup\cdots\cup C_{n}$ with $C_{i}\cap C_{i+1}$ a point and otherwise disjoint.

Brieskorn \cite{Brieskorn}, following suggestions of Grothendieck, related the singularity of the nilpotent cone along the subregular orbits with Kleinian singularities, and he also realized the semi-universal deformation of this singularity inside $\frg$. Assume $G$ is simply-laced. One can construct a transversal slice $S_{e}\subset\sN$ through $e$ of dimension two such that $S_{e}$ consists of regular elements except $e$. Then $S_{e}$ is a normal surface with a Kleinian singularity at $e$ of the same type as the Dynkin diagram of $G$. \footnote{A Kleinian singularity is a surface singularity analytically isomorphic to the singularity at $(0,0)$ of the quotient of $\AA^{2}$ by a finite subgroup of $\SL_{2}$.}.  The preimage $\tilS_{e}:=\pi^{-1}(S_{e})\subset\tN$ turns out to be the minimal resolution of $S_{e}$, and hence $\sB_{e}$ is the union of exceptional divisors. Upon identifying the components of $\sB_{e}$ with simple roots of $G$, the intersection matrix of the exceptional divisors is exactly the negative of the Cartan matrix of $G$. Slodowy \cite{Slodowy} extended the above picture to non-simply-laced groups, and we refer to his book \cite{Slodowy} for a beautiful account of the connection between the subregular orbit and Kleinian singularities.

\subsection{Geometric Properties of Springer fibers}

\subsubsection{Connectivity} The Springer fibers $\sB_{e}$ are connected. See \cite[p.131, Prop 1]{StConj}, \cite[Ch II, Cor 1.7]{Spaltenstein}, and Exercise \ref{exI:conn}.

\subsubsection{Equidimensionality} Spaltenstein \cite{Spa77}, \cite[Ch II, Prop 1.12]{Spaltenstein} showed that all  irreducible components of $\sB_{e}$ have the same dimension. 

\subsubsection{The dimension formula} Let $d_{e}$ be the dimension of $\sB_{e}$. Steinberg \cite[Thm 4.6]{St} and Springer \cite{Springer} showed that
\begin{equation}\label{Spr dim}
d_{e}=\frac{1}{2}(\dim\sN-\dim G\cdot e)=\frac{1}{2}(\dim C_{G}(e)-r).
\end{equation}

\subsubsection{Centralizer action} The group $G\times \Gm$ acts on $\sN$ (where $\Gm$ by dilation). Let $\tilG_{e}=\Stab_{G\times \Gm}(e)$ be the stabilizer. Then $\tilG_{e}$ acts on $\sB_{e}$. Note that $\tilG_{e}$ always surjective onto $\Gm$ with kernel $C_{G}(e)$.

The action of $\tilG_{e}$ on $\cohog{*}{\sB_{e}}$ factors through the finite group $A_{e}=\pi_{0}(\tilG_{e})=\pi_{0}(G_{e})$. Note that $A_{e}$ depends not only on the isogeny class of $G$, but on the isomorphism class of $G$. For example, for $e$ regular, $A_{e}=\pi_{0}(ZG)$ where $ZG$ is the center of $G$. The action of $A_{e}$ on $\cohog{*}{\sB_{e}}$ further factors through the image of $A_{e}$ in the adjoint group $G^{\ad}$.

\subsubsection{Purity} Springer \cite{SpringerPure} proved that the cohomology of $\sB_{e}$ is always pure (in the sense of Hodge theory when $k=\CC$, or in the sense of Frobenius weights when $k$ is a finite field).

\subsubsection{} Let $e\in\sN$. Consider the restriction map $i^{*}_{e}: \cohog{*}{\sB}\to \cohog{*}{\sB_{e}}$ induced by the inclusion $\sB_{e}\incl\sB$. Then the image of $i^{*}_{e}$ is exactly the invariants of $\cohog{*}{\sB_{e}}$ under $A_{e}$. This is a theorem of Hotta and Springer \cite[Theorem 1.1]{HottaSpringer}. In particular, when $G$ is of type $A$, such restriction maps are always surjective.

\subsubsection{Parity vanishing} De Concini, Lusztig and Procesi \cite{DLP} proved that $\cohog{i}{\sB_{e}}$ vanishes for all odd $i$ and any $e\in\sN$. When $k=\CC$, they prove a stronger statement: $\cohog{i}{\sB_{e},\ZZ}$ vanishes for odd $i$ and is torsion-free for even $i$.

\subsection{The Springer correspondence}\label{ss:Spr corr} Let $W$ be the Weyl group of $G$. In 1976, Springer \cite{Springer} made the fundamental observation that there is natural $W$-action on $\cohog{*}{\sB_{e}}$, even though $W$ does not act on $\sB_{e}$ as automorphisms of varieties.

\begin{theorem}[Springer {\cite[Thm 6.10]{Springer}}]\label{th:Spr corr} 
\begin{enumerate}
\item For each nilpotent element $e$, there is a natural graded action of $W$ on $\cohog{*}{\sB_{e}}$ that commutes with the action of $A_{e}$.
\item For each nilpotent element $e$ and each irreducible representation $\rho$ of $A_{e}$, the multiplicity space $M(e,\rho):=\Hom_{A_{e}}(\rho, \cohog{2d_{e}}{\sB_{e}})$ is either zero or an irreducible representation of $W$ under the action in part (1).
\item Each irreducible representation $\chi$ of $W$ appears as $M(e,\rho)$ for a unique pair $(e,\rho)$ up to $G$-conjugacy. The assignment $\chi\mapsto (e,\rho)$ thus gives an injection
\begin{equation}\label{Spr corr}
\Irr(W)\incl \{(e,\rho)\}/G.
\end{equation}
\end{enumerate}
\end{theorem}

\subsubsection{Convention}\label{sss:two actions} In fact there are two natural actions of $W$ on $\cohog{*}{\sB_{e}}$ that differ by tensoring with the sign representation of $W$. In these notes we use the action that is normalized by the following properties. The trivial representation of $W$ corresponds to regular nilpotent $e$ and the trivial $\rho$. The sign representation of $W$ corresponds to $e=0$. Note however that Springer's original paper \cite{Springer} uses the other action.

\subsubsection{The case $e=0$} Taking $e=0$, Springer's theorem gives a graded action of $W$ on $\cohog{*}{\calB}$. What is this action? First, this action can be seen geometrically by considering $G/T$ instead of $\sB=G/B$. In fact, since $N_{G}(T)/T=W$, the right action of $N_{G}(T)$ on $G/T$ induces an action of $W$ on $G/T$, which then induces an action of $W$ on $\cohog{*}{G/T}$. Since the projection $G/T\to G/B$ is an affine space bundle, it follows that $\cohog{*}{\sB}\cong\cohog{*}{G/T}$. It can be shown that under this isomorphism, the action of $W$ on $\cohog{*}{G/T}$ corresponds exactly to Springer's action on $\cohog{*}{\calB}$.

Let $S=\Sym(\xch(T)\otimes\Ql)$ be the graded symmetric algebra where $\xch(T)$ has degree $2$. The reflection representation of $W$ on $\xch(T)$ then induces a graded action of $W$ on $S$. Recall Borel's presentation of the cohomology ring of the flag variety
\begin{equation}\label{Borel HB}
\cohog{*}{\calB,\Ql}\cong S/(S^{W}_{+})
\end{equation}
where $S_{+}\subset S$ is the ideal spanned by elements of positive degree, and $(S^{W}_{+})$ denotes the ideal of $S$ generated by $W$-invariants on $S_{+}$. Then \eqref{Borel HB} is in fact an isomorphism of $W$-modules (see \cite[Prop 7.2]{Springer}). By a theorem of Chevalley, $S/(S^{W}_{+})$ is isomorphic to the regular representation of $W$, therefore, as a $W$-module, $\cohog{*}{\sB}$ is also isomorphic to the regular representation of $W$.  

\begin{remark} The target set in \eqref{Spr corr} can be canonically identified with the set of isomorphism classes of irreducible $G$-equivariant local systems on nilpotent orbits. In fact, for an irreducible $G$-equivariant local systems $\calL$ on a nilpotent orbit $\calO\subset\sN$, its stalk at $e\in\calO$ gives an irreducible representation $\rho$ of the centralizer $C_{G}(e)$ which factors through $A_{e}$. Note that the notion of $G$-equivariance changes when $G$ varies in a fixed isogeny class.  It is possible to extend the above injection \eqref{Spr corr} into a bijection by supplementing $\Irr(W)$ with $\Irr(W')$ for a collection of smaller Weyl groups. This is called the {\em generalized Springer correspondence} discovered by Lusztig \cite{LuIC}.
\end{remark}

Springer's original proof of Theorem \ref{th:Spr corr} uses trigonometric sums over $\frg(\FF_{q})$ and, when $k$ has characteristic zero, his proof uses reduction to finite fields. The following theorem due to Borho and MacPherson \cite{BM} can be used to give a direct proof of the Springer correspondence for all base fields $k$ of large characteristics or characteristic zero. To state it, we need to use the language of constructible (complexes of) $\Ql$-sheaves and perverse sheaves, for which we refer to the standard reference \cite{BBD} and de Cataldo's lectures \cite{dC} in this volume. 

\begin{theorem}\label{th:BM} The complex $\calS:=\bR\pi_{*}\Ql[\dim\sN]$ is a perverse sheaf on $\sN$ whose endomorphism ring is canonically isomorphic to the group algebra $\Ql[W]$. In particular, $W$ acts on the stalks of $\bR\pi_{*}\Ql$, i.e., $W$ acts on $\cohog{*}{\sB_{e}}$ for all $e\in\sN$.
\end{theorem}

We sketch three constructions of the $W$-action on $\bR\pi_{*}\Ql[\dim\sN]$.

\subsubsection{Construction via middle extension}\label{sss:mid} This construction (or rather the version where $\frg$ is replaced by $G$) is due to Lusztig \cite[\S3]{LuGreen}. The dimension formula for Springer fibers \eqref{Spr dim} imply that 
\begin{itemize}
\item The map $\pi:\tN\to\sN$ is semismall. \footnote{A proper surjective map $f:X\to Y$ of irreducible varieties is called {\em semismall} (resp. {\em small}) if for any $d\geq 1$, $\{y\in Y|\dim f^{-1}(y)\geq d\}$ has codimension at least $2d$ (resp. $2d+1$) in $Y$.}
\end{itemize}
There is an extension of the dimension formula \eqref{Spr dim} for the dimension of $\tB_{X}$ valid for all elements $X\in \frg$. Using this formula one can show that
\begin{itemize}
\item The map $\pi_{\frg}:\tg\to\frg$ is small.
\end{itemize}
As a well-known fact in the theory of perverse sheaves, the smallness of $\pi_{\frg}$ (together with the fact that $\tg$ is smooth) implies that $\calS_{\frg}:=\bR\pi_{\frg,*}\Ql[\dim\frg]$ is a perverse sheaf which is  the middle extension of its restriction to any open dense subset of $\frg$. Over the regular semisimple locus $\frg^{\rs}\subset \frg$, $\pi_{\frg}$ is a $W$-torsor, therefore $\calS_{\frg}|_{\frg^{rs}}$ is a local system shifted in degree $-\dim\frg$ that admits an action of $W$. By the functoriality of middle extension, $\calS_{\frg}$ admits an action of $W$. Taking stalks of $\calS_{\frg}$, we get an action of $W$ on $\cohog{*}{\tB_{X}}$ for all $X\in\frg$. 

In particular, for a nilpotent element $e$, we get an action of $W$ on $\cohog{*}{\sB_{e}}=\cohog{*}{\tB_{e}}$, because $\tB_{e}$ and $\sB_{e}$ have the same reduced structure. This is the action defined by Springer in his original paper \cite{Springer}, {\em which differs from our action by tensoring with the sign character of $W$}.

\subsubsection{Construction via Fourier transform}\label{sss:Four} By the semismallness of $\pi$,  the complex $\calS=\bR\pi_{*}\Ql[\dim\sN]$ is also a perverse sheaf. However, it is not the middle extension from an open subset of $\calN$. There is a notion of Fourier transform for $\Gm$-equivariant sheaves on affine spaces \cite{Laumon}. One can show that $\calS$ is isomorphic to the Fourier transform of $\calS_{\frg}$ and vice versa. The $W$-action on $\calS_{\frg}$ then induces an action of $W$ on $\calS$ by the functoriality of Fourier transform. Taking the stalk of $\calS$ at $e$ we get an action of $W$ on $\cohog{*}{\sB_{e}}$. This action is normalized according to our convention in \S\ref{sss:two actions}.

\subsubsection{Construction via correspondences}  Consider the Steinberg variety $\St_{\frg}=\tg\times_{\frg}\tg$ which classifies triples $(X,B_{1},B_{2})$, where $B_{1},B_{2}$ are Borel subgroups of $G$ and $X\in\Lie B_{1}\cap\Lie B_{2}$. The irreducible components of $\St_{\frg}$ are indexed by elements in the Weyl group: for $w\in W$, letting $\St_{w}$ be the closure of the graph of the $w$-action on $\tg^{\rs}$, then $\St_{w}$ is  an irreducible component of $\St$ and these exhaust all irreducible components of $\St$. The formalism of cohomological correspondences allows us to get an endomorphism of the complex $\calS_{\frg}=\bR\pi_{\frg,*}\Ql[\dim\frg]$ from each $\St_{w}$. It is nontrivial to show that these endomorphisms together form an action of $W$ on $\calS_{\frg}$. The key ingredient in the argument is still the smallness of the map $\pi_{\frg}$. After the $W$-action on $\calS_{\frg}$ is defined, one then define the Springer action on $\cohog{*}{\sB_{e}}$ by either twisting the action of $W$ on the stalk $\calS_{\frg, e}$ by the sign representation as in \S\ref{sss:mid}, or by using Fourier transform as in \S\ref{sss:Four}. We refer to \cite[Remark 3.3.4]{GS} for some discussion of this construction. See also \cite[\S3.4]{CG} for a similar but different construction using limits of $\St_{w}$ in the nilpotent Steinberg variety $\tN\times_{\sN}\tN$. 

Note that the above three constructions all allow one to show that $\End(\calS)\cong\Ql[W]$, hence giving a proof of Theorem \ref{th:BM}.

\subsubsection{Construction via monodromy} We sketch a construction of Slodowy \cite[\S4]{Slo4} which works for $k=\CC$. This construction was conjectured to give the same action of $W$ on $\cohog{*}{\sB_{e}}$ as the one in Theorem \ref{th:Spr corr}.  A similar construction by Rossmann appeared in \cite[\S2]{Rossmann}, in which the author identified his action with that constructed by Kazhdan and Lusztig in \cite{KL80}, and the latter was known to be the same as Springer's action. Thus all these constructions give the same $W$-action as in Theorem \ref{th:Spr corr}.

Let $e\in\calN$ and let $S_{e}\subset \frg$ be a transversal slice to the orbit of $e$. Upon choosing an $\sl_{2}$-triple $(e,h,f)$ containing $e$, there is a canonical choice of such a transversal slice $S_{e}=e+\frg_{f}$, where $\frg_{f}$ is the centralizer of $f$ in $\frg$. Now consider the following diagram where the squares are Cartesian except for the rightmost one
\begin{equation}\label{Slo}
\xymatrix{\sB_{e}\ar[d]\ar@{^{(}->}[r] & \tilS^{\nil}_{e}\ar[d]\ar@{^{(}->}[r] & \tilS_{e}\ar[d]\ar@{^{(}->}[r] & \tg\ar[d]^{\pi_{\frg}}\ar[r]^{\wt\chi} & \frt\ar[d]\\
\{e\}\ar@{^{(}->}[r] & S^{\nil}_{e}\ar@{^{(}->}[r] & S_{e}\ar@{^{(}->}[r] & \frg\ar[r]^{\chi} & \frc}
\end{equation}
Here the rightmost square is \eqref{tgg}, $S^{\nil}_{e}=S_{e}\cap\sN$ and $\tilS_{e}$ and $\tilS^{\nil}_{e}$ are the preimages of $S_{e}$ and $S^{\nil}_{e}$ under $\pi_{\frg}$. Let $V_{e}\subset S_{e}$ be a small ball around $e$ and let  $V_{0}\subset \frc$ be an even smaller ball around $0\in\frc$. Let $U_{e}=V_{e}\cap \chi^{-1}(V_{0})\subset S_{e}$. Then the diagram \eqref{Slo} restricts to a diagram
\begin{equation}\label{Slo U}
\xymatrix{\sB_{e}\ar[d]\ar@{^{(}->}[r] & \tilU^{\nil}_{e}\ar[d]\ar@{^{(}->}[r] & \tilU_{e}\ar[d]\ar[r]^{\wt\chi_{e}} & \tilV_{0}\ar[d]\\
\{e\}\ar@{^{(}->}[r] & U^{\nil}_{e}\ar@{^{(}->}[r] & U_{e}\ar[r] & V_{0}}
\end{equation} 
Here $U^{\nil}_{e}=U_{e}\cap\sN$, and $\tilV_{0}, \tilU_{e}$ and $\tilU^{\nil}_{e}$ are the preimages of $V_{0}, U_{e}$ and $U^{\nil}_{e}$ under the vertical maps in \eqref{Slo}. The key topological facts here are
\begin{itemize}
\item The inclusion $\sB_{e}\incl \tilU^{\nil}_{e}$ admits a deformation retract, hence it is a homotopy equivalence;
\item The map $\wt\chi_{e}: \tilU_{e}\to \tilV_{0}$ is a trivializable fiber bundle (in the sense of differential topology).
\end{itemize}
Now a general fiber of $\wt\chi_{e}$ admits a homotopy action of $W$ by the second point above because the rightmost  square in \eqref{Slo U} is Cartesian over $V_{0}\cap\frc^{\rs}$ and the map $\tilV_{0}\to V_{0}$ is a $W$-torsor over $V_{0}\cap\frc^{\rs}$. By the first point above, $\sB_{e}$ has the same homotopy type with $\tilU^{\nil}_{e}=\wt\chi^{-1}_{e}(0)$, hence $\sB_{e}$ also has the same homotopy type as a general fiber of $\wt\chi_{e}$ because $\wt\chi_{e}$ is a fiber bundle. Combining these facts, we get an action of $W$ on the {\em homotopy type of $\sB_{e}$}, which is a stronger structure than an action of $W$ on the cohomology of $\sB_{e}$. A consequence of this construction is that the $W$-action on $\cohog{*}{\sB_{e}}$ in Theorem  \ref{th:Spr corr} preserves the ring structure.

\subsubsection{Proof of Theorem \ref{th:Spr corr} assuming Theorem \ref{th:BM}}
We decompose the perverse sheaf $\calS$ into isotypical components under the $W$-action
\begin{equation*}
\calS=\bigoplus_{\chi\in\Irr(W)}V_{\chi}\otimes\calS_{\chi}
\end{equation*}
where $V_{\chi}$ is the space on which $W$ acts via the irreducible representation $\chi$, and $\calS_{\chi}=\Hom_{W}(V_{\chi},\calS)$ is a perverse sheaf on $\sN$. Since $\End(\calS)\cong\Ql[W]$, we conclude that
each $\calS_{\chi}$ is nonzero and that
\begin{equation}\label{HomSS}
\Hom(\calS_{\chi},\calS_{\chi'})=\begin{cases}\Ql & \chi=\chi' \\ 0 & \textup{otherwise}\end{cases}
\end{equation}
The decomposition theorem \cite[Th 6.2.5]{BBD} implies that each $\calS_{\chi}$ is a semisimple perverse sheaf. Therefore \eqref{HomSS} implies that $\calS_{\chi}$ is simple. Hence $\calS_{\chi}$ is of the form $\IC(\ov\calO, \calL)$ where $\calO\subset\sN$ is a nilpotent orbit and $\calL$ is an irreducible  $G$-equivariant local system on $\calO$. Moreover, since $\Hom(\calS_{\chi},\calS_{\chi'})=0$ for $\chi\neq\chi'$, the simple perverse sheaves $\{\calS_{\chi}\}_{\chi\in\Irr(W)}$ are non-isomorphic to each other. This proves part (3) of Theorem \ref{th:Spr corr} by interpreting the right side of \eqref{Spr corr} as the set of isomorphism classes of irreducible $G$-equivariant local systems on nilpotent orbits. If $\calS_{\chi}=\IC(\ov\calO,\calL)$ and $e\in\calO$, the semismallness of $\pi$ allows us to identify the stalk $\calL_{e}$ with an $A_{e}$-isotypic subspace of $\cohog{2d_{e}}{\sB_{e}}$. This proves part (2) of Theorem \ref{th:Spr corr}. \qed

We give some further examples of the Springer correspondence.

\subsubsection{Type $A$}\label{sss:A Spr corr} When $G=\SL_{n}$, all $A_{e}$ are trivial. The Springer correspondence sets a bijection between irreducible representations of $W=S_{n}$ and nilpotent orbits of $\frg=\sl_{n}$, both parametrized by partitions of $n$. In \S\ref{sss:Spr type A} we have seen that if $e$ has Jordan type $\l$, the top dimensional cohomology $\cohog{2d_{e}}{\sB_{e}}$ has a basis indexed by the standard Young tableaux of $\l$, the latter also indexing a basis of the irreducible representation of $S_{n}$ corresponding to the partition $\l$. For example, for $G=\SL_{3}$, the Springer correspondences reads
\begin{itemize}
\item trivial representation $\bij$ regular orbit, partition $3=3$;  
\item two-dimensional representation $\bij$ subregular orbit, partition $3=2+1$; 
\item sign representation $\bij$ \{0\}, partition $3=1+1+1$. 
\end{itemize}

\subsubsection{The subregular orbit and the reflection representation}\label{sss:subreg Spr corr} Consider the case $e$ is a subregular nilpotent element. In this case, the component group $A_{e}$ can be identified with the automorphism group of the Dynkin diagram of $G'$ introduced in \S\ref{sss:subreg} (see \cite[\S7.5, Proposition]{Slodowy}). After identifying the irreducible components of $\sB_{e}$ with the vertices of the Dynkin diagram of $G'$, the action of $A_{e}$ on $\cohog{2}{\sB_{e}}$ is by permuting the basis given by irreducible components in the same way as its action on the Dynkin diagram of $G'$. For example, when $G=G_{2}$,  we may write $\sB_{e}=C_{1}\cup C_{2}\cup C_{3}\cup C_{4}$ with $C_{1}$, $C_{2}$, $C_{3}$ each intersecting $C_{4}$ in a point and otherwise disjoint. The group $A_{e}$ is isomorphic to $S_{3}$, and its action on $\cohog{2}{\sB_{e}}$ fixes the fundamental class of $C_{4}$ and permutes the fundamental classes of $C_{1}, C_{2}$ and $C_{3}$.

Note that $\cohog{2}{\sB_{e}}^{A_{e}}$ always has dimension $r$, the rank of $G$. In fact, as a $W$-module, $\cohog{2}{\sB_{e}}^{A_{e}}$ is isomorphic to the reflection representation of $W$ on $\frt^{*}$. In other words, under the Springer correspondence, the pair $(e=\textup{subregular}, \rho=1)$ corresponds to the reflection representation of $W$.

\subsection{Comments and generalizations} 

\subsubsection{Extended symmetry} The $W$-action on $\cohog{*}{\sB_{e}}$ can be extended to an action of a larger algebra in various ways, if we use more sophisticated cohomology theories. On the equivariant cohomology $\upH^{*}_{\tilG_{e}}(\sB_{e})$, there is an action of the {\em graded affine Hecke algebra} (see Lusztig \cite{LuGrHk}). On the $\tilG_{e}$-equivariant $K$-group of $\sB_{e}$, there is an action of the {\em affine Hecke algebra} (see Kazhdan-Lusztig \cite{KLAffHk} and Chriss-Ginzburg \cite{CG}). 

\subsubsection{The group version} There are obvious analogs of the Springer resolution and the Grothendieck alteration when $\sN$ and $\frg$ are replaced with the unipotent variety $\sU\subset G$ and $G$ itself. When $\chark$ is large, the exponential map identifies $\sN$ with $\sU$ in a $G$-equivariant manner, hence the theories of Springer fibers for nilpotent elements and unipotent elements are identical. The group version of the perverse sheaf $\calS_{\frg}$ and its irreducible direct summands are precursors of {\em character sheaves}, a theory developed by Lusztig (\cite{LuChShI}, \cite{LuChShII}, \cite{LuChShIV} and \cite{LuChShV}) to study characters of the finite groups  $G(\FF_{q})$.

\subsubsection{Partial Springer resolutions}\label{sss:partial Spr} We may define analogs of $\sB_{e}$ in partial flag varieties.  Let $\sP$ be a partial flag variety of $G$ classifying parabolic subgroups $P$ of $G$ of a given type. There are two analogs of the map $\pi:\tN\to\sN$ one may consider. 

First,  instead of considering $\tN=T^{*}\sB$, we may consider $T^{*}\sP$, which classifies pairs $(e, P)\in\sN\times\sP$ such that $e\in \Lie\frn_{P}$, where $\frn_{P}$ is the nilpotent radical of $\Lie P$. Let $\tau_{\sP}: T^{*}\sP\to\sN$ be the first projection. In general this map is not surjective, its image is the closure of a nilpotent orbit $\calO_{\calP}$. The orbit $\calO_{P}$ is characterized by the property that its intersection with $\frn_{P}$ is dense in $\frn_{P}$, for any $P\in\sP$. This is called the {\em Richardson class} associated to parabolic subgroups of type $\sP$. When $G$ is of type $A$, each nilpotent class $\calO$ is the Richardson class associated to parabolic subgroups of some type $\sP$ (not unique in general). The map $T^{*}\sP\to \ov\calO$ is a resolution of singularities. For general $G$, not every nilpotent orbit is Richardson.

Second, we may consider the subscheme $\tN_{\sP}\subset\sN\times \sP$ classifying pairs $(e,P)$ such that $e\in\sN_{P}$, where $\sN_{P}\subset\Lie P$ is the nilpotent cone of $P$.  The projection $\pi_{\sP}: \tN_{\sP}\to \sN$ is now surjective, and is a partial resolution of singularities. The Springer resolution $\pi$ can be factored as
\begin{equation*}
\pi: \tN=\tN_{\sB}\xrightarrow{\nu_{\sP}}\tN_{\sP}\xrightarrow{\pi_{\sP}} \sN. 
\end{equation*}

We have an embedding $T^{*}\sP\incl \tN_{\sP}$.  We may consider the fibers of either $\tau_{\sP}$ or $\pi_{\sP}$ as parabolic analogs of Springer fibers.  We call them {\em partial Springer fibers}. The Springer action of $W$ on the cohomology of $\sB_{e}$ has an analog for partial Springer fibers. For more information, we refer the readers to \cite{BM}.

\subsubsection{Hessenberg varieties} The Grothendieck alteration $\pi_{\frg}:\tg\to \frg$ admits a generalization where $\frg$ is replaced with another linear representation of $G$.

Fix a Borel subgroup $B$ of $G$. Let $(V,\rho)$ be a representation of $G$ and $V^{+}\subset V$ be a subspace which is stable under $B$. Now we use the pair $(V,V^{+})$ instead of the pair $(\frg,\frb)$, we get a generalization of the Grothendieck alteration. More precisely, let $\tilV\subset V\times \calB$ be the subscheme consisting of pairs $(v, gB)\in V\times\sB$ such that $v\in \rho(g)V^{+}$. Let $\pi_{V}:\tilV\to V$ be the first projection. The fibers of $\pi_{V}$ are called {\em Hessenberg varieties}.

Hessenberg varieties appear naturally in the study of certain affine Springer fibers, as we will see in \S\ref{ss:aff examples}. For more information on Hessenberg varieties, see \cite{GKM} and \cite{OY}.

\subsection{Exercises}

\subsubsection{} For $G=\SL_{n}$, determine the sizes of the Jordan blocks of a regular and subregular nilpotent element of $\frg$.
 
\subsubsection{}\label{exI:nonred Spr} For $G=\SL_{2}$ and $\SL_{3}$, calculate the coordinate ring of the {\em non-reduced} Springer fiber $\tB_{e}$ for a regular nilpotent element $e$. Show also that the Springer fiber $\sB_{e}$ is indeed a reduced point. 

{\em Hint}: if you write $e$ as an upper triangular matrix, then $\tB_{e}$ lies in the big Bruhat cell of the flag variety $\sB$, from which you get coordinates for your calculation.

\subsubsection{}\label{exI:211Sp4} Verify the statement in \S\ref{sss:Sp4}: consider $G=\Sp(V)$, $V=\jiao{v_{1},v_{2},v_{3},v_{4}}$ with the symplectic form given by $\jiao{v_{i},v_{5-i}}=1$ if $i=1,2$ and $\jiao{v_{i},v_{j}}=0$ for $i+j\neq 5$. Let $e: v_{4}\mapsto v_{1}\to 0, v_{2}\mapsto 0,v_{3}\mapsto 0$. Then any flag $0\subset V_{1}\subset V_{2}\subset V_{1}^{\bot}\subset V$ in $\sB_{e}$ must satisfy
\begin{equation*}
\jiao{v_{1}}\subset V_{2}\subset\jiao{v_{1},v_{2},v_{3}}.
\end{equation*}
Moreover, this is the only condition for a flag $0\subset V_{1}\subset V_{2}\subset V_{1}^{\bot}\subset V$ to lie in $\sB_{e}$.

\subsubsection{}\label{exI:inv} Let $e\in\sN$.  Let $B\subset G$ be a Borel subgroup.
\begin{enumerate}
\item  Let $\alpha$ be a simple root. Let $P_{\alpha}\supset B$ be a parabolic subgroup whose Levi factor has semisimple rank one with roots $\pm\alpha$. Let $\sP_{\alpha}\cong G/P_{\alpha}$ be the partial flag variety of $G$ classifying parabolics conjugate to $P_{\alpha}$. Restricting the projection $\sB\to \sP_{\alpha}$  to $\sB_{e}$, we get a map
\begin{equation*}
\pi_{\alpha}: \sB_{e}\to \pi_{\alpha}(\sB_{e}).
\end{equation*}
Show that the pullback $\pi^{*}_{\alpha}$ on cohomology induces an isomorphism
\begin{equation}\label{sa inv}
\cohog{*}{\pi_{\alpha}(\sB_{e})}\cong \cohog{*}{\sB_{e}}^{s_{\alpha}}
\end{equation}
where $s_{\alpha}\in W$ is the simple reflection associated with $\alpha$, which acts on $\cohog{*}{\sB_{e}}$ via Springer's action.
\item Can you generalize \eqref{sa inv} to other partial flag varieties?
\end{enumerate}

\subsubsection{}\label{exI:SL3}  Let $G=\SL_{3}$ and let $e\in\sN$ be a subregular element. Calculate the action of $S_{3}$ on the two dimensional $\cohog{2}{\sB_{e}}$ in terms of the basis given by the fundamental classes of the two irreducible components (see \S\ref{sss:SL3 subreg}).

{\em Hint}: for this and the next problem, you may find Exercise \ref{exI:inv} useful.

\subsubsection{}\label{exI:Sp4} Describe the Springer fibers for $\Sp_{4}$. Calculate the Springer correspondence for $G=\Sp_{4}$ explicitly.

\subsubsection{} Using the dimension formula for $\sB_{e}$, verify that the Springer resolution $\pi$ is semismall.

\subsubsection{}\label{exI:dim BX} Let $X\in\frg$ and let $X=X_{s}+X_{n}$ be the Jordan decomposition of $X$. Let $H=C_{G}(X_{s})\subset G$. This is a Levi subgroup of $G$. Let $\sB^{H}$ be the flag variety of $H$ and let $\sB^{H}_{X_{n}}$ be the Springer fiber of $X_{n}$ viewed as a nilpotent element in $\Lie H$. Show that $\dim\tB_{X}=\dim \sB^{H}_{X_{n}}$.

\subsubsection{} Use Exercise \ref{exI:dim BX} and the dimension formula \eqref{Spr dim} to derive a formula for the dimension of $\tB_{X}$ for all elements $X\in\frg$. Use your formula to prove that the Grothendieck alteration $\pi_{\frg}$ is small.

\subsubsection{} Show that $\sN$ is rationally smooth; i.e., its intersection cohomology complex is isomorphic to the shifted constant sheaf $\Ql[\dim\sN]$.

{\em Hint}: the largest direct summand in $\bR\pi_{*}\Ql$ is the shifted IC sheaf of $\sN$, and it is also the restriction of a direct summand of $\bR\pi_{\frg,*}\Ql$.

\subsubsection{}\label{exI:conn} Show that the Springer fibers $\sB_{e}$ are connected. 

{\em Hint}: the $\upH^{0}$ of the Springer fibers are packed in some sheaf.

\subsubsection{} Denote the simple roots of $G$ by $\{\alpha_{1},\cdots,\alpha_{r}\}$. A parabolic subgroup $P\subset G$ is called of type $i$ if the roots of its Levi quotient $L_{P}$ are $\pm\alpha_{i}$.
\begin{enumerate}
\item Let $1\leq i\leq r$ and let $P$ be a parabolic subgroup of type $i$. Let $\frn_{P}$ be the nilpotent radical of $\Lie P$. Show that $\frn_{P}\cap \calO_{\subreg}$ is dense in $\frn_{P}$.
\item Let $e\in\calO_{\subreg}$. Show that for each $i$, there are finitely many parabolics $P$ of type $i$ such that $e\in\frn_{P}$. For each such $P$, the subvariety $C_{P}:=\{B\in\sB_{e}|B\subset P\}$ of $\sB_{e}$ is isomorphic to $\PP^{1}$, and is called a curve of type $i$.
\item For parabolics $P\neq Q$ of type $i$ and $j$, show that $C_{P}\cap C_{Q}$ is either empty or a point.
\item Let $P$ be a parabolic subgroup of type $i$ such that $e\in\frn_{P}$. For any $1\leq j\leq r, i\neq j$, $C_{P}$ intersects exactly $-\jiao{\alpha^{\vee}_{i}, \alpha_{j}}$ curves of type $j$.
\item Show that $\sB_{e}$ has the configuration described in \S\ref{sss:subreg}.
\item Use Exercise \ref{exI:inv} to calculate the Springer action of $W$ on $\cohog{2}{\sB_{e}}$.
\end{enumerate}


\section{Lecture II: Affine Springer fibers}\label{s:ASF}

Affine Springer fibers are analogs of Springer fibers for loop groups. They were introduced by Kazhdan and Lusztig \cite{KL}. Roughly speaking, in the case of classical groups, instead of classifying flags in a $k$-vector space fixed by a $k$-linear transformation, affine Springer fibers classify (chains of) lattices in an $F$-vector space fixed by an $F$-linear transformation, where $F=k\lr{t}$. The cohomology groups of affine Springer fibers carry actions of the affine Weyl group.

The setup in this section is the same as in \S\ref{ss:notation}.

\subsection{Loop group, parahoric subgroups and the affine flag variety}
Let $F=k\lr{t}$ be the field of formal Laurent series in one variable $t$.  Then $F$ has a discrete valuation $\val_{F}:F^{\times}\to\ZZ$ such that $\val_{F}(t)=1$ and its valuation ring is $\calO_{F}=k\tl{t}$.

\subsubsection{The loop group} The {\em positive loop group} $L^{+}G$ is a group-valued functor on $k$-algebras. For any $k$-algebra $R$, we define $LG(R):=G(R\tl{t})$. It turns out that $LG$ is representable by a scheme over $k$ which is not of finite type. 

For example, when $G=\GL_{n}$, an element in $LG(R)=\GL_{n}(R\tl{t})$ is given by $n^{2}$ formal Laurent series $a_{ij}=\sum_{s\geq 0}a^{(s)}_{ij}t^{s}$ ($1\leq i,j\leq n$), with $a^{(s)}_{ij}\in R$, subject to one condition that $\det((a_{ij}))$ (which is a polynomial in the $a^{(s)}_{ij}$ of degree $n$) is invertible in $R$. Therefore in this case $LG$ is an open subscheme in the infinite dimensional affine space with coordinates $a^{(s)}_{ij}$, $1\leq i,j\leq n$ and $s\geq 0$. 

Similarly we may define the {\em loop group} $LG$ to be the functor $LG(R)=G(R\lr{t})$ on $k$-algebras $R$. The functor $LG$ is no longer representable by a scheme, but rather by an {\em ind-scheme}. An ind-scheme is an inductive limit $\varinjlim_{m} X_{m}$ in the category of schemes, i.e., $\{X_{m}\}$ form an inductive system of schemes over $k$, and $\varinjlim_{m} X_{m}$ is the functor $R\mapsto \varinjlim_{m} X_{m}(R)$.  When $G=\GL_{n}$, we may  define $X_{m}$ to be the subfunctor of $LG$ such that $X_{m}(R)$ consists of $n$-by-$n$ invertible matrices with entries in $t^{-m}R\tl{t}\subset R\lr{t}$. Then the same argument as in the  $L^{+}G$ case shows that $X_{m}$ is representable by a scheme over $k$. For $m<m'$, we have a natural closed embedding $X_{m}\incl X_{m'}$, and $LG$ in this case is the inductive limit $\varinjlim_{m}X_{m}$.  For general $G$, see \cite[\S1]{BL} and \cite[\S2]{Faltings}.

\subsubsection{The affine Grassmannian}\label{sss:Gr} The affine Grassmannian $\Gr_{G}$ of $G$ is defined as the sheafification of the presheaf $R\mapsto LG(R)/L^{+}G(R)$ in the category of $k$-algebras under the fpqc topology.
In particular, we have $\Gr_{G}(k)=G(F)/G(\calO_{F})$. 

When $G=\GL_{n}$, the affine Grassmannian $\Gr_{G}$ can be identified with the moduli space of projective $R\tl{t}$-submodules $\L\subset R\lr{t}^{n}$ such that
\begin{equation}\label{lattice range}
(t^{m}R\tl{t})^{n}\subset \L\subset (t^{-m}R\tl{t})^{n}
\end{equation}
for some $m\geq 0$. Such an $R\tl{t}$-module $\L$ is called a {\em lattice} in $R\lr{t}^{n}$.  For fixed $m$, let $X_{m}$ be the subfunctor of $\Gr_{G}$ classifying those $\L$ such that \eqref{lattice range} holds, then $X_{m}$ is representable by a projective scheme over $k$. The natural closed embeddings $X_{m}\incl X_{m+1}$ make $\{X_{m}\}$ into an inductive system of projective schemes, and $\Gr_{G}$ is representable by the ind-scheme $\varinjlim_{m}X_{m}$.

Let us elaborate on the bijection between $\Gr_{G}(k)$ and the lattices in $F^{n}=k\lr{t}^{n}$. Let $\calO_{F}^{n}\subset F^{n}$ be the standard lattice. Let $\frL_{n}$ be the set of lattices in $F^{n}$ (in the case $R=k$ a lattice is simply an $\calO_{F}$-submodules of $F^{n}$ of rank $n$). The group $LG(k)=G(F)$ acts on $\frL_{n}$ by $LG\ni g:\L\mapsto g\L$. This action is transitive and the stabilizer of the standard lattice $\calO^{n}_{F}$ is $L^{+}G(k)=G(\calO_{F})$. Therefore this action induces a $G(F)$-equivariant bijection
\begin{equation}\label{Gr lattices}
\Gr_{G}(k)=G(F)/G(\calO_{F})\isom \frL_{n}.
\end{equation}

For general $G$, $\Gr_{G}$ is always representable by an ind-scheme of the form $\varinjlim_{m}X_{m}$ where $X_{m}$ are projective schemes over $k$, and the transition maps $X_{m}\incl X_{m+1}$ are closed embeddings.  We have a canonical exhaustion of $\Gr_{G}$ by projective schemes given by the affine Schubert stratification, which we now recall. The action of $L^{+}G$ on $\Gr_{G}$ by left translation has orbits indexed by dominant cocharacters $\l\in\xcoch(T)^{+}$. We denote by $\Gr_{G,\l}$ the $L^{+}G$-orbit through $\l(t)$. Let $\Gr_{G,\leq \l}$ be the closure of $\Gr_{G,\l}$. Then $\Gr_{G,\leq\l}$ is a projective scheme and $\Gr_{G}$ is the union of $\Gr_{G,\leq\l}$. For more details on the affine Grassmannian, we refer to \cite[\S2]{BL}, \cite[\S2]{Faltings} and Zhu's lectures \cite{Zhu}.

\subsubsection{Parahoric subgroups} The subgroup $L^{+}G$ of $LG$ is an example of a class of subgroups called {\em parahoric subgroups}.  Fix a Borel subgroup $B\subset G$ and let $\bI\subset L^{+}G$ be the preimage of $B$ under the map $L^{+}G\to G$ given by reduction modulo $t$. Then $\bI$ is an example of an {\em Iwahori subgroup} of $LG$. General Iwahori subgroups are conjugates of $\bI$ in $LG$. Like $L^{+}G$, Iwahori subgroups are group subschemes of $LG$ of infinite type. Parahoric subgroups are connected group subschemes of $LG$ containing an Iwahori subgroup with finite codimension. A precise definition of parahoric subgroups involves a fair amount of Bruhat-Tits theory, which we refer the readers to the original papers of Bruhat and Tits \cite{BT1}, and the survey paper \cite{Tits}. 

Just as the conjugacy classes of parabolic subgroups of $G$ are in bijection with subsets of the Dynkin diagram of $G$, the $LG$-conjugacy classes of parahoric subgroups of $LG$ are in bijection with {\em proper subsets} of the vertices of the {\em extended Dynkin diagram} $\wt{\Dyn}(G)$ of $G$, which has one more vertex than the Dynkin diagram of $G$. See Kac's book \cite[\S4.8]{Kac}, Bourbaki \cite[Ch VI]{Bourbaki} for extended Dynkin diagrams and the expository paper of Gross \cite{Gross} for connection with parahoric subgroups. 

Each $\bP$ admits a canonical exact sequence of group schemes
\begin{equation*}
1\to \bP^{+}\to \bP\to L_{\bP}\to1
\end{equation*}
where $\bP^{+}$ is the {\em pro-unipotent radical} of $\bP$ and $L_{\bP}$ is a reductive group over $k$, called the {\em Levi quotient} of $\bP$. If $\bP$ corresponds to a subset $J$ of the vertices of $\wt{\Dyn}(G)$, then the Dynkin diagram of the Levi quotient $L_{\bP}$ is the sub-diagram of $\wt{\Dyn}(G)$ spanned by $J$.

\subsubsection{Affine flag varieties} For each parahoric subgroup $\bP\subset LG$ we may define the corresponding {\em affine partial flag variety} $\Fl_{\bP}$ as the fpqc sheafification of the functor $R\mapsto LG(R)/\bP(R)$ on the category of $k$-algebras. This functor is also representable by an ind-scheme $\varinjlim_{m}X_{m}$ where each $X_{m}$ is a projective scheme over $k$ and the transition maps are closed embeddings.  The affine Grassmannian $\Gr_{G}$ is a special case of $\Fl_{\bP}$ for $\bP=L^{+}G$.

Consider the special case $\bP=\bI$ is an Iwahori subgroup of $LG$. When $G$ is simply-connected, $\bI$ is its own normalizer, and we may identify $\Fl_{\bI}$ as the moduli space of Iwahori subgroups of $LG$, hence giving an intrinsic definition of the affine flag variety. We usually denote $\Fl_{\bI}$ by $\Fl$ or $\Fl_{G}$ and call it {\em the affine flag variety of $G$}, with the caveat that $\Fl_{\bI}$ is canonically independent of the choice of $\bI$ only when $G$ is simply-connected.

Let $\bP\subset\bQ$ be two parahoric subgroups of $LG$. Then we have a natural projection $\Fl_{\bP}\to \Fl_{\bQ}$. The fibers of this projection are isomorphic to the partial flag variety of $L_{\bQ}$ corresponding to its parabolic subgroup given by the image of $\bP\to L_{\bQ}$. In particular, there is a natural projection $\Fl_{G}\to \Gr_{G}$ whose fibers are isomorphic to the flag variety $\sB$.

\subsubsection{The case of $\SL_{n}$}\label{sss:GLn parahoric} We have seen in \S\ref{sss:Gr} that the affine Grassmannian of $\GL_{n}$ has an interpretation as the moduli space of lattices. In fact, parahoric subgroups of $LG$ and the associated affine partial flag varieties can also be described using lattices. Here we consider the case $G=\SL_{n}$.  

Recall that the set of lattices in $F^{n}$ is denoted by $\frL_{n}$. For any two lattices $\L_{1},\L_{2}\in\frL_{n}$ we may define their relative length to be the integer
\begin{equation*}
[\L_{1}:\L_{2}]:=\dim_{k}(\L_{1}/\L_{1}\cap\L_{2})-\dim_{k}(\L_{2}/\L_{1}\cap\L_{2}).
\end{equation*}
Let $J\subset \ZZ/n\ZZ$ be a non-empty subset. Let $\tilJ$ be the preimage of $J$ under the projection  $\ZZ\to\ZZ/n\ZZ$. A {\em periodic $J$-chain of lattices} is a function $\L:\tilJ\to \frL_{n}$ sending each $i\in \tilJ$ to a lattice $\L_{i}\in\frL_{n}$ such that
\begin{itemize}
\item $[\L_{i}:\calO_{F}^{n}]=i$  for all $i\in \tilJ$;
\item $\L_{i}\subset \L_{j}$ for $i<j$ in $\tilJ$;
\item $\L_{i}=t\L_{i+n}$ for all $i\in \tilJ$.
\end{itemize}
Let $\frL_{J}$ be the set of periodic $J$-chains of lattices. For each $\{\L_{i}\}_{i\in \tilJ}\in\frL_{J}$, let $\bP_{\{\L_{i}\}_{i\in J}}\subset LG$ to be the simultaneous stabilizers of all $\L_{i}$'s. Then $\bP_{\{\L_{i}\}_{i\in J}}$ is a parahoric subgroup of $LG$. We call such a parahoric subgroup of type $J$. In fact all parahoric subgroups of $LG$ arise from a unique periodic $J$-chain of lattices, for a unique non-empty $J\subset \ZZ/n\ZZ$.  Therefore we get a bijection between $\sqcup_{J}\frL_{J}$ and the set of parahoric subgroups of $LG$.  In particular, $L^{+}G$ is the parahoric subgroup corresponding to the periodic $\{0\}$-chain of lattices given by $\L_{i}=t^{i/n}\calO^{n}_{F}$, where $i\in \tilJ=n\ZZ$. 

The extended Dynkin diagram of $G$ is a loop with $n$ nodes which we index cyclically by the set $\ZZ/n\ZZ$, such that $0$ corresponds to the extra node compared to the usual Dynkin diagram. Parahoric subgroups of type $J\neq\varnothing$ corresponds to the proper subset $\ZZ/n\ZZ-J$ of the nodes of the extended Dynkin diagram.
 
One can find the moduli space $\Fl_{J}$ of periodic $J$-chains of lattices such that  $\Fl_{J}(k)=\frL_{J}$. Fixing any parahoric subgroup $\bP$ of type $J$, $\Fl_{J}$ can be identified with the affine partial flag variety $\Fl_{\bP}$.  In particular, the affine flag variety $\Fl$ for $G=\SL_{n}$ can be identified with the moduli space of {\em periodic full chains of lattices}, i.e., a sequence of lattices $\cdots \L_{-1}\subset \L_{0}\subset \L_{1}\cdots$ in $F^{n}$ with $[\L_{i}:\calO^{n}_{F}]=i$ and $\L_{i}=t\L_{i+n}$ for all $i\in\ZZ$.

\subsubsection{The case of $\Sp_{2n}$}\label{sss:Sp parahoric} Now consider $G=\Sp_{2n}=\Sp(V)$, where $V=k^{2n}$ is equipped with a symplectic form. We extended the symplectic form on $V$ $F$-linearly to a symplectic form $\jiao{-,-}$ on $V\otimes_{k}F$. For a lattice $\L\in \frL_{2n}$,  define its symplectic dual to be the set $\L^{\vee}:=\{v\in V\otimes_{k}F|\jiao{v,\L}\subset \calO_{F}\}$. Then $\L^{\vee}$ is again a lattice in $V\otimes_{k}F$. The operation $\L\mapsto \L^{\vee}$ defines an involution on $\frL_{2n}$.  

Let $J\subset\ZZ/2n\ZZ$ be a non-empty subset stable under multiplication by $-1$. Let $\tilJ\subset \ZZ$ be the preimage of $J$ under the natural projection $\ZZ\to\ZZ/2n\ZZ$. A {\em periodic self-dual $J$-chain of lattices} in $V\otimes_{k}F$ is a periodic $J$-chain of lattices (i.e., an element in $\frL_{J}$ in the notation of \S\ref{sss:GLn parahoric}) satisfying the extra condition that
\begin{equation*}
\L^{\vee}_{i}=\L_{-i}, \textup{ for all } i\in \tilJ.
\end{equation*}
Denote the set of periodic self-dual $J$-chains of lattices in $V\otimes_{k}F$ by $\frL^{\Sp(V)}_{J}$. This is a set with an action of $G(F)=\Sp(V\otimes_{k}F)$. For any $\{\L_{i}\}_{i\in\tilJ}\in \frL^{\Sp(V)}_{J}$, the simultaneous stabilizer of the $\L_{i}$'s is a parahoric subgroup of $LG$, and every parahoric subgroup of $LG$ arises this way. For a parahoric subgroup $\bP$ of type $J$, the corresponding affine partial flag variety $\Fl_{\bP}$ can be identified with the moduli space of periodic self-dual $J$-chains of lattices so that $\Fl_{\bP}(k)\cong\frL^{\Sp(V)}_{J}$ as $G(F)$-sets. The readers are invited to work out the similar story for orthogonal groups, see Exercise \ref{exII:SO parahoric}.

\subsection{Affine Springer fibers} 

\subsubsection{Affine Springer fibers in the affine Grassmannian} For any $k$-algebra $R$, we denote $\frg\otimes_{k}R$ by $\frg(R)$. In particular, $\frg(F)=\frg\otimes_{k}F$ is the Lie algebra of the loop group $LG$. For $g\in LG$, let $\Ad(g)$ denote its adjoint action on $\frg(F)$. 

Let $\gamma\in\frg(F):=\frg\otimes_{k}F$ be a regular semisimple element \footnote{Here we are dealing with a Lie algebra $\frg$ over the non-algebraically-closed field $F$. An element  $\gamma\in\frg(F)$ is regular semisimple if it is regular semisimple as an element in $\frg(\overline{F})$, see the footnote in \S\ref{sss:Gr alt}. Equivalently, $\gamma$ is regular semisimple if its image in $\frc(F)$ lies in $\frc^{\rs}(F)$.}.  We consider the subfunctor of $\Gr_{G}$ whose value on a $k$-algebra $R$ is given by 
\begin{equation}\label{aff Spr}
\tX_{\gamma}(R)=\{[g]\in \Gr_{G}(R)|\Ad(g^{-1})\gamma\in\frg(R\tl{t})\}.
\end{equation} 
Then $\tX_{\gamma}$ is a closed sub-ind-scheme of $\Gr_{G}$. Let $\sX_{\gamma}=\tX^{\red}_{\gamma}$ be the underlying reduced ind-scheme of $\tX_{\gamma}$. We call $\sX_{\gamma}$ the {\em affine Springer fiber} of $\gamma$ in the affine Grassmannian $\Gr_{G}$.

\subsubsection{Alternative definition in terms of lattices}\label{sss:ASF Gr GLn} We consider the case $G=\GL_{n}$. Let $\gamma\in\frg(F)=\gl_{n}(F)$ be a regular semisimple matrix. As in \ref{sss:Gr} we identify $\Gr_{G}$ with the moduli space of lattices in $F^{n}$, or more precisely $\Gr_{G}(R)$ is the set of lattices in $R\lr{t}^{n}$. Then $\tX_{\gamma}(R)$ can be identified with those lattices $\L\subset R\lr{t}^{n}$ such that $\gamma \L\subset \L$, i.e., those stable under the endomorphism of $R\lr{t}^{n}$ given by $\gamma$. 

When $G=\SL_{n}$, $\Gr_{G}(R)$ classifies lattices $\L$ in $R\lr{t}^{n}$ such that $[\L:R\tl{t}^{n}]=0$. The affine Springer fiber $\tX_{\gamma}$ in this case is cut out by the same condition $\gamma \L\subset \L$.

When $G=\Sp_{2n}$, $\Gr_{G}(R)$ classifies lattices $\L$ in $R\lr{t}^{2n}$ such that $\L^{\vee}=\L$, see \S\ref{sss:Sp parahoric}. The affine Springer fiber $\tX_{\gamma}$ in this case is also cut out by the same condition $\gamma \L\subset \L$.

We give the simplest examples of affine Springer fibers. 

\subsubsection{}\label{sss:reduction rs} Let $\gamma\in\frt(\calO_{F})$  such that the reduction $\ov\gamma\in\frt$ is regular semisimple. For each cocharacter $\l:\Gm\to T$, the element $t^{\l}:=\l(t)\in T(F)$ gives a point $[t^{\l}]\in \Gr_{G}(k)$ which lies in $\tX_{\gamma}$ since $\Ad(t^{-\l})\gamma=\gamma\in \frg(\calO_{F})$. The reduced ind-scheme $\sX_{\gamma}$ is in fact the discrete set $\{[t^{\l}]\}$ which is in bijection with $\xcoch(T)$. More canonically, there is an action of the loop group $LT$ on $\tX_{\gamma}$ given by its left translation action on $\Gr_{G}$. This action factors through the quotient $\Gr_{T}=LT/L^{+}T$ and realizes $\tX_{\gamma}$ as a $\Gr_{T}$-torsor.

\subsubsection{}\label{P1 chain} Consider the case $G=\SL_{2}$ and $\gamma=\mat{t}{0}{0}{-t}$. Then $\sX_{\gamma}$ is an infinite chain of $\PP^{1}$'s. More precisely, for each $n\in\ZZ$, we consider the subscheme $C_{n}$ of $\Gr_{G}$ classifying lattices $\L\subset F^{2}$ such that $t^{n}\calO_{F}\oplus t^{-n+1}\calO_{F} \subset \L\subset t^{n-1}\calO_{F}\oplus t^{-n}\calO_{F}$. Then $C_{n}\cong\PP^{1}$. We have $\sX_{\gamma}=\cup_{n\in\ZZ}C_{n}$ is an infinite chain of $\PP^{1}$'s. The components  $C_{n}$ and $C_{n+1}$ intersect at one point $t^{-n}\calO_{F}\oplus t^{n}\calO_{F}$ and otherwise the components are disjoint.

Here is a way to calculate the $k$-points of $\sX_{\gamma}$. We use the {\em Iwasawa decomposition} for $G(F)$: 
\begin{equation*}
G(F)=\bigsqcup_{n\in\ZZ}N(F)\mat{t^{n}}{0}{0}{t^{-n}}G(\calO_{F}).
\end{equation*}
According to this decomposition, any point in $\Gr_{G}$ can be represented by
\begin{equation}\label{Iwasawa elt}
g=\mat{1}{x}{0}{1}\mat{t^{n}}{0}{0}{t^{-n}}
\end{equation}
for some  $x\in F$ and a unique $n\in\ZZ$, and $x$ has a well-defined image in $F/t^{2n}\calO_{F}$. Since
\begin{equation*}
\Ad(g^{-1})\gamma=g^{-1}\gamma g=\mat{t}{2t^{1-2n}x}{0}{-t}
\end{equation*}
the condition $\Ad(g^{-1})\gamma\in\frg(\calO_{F})$ is the same as requiring $x\in t^{2n-1}\calO_{F}$. Therefore $\sX_{\gamma}(k)=\sqcup_{n\in\ZZ}Y_{n}$ where $Y_{n}$ consists of elements of the form \eqref{Iwasawa elt}  with $x\in t^{2n-1}\calO_{F}/t^{2n}\calO_{F}$. Therefore each $Y_{n}$ can be identified with $k$. It is easy to check that $Y_{n}\subset C_{n}(k)$.

\subsubsection{}\label{SL2 ram dim 1} Consider the case $G=\SL_{2}$ and $\gamma=\mat{0}{t^{2}}{t}{0}$. Then $\sX_{\gamma}$ consists exactly of those lattices $\L\in\Gr_{G}$ such that $t\calO_{F}\oplus \calO_{F}\subset\L\subset\calO_{F}\oplus t^{-1}\calO_{F}$. Therefore $\sX_{\gamma}\cong\PP^{1}$. Details of these calculations are left to the reader, see Exercises \ref{exII:SL2}.

\subsubsection{Invariance under conjugation} If $\gamma,\gamma'\in\frg(F)$ are related by $\Ad(g)\gamma=\gamma'$ for some $g\in G(F)$, then the left multiplication by $g$ on $\Gr_{G}$ restricts to an isomorphism $\tX_{\gamma}\cong\tX_{\gamma'}$, hence also $\sX_{\gamma}\cong\sX_{\gamma'}$. Therefore the isomorphism type of $\sX_{\gamma}$ is invariant under $G(F)$-conjugation on $\gamma$.  Recall we have map $\chi:\frg\to \frc:=\frg\sslash G\cong\frt\sslash W$. For a regular semisimple point $a\in\frc^{\rs}(F)$, the fiber $\chi^{-1}(a)$ is a single $G(F)$-conjugacy class (here we are using that the residue field of $F$ is algebraically closed). Therefore, the isomorphism type of $\sX_{\gamma}$  depends only on $a=\chi(\gamma)\in\frc^{\rs}(F)$. 

Unlike Springer fibers, $\sX_{\gamma}$ can be empty for certain $\gamma$. The affine Springer fiber $\sX_{\gamma}$ is nonempty if and only if $a=\chi(\gamma)\in\frc(\calO_{F})$. In fact, if $gG(\calO_{F})\in\sX_{\gamma}(k)$, then $\Ad(g^{-1})\gamma\in\frg(\calO_{F})$ hence $\chi(\gamma)=\chi(\Ad(g^{-1})\gamma)\in\frc(\calO_{F})$. Conversely, we have a Kostant section $\ep: \frc\to \frg$ of $\chi$ which identifies $\frc$ with $e+\frg_{f}$, where $(e,h,f)$ is a regular $\sl_{2}$-triple in $\frg$. Therefore, for any $a\in\frc(\calO_{F})\cap\frc^{\rs}(F)$, $\ep(a)\in\frg(\calO_{F})$, and $\sX_{\ep(a)}$ contains the unit coset in $\Gr_{G}$ hence nonempty; since $\sX_{\gamma}$ is isomorphic to $\sX_{\ep(\chi(\gamma))}$, it is also nonempty.

For $a\in\frc(\calO_{F})\cap\frc^{\rs}(F)$, we also write $\sX_{a}$ for $\sX_{\ep(a)}$. The above discussion shows that all Springer fibers $\sX_{\gamma}$ are isomorphic to $\sX_{a}$ for $a=\chi(\gamma)$.

\subsubsection{Parahoric versions} For each parahoric subgroup $\bP\subset LG$, we may similarly define the closed sub-indscheme $\tX_{\bP,\gamma}\subset\Fl_{\bP}$ using the analog of the condition \eqref{aff Spr} with $\frg(R\tl{t})$ replaced by $(\Lie\bP)\widehat{\otimes}_{k}R$. The reduced ind-scheme $\sX_{\bP,\gamma}=\tX_{\bP,\gamma}^{\red}\subset\Fl_{\bP}$ is called the {\em affine Springer fiber of $\gamma$ of type $\bP$}. In particular, when $\bI$ is an Iwahori subgroup of $LG$, we denote $\sX_{\bI,\gamma}$ by $\sY_{\gamma}$.

For $\bP\subset\bQ$ two parahoric subgroups of $LG$, the natural projection $\Fl_{\bP}\to \Fl_{\bQ}$ induces a map $\sX_{\bP,\gamma}\to \sX_{\bQ,\gamma}$. In particular we have a map $\sY_{\gamma}\to \sX_{\gamma}$.

\subsection{Symmetry on affine Springer fibers}
For the Springer fiber $\sB_{e}$, the centralizer $C_{G}(e)$ acts on it. In this subsection we investigate a similar structure for affine Springer fibers. 

\subsubsection{Centralizer action}\label{sss:centralizer} Let $G_{\gamma}$ be the centralizer of $\gamma$ in $G_{F}$ (the algebraic group over $F$ obtained from $G$ by base change). Then $G_{\gamma}$ is an algebraic group over $F$. Since $\gamma$ is regular semisimple, $G_{\gamma}$ is a torus over $F$. One can define the loop group $LG_{\gamma}$ of $G_{\gamma}$ as the functor $R\mapsto G_{\gamma}(R\lr{t})$ on $k$-algebras.

We claim that $LG_{\gamma}$ acts on the ind-scheme $\tX_{\gamma}$. This can be seen on the level of $R$-points. Suppose $h\in LG_{\gamma}(R)=G_{\gamma}(R\lr{t})$ and $[g]\in \tX_{\gamma}(R)$. Then the coset $[hg]\in \Gr_{G}(R)$ still satisfies
\begin{equation*}
\Ad((hg)^{-1})\gamma=\Ad(g^{-1})\Ad(h^{-1})\gamma=\Ad(g^{-1})\gamma\in\frg(R\tl{t})
\end{equation*}
using that $h$ is in the centralizer of $\gamma$. Therefore $[hg]\in \tX_{\gamma}(R)$. The assignment $[g]\mapsto [hg]$ for $h\in LG_{\gamma}$ and $[g]\in\tX_{\gamma}$ defines an action of   $LG_{\gamma}$ on $\tX_{\gamma}$. It induces an action of $LG_{\gamma}$ on the reduced structure $\sX_{\gamma}$.

\subsubsection{The split case} We consider the case where $\gamma\in\frt^{\rs}(F)$. In this case $G_{\gamma}=T\otimes_{k}F$, and
\begin{equation*}
LG_{\gamma}=LT\cong\xcoch(T)\otimes_{\ZZ}L\Gm
\end{equation*}
where $L\Gm$ is the loop group of the multiplicative group $\Gm$. For any $k$-algebra $R$, $L\Gm(R)=R\lr{t}^{\times}$. It is easy to see that an element $a=\sum_{i}a_{i}t^{i}\in R\lr{t}$ is invertible if and only if $a$ starts with finitely many nilpotent coefficients and the first non-nilpotent coefficient is invertible in $R$. When $R$ is reduced, the leading coefficient of $a$ must be invertible in $R$, which implies $R\lr{t}^{\times}=t^{\ZZ}\cdot R\tl{t}^{\times}$, and $R\tl{t}^{\times}=\{(a_{0},a_{1},\cdots)|a_{0}\in R^{\times}, a_{i}\in R, \forall i\geq 1\}$. We see that the reduced ind-scheme $(L\Gm)^{\red}\cong \ZZ\times L^{+}\Gm$, and that $L^{+}\Gm\cong \Gm\times\AA^{\NN}$ as schemes, where $\AA^{\NN}=\Spec k[x_{1},x_{2},\cdots]$. Therefore, when  $\gamma\in\frt^{\rs}(F)$, we have $(LT)^{\red}\cong\xcoch(T)\times L^{+}T$, and $L^{+}T$ is an affine scheme of infinite type.

\subsubsection{The lattice $\L_{\gamma}$}\label{sss:L gamma} For a general regular semisimple $\gamma\in\frg(F)$, let $\xcoch(G_{\gamma}):=\Hom_{F}(\Gm,G_{\gamma})$ be the $F$-rational cocharacter lattice of the torus $G_{\gamma}$. For each $\l\in\xcoch(G_{\gamma})$ viewed as a homomorphism $\Gm\to G_{\gamma}$ defined over $F$, we may consider the element $\l(t)$. The assignment $\l\mapsto \l(t)$ defines an injective homomorphism
\begin{equation*}
\xcoch(G_{\gamma})\incl G_{\gamma}(F).
\end{equation*}
whose image is denoted by $\L_{\gamma}$. It can be shown that the quotient $\L_{\gamma}\bs (LG_{\gamma})^{\red}$ is an affine scheme that is a finite disjoint union of $\Gm^{a}\times \AA^{\NN}$ for some integer $a$.

\subsubsection{The case $G=\GL_{n}$}\label{sss:GLn centralizer} We continue with the setup of \S\ref{sss:ASF Gr GLn}. We assume that $\chark>n$. Then the characteristic polynomial $P(x)=x^{n}+a_{1}x^{n-1}+\cdots+a_{n}\in F[x]$ of $\gamma$ is separable. The $F$-algebra $F[x]/(P(x))$ is then a product of fields $F_{1}\times \cdots \times F_{m}$, with $\sum_{i=1}^{m}[F_{i}:F]=n$. Each field extension $F_{i}/F$ is obtained by adjoining a root of an irreducible factor $P_{i}(x)$ of $P(x)$, and $F_{i}$ is necessarily of the form $k\lr{t^{1/e_{i}}}$ since $\chark>n$. Then the centralizer $G_{\gamma}$ is isomorphic to the product of {\em induced tori}
\begin{equation*}
G_{\gamma}\cong \prod_{i=1}^{m}\Res^{F_{i}}_{F}\Gm.
\end{equation*} 
We have $\xcoch(G_{\gamma})\cong \ZZ^{m}$, and the map $\xcoch(G_{\gamma})\to G_{\gamma}(F)$ is given by
\begin{equation*}
\ZZ^{m}\ni(d_{1},\cdots, d_{m})\mapsto (t^{d_{1}},\cdots, t^{d_{m}})\in F^{\times}_{1}\times\cdots\times F^{\times}_{m}
\end{equation*}
The quotient $\L_{\gamma}\bs G_{\gamma}(F)$ is isomorphic to $\prod_{i=1}^{m}F^{\times}_{i}/t^{\ZZ}$. Since each $F_{i}$ is isomorphic to $k\lr{t^{1/e_{i}}}$, we have  an exact sequence $1\to \calO^{\times}_{F_{i}}\to F^{\times}_{i}/t^{\ZZ}\to \ZZ/e_{i}\ZZ\to 0$, and hence the quotient $\L_{\gamma}\bs (LG_{\gamma})^{\red}$  contains the group scheme $\prod_{i=1}^{m}L^{+}_{F_{i}}\Gm$ with finite index. Here $L^{+}_{F_{i}}\Gm$ is isomorphic to $L^{+}\Gm$ as a scheme, except that we are renaming the uniformizer $t^{1/e_{i}}$. 

Alternatively, we may fix a uniformizer $t_{i}\in F_{i}$ (for example take $t_{i}=t^{1/e_{i}}$) and let $\wt\L_{\gamma}=t_{1}^{\ZZ}\times\cdots\times t^{\ZZ}_{m}\subset \prod_{i}F_{i}^{\times}=G_{\gamma}(F)$. The lattice $\wt\L_{\gamma}$ will be useful in calculating orbital integrals, see \S\ref{sss:GLn orb coho}.

\subsubsection{The case $G=\SL_{2}$} Let $G=\SL_{2}$ and $\gamma=\mat{0}{t^{n}}{1}{0}$ where $n\geq 1$ is odd. Then $G_{\gamma}(F)$ consists of matrices $\mat{a}{bt^{n}}{b}{a}$ with $a,b\in F$ and $a^{2}-t^{n}b^{2}=1$. Note that this equation forces $a,b\in\calO_{F}$, hence $G_{\gamma}(F)=G_{\gamma}(\calO_{F})$. The torus $G_{\gamma}$ is non-split and splits over the quadratic extension $E=F(t^{1/2})$. The lattice $\L_{\gamma}=\Hom_{F}(\Gm,G_{\gamma})=0$. Writing $a=\sum_{i\geq0}a_{i}t^{i}$ and $b_{i}=\sum_{i\geq 0}b_{i}t^{i}$, we see that $a_{0}=\pm1$, and once $b$ and $a_{0}$ are fixed, the higher coefficients of $a$ can be solved uniquely using the Taylor expansion of $(1+t^{n}b^{2})^{1/2}$. Therefore $(LG_{\gamma})^{\red}\cong L^{+}G_{\gamma}$ is isomorphic to $\{\pm1\}\times \AA^{\NN}$, given by $(a,b)\mapsto (a_{0},b_{0},b_{1},\cdots)$.

\subsubsection{Symmetry on affine Springer fibers}\label{sss:sym ASF} Ng\^o has found a more precise statement about the action of $LG_{\gamma}$ on $\tX_{\gamma}$, namely the action factors through a canonical finite-dimensional quotient. We sketch the story following \cite[\S3.3]{NgoFL}. Let $a=\chi(\gamma)\in\frc(F)^{\rs}$ be the image of $\gamma$ under $\chi:\frg\to \frc$. We assume $a\in\frc(\calO_{F})$ for otherwise $\sX_{\gamma}$ is empty. 

There is a smooth affine group scheme $J$ over $\frc$ called the {\em regular centralizer group scheme}. It is characterized by the property that its pullback to $\frg$ via $\chi$, denoted $\chi^{*}J$, maps into the universal centralizer group scheme $I$ over $\frg$, and this map is an isomorphism over the regular locus $\frg^{\reg}$. Let $J_{a}$ be pullback of $J$ under the map $a:\Spec\calO_{F}\to \frc$. Then $J_{a}$ is a smooth affine group scheme over $\calO_{F}$ whose $F$-fiber is the torus $G_{\gamma}$ (i.e., $J_{a}$ is an integral model of $G_{\gamma}$ over $\calO_{F}$). We may form the positive loop group $L^{+}J_{a}$ of $J_{a}$ as well as its affine Grassmannian $P_{a}:=LG_{\gamma}/L^{+}J_{a}$ (also called the {\em local Picard group}). The reduced group scheme $P^{\red}_{a}$ is finite-dimensional and locally of finite type. Ng\^o showed that the action of $LG_{\gamma}$ on $\tX_{\gamma}$ (and hence on $\sX_{\gamma}$) factors through the local Picard group $P_{a}$, and it does not factor through any further quotient. For related statement, see \S\ref{dim formula discussion}.

\subsection{Further examples of affine Springer fibers}\label{ss:aff examples}
In this subsection we give more examples illustrating the rich geometry of affine Springer fibers. We omit the calculations that lead to the geometric descriptions. 

In all examples below, the affine Springer fibers are {\em homogeneous} in the sense that $\sX_{\gamma}$ is equipped an extra $\Gm$-action coming from the {\em loop rotation} on $\Gr_{G}$ by dilation on the uniformizer $t$. For more information on homogeneous affine Springer fibers and their application to representation theory, see \cite{OY}.

\subsubsection{}\label{SL2 ram} When $G=\SL_{2}$ and $\gamma=\mat{0}{t^{m+1}}{t^{m}}{0}$.  When $m$ is even, $\sX_{\gamma}$ consists exactly of those lattices $\L$ such that
\begin{equation*}
t^{m/2}\calO_{F}\oplus t^{m/2}\calO_{F}\subset\L\subset t^{-m/2}\calO_{F}\oplus t^{-m/2}\calO_{F}.\end{equation*}
In this case, $\sX_{\gamma}$ coincides with the closure of the $L^{+}G$-orbit in $\Gr_{G}$ corresponding to the coweight $m\alpha^{\vee}/2$. 

When $m$ is odd, $\sX_{\gamma}$ consists exactly of those lattices $\L$ such that
\begin{equation*}
t^{(m+1)/2}\calO_{F}\oplus t^{(m-1)/2}\calO_{F}\subset\L\subset t^{-(m-1)/2}\calO_{F}\oplus t^{-(m+1)/2}\calO_{F}.
\end{equation*}
In this case, consider instead the affine Grassmannian $\Gr_{G'}$ of $G'=\PGL_{2}$, which contains $\Gr_{G}$ as a component. Then $\sX_{\gamma}$ can be identified with the closure of the $L^{+}G'$-orbit in $\Gr_{G'}$ corresponding to the coweight $m\alpha^{\vee}/2$. 

In either case, we have $\dim\sX_{\gamma}=m$.

\subsubsection{The Lusztig-Smelt examples}\label{sss:LS} Let $G=\GL_{n}$ and $\gamma\in\frg(F)$ with characteristic polynomial $P(x)=x^{n}-t^{m}=0$, where $(m,n)=1$. If a lattice $\L\subset F^{n}$ is stable under $\gamma$, it carries an action of the ring $R=\calO_{F}[X]/(X^{n}-t^{m})$. Let $K=\Frac(R)$, then $K= k\lr{s}$ with $x=s^{m}$ and $t=s^{n}$. Then the integral closure of $R$ in $K$ is $\tilR:=k\tl{s}$. The action of $\gamma$ on $F^{n}$ makes it a one-dimensional $K$-vector space. We fix a $K$-linear isomorphism $F^{n}\cong K$, under which a lattice in $F^{n}$ stable under $R$ is simply a {\em fractional $R$-ideal}, i.e., a finitely generated $R$-submodule of $K=\Frac(R)$.  Then $\sX_{\gamma}$ can be identified with the moduli space of fractional $R$-ideals.

The centralizer $G_{\gamma}(F)$ is simply $K^{\times}$, which acts on the set of fractional $R$-ideals by multiplication. This action clearly factors through $K^{\times}/R^{\times}$, which is the group of $k$-points of the local Picard group scheme $P_{\gamma}$.

There is an action of $\Gm$ on $K$ (by field automorphisms) given by scaling $s$ (so $s^{i}$ gets weight $i$ under this action). This induces an action of $\Gm$ on $\sX_{\gamma}$. The fixed points of $\Gm$ on $\sX_{\gamma}$ correspond to fractional ideals generated by monomials of $s$. More precisely, if a fractional $R$-ideal $\L\subset K=k\lr{s}$ is fixed by $\Gm$,  define $M_{\L}=\{i\in\ZZ|s^{i}\in \L\}$ which is a subset of $\ZZ$ stable under adding $m$ and $n$, because $\L$ is an $R=k\tl{s^{m},s^{n}}$-module. Therefore $M_{\L}\subset\ZZ$ is a finitely generated module for the monoid $A_{m,n}:=\ZZ_{\geq 0}m+\ZZ_{\geq0}n\subset\ZZ_{\geq 0}$. The assignment $\L\mapsto M_{\L}$ gives a bijection
\begin{equation*}
\sX_{\gamma}^{\Gm}(k)\isom\{\mbox{$A_{m,n}$-submodules $M\subset \ZZ$}\}.
\end{equation*}
Any $A_{m,n}$-submodule of $\ZZ$ contains all sufficiently large integers. Therefore any two such $A_{m,n}$-module $M$ and $M'$ differ by finitely many elements, and we can define $[M:M']=\#(M\bs M')-\#(M'\bs M)$. Fox any fixed $i\in\ZZ$, we have a total of $\frac{1}{n+m}\binom{n+m}{n}$ fixed points with $[M:\ZZ_{\geq0}]=i$. For a fixed point $p_{M}$ corresponding to an $A_{m,n}$-module $M$, consider the subvariety $C_{M}=\{p\in\sX_{\gamma}|\lim_{\Gm\ni z\to0}z\cdot p=p_{M}\}$. Then $C_{M}$ is isomorphic to an affine space whose dimension can be expressed combinatorially in terms of $M$. The cells $C_{M}$ give a stratification of $\sX_{\gamma}$. This gives a way to compute the Poincar\'e polynomial of connected components of $\sX_{\gamma}$. For more details, and the similar picture for $\sY_{\gamma}$, see \cite{LS}.

\subsubsection{} We look at the geometry of $\sY_{\gamma}$ in a special case of \S\ref{sss:LS}. We consider the case $G=\GL_{3}$ and $\gamma^{3}-t^{2}=0$. Introducing the variable $s$ with $t=s^{3}$ and $\gamma=s^{2}$ as before, then $R=k\tl{s^{2},s^{3}}$ with fraction field $K=k\lr{s}$. The affine Springer fiber  $\sY_{\gamma}$ classifies a chain of fractional $R$-ideals $\L_{0}\subset \L_{1}\subset \L_{2}\subset s^{-3}\L_{0}$.  We consider a component of $\sY^{0}_{\gamma}\subset\sY_{\gamma}$, classifying chains $\{\L_{i}\}$ as above with $[\L_{i}:k\tl{s}]=i$, $0\leq i\leq 2$. 

We first study the $\Gm$-fixed points on $\sY^{0}_{\gamma}$. For each $\Gm$-fixed $R$-submodule of $k\lr{s}$, we denote its associated module for the monoid $A_{3,2}$ (see \S\ref{sss:LS}) by a sequence of integers. For example, $(0,1,2,\cdots)$ stands for the standard lattice $k\tl{s}$. Note that the sequence for any $\Gm$-fixed $R$-fractional ideal is either consecutive $(n,n+1,n+2,\cdots)$ or has at most one gap at the second place $(n,n+2,n+3,\cdots)$. We have the following four fixed points in $\sY^{0}_{\gamma}$:
\begin{enumerate}
\item $q: (0,1,2,\cdots)\subset (-1,0,1,\cdots)\subset (-2,-1,0,\cdots)$;
\item $p_{0}: (-1,1,2,\cdots)\subset (-1,0,1,\cdots)\subset (-2,-1,0,\cdots)$;
\item $p_{1}: (0,1,2,\cdots)\subset (-2,0,1,\cdots)\subset (-2,-1,0,\cdots)$;
\item $p_{2}:(0,1,2,\cdots)\subset (-1,0,1,\cdots)\subset (-3,-1,0,\cdots)$. 
\end{enumerate}
Our indexing scheme is that $p_{i}$ is obtained from $q$ by changing the lattice $\L_{i}$. 

Then $\sY^{0}_{\gamma}$ is the union of three irreducible components $C_{0}\cup C_{1}\cup C_{2}$, and each component $C_{i}$ is isomorphic to $\PP^{1}$. They all contain $q$ and that is the only intersection between any two of them. We have $p_{i}\in C_{i}$ for $i=0,1,2$.

There is a natural way to index affine partial flag varieties  of $G$ by subsets $J\subset \{0,1,2\}$, as we saw in \S\ref{sss:GLn parahoric}. Let $\sX_{J,\gamma}$ be the affine Springer fiber of $\gamma$ in $\Fl_{J}$. Under the projection $\sY_{\gamma}\to \sX_{J,\gamma}$, the curves $C_{i}$ for $i\notin J$ collapse to a point, and the other curves map isomorphically onto their images.

\subsubsection{}\label{ex:Sp} Let $G=\Sp(V)$ where $V$ is a symplectic space of dimension $2n$ over $k$, and assume $\chark\neq2$. Fix a decomposition $V=U\oplus U^{*}$ into Lagrangian subspaces of $V$, such that the symplectic form restricts to the natural pairing on $U\times U^{*}$. Consider $\gamma=\mat{0}{tX}{Y}{0}\in\frg(F)$ where $X\in\Sym^{2}(U)$ (viewed as a self-adjoint map $\xi:U^{*}\to U$) and $Y\in\Sym^{2}(U^{*})$ (viewed as a self-adjoint map $\eta:U\to U^{*}$). The condition that $\gamma$ is regular semisimple is equivalent to that: (1) both $\xi$ and $\eta$ are isomorphisms; (2) $\xi\eta\in\GL(U)$ is regular semisimple, or equivalently $\eta\xi\in\GL(U^{*})$ is regular semisimple.

The affine Springer fiber $\sX_{\gamma}$ classifies self-dual lattices $\L\subset V\otimes_{k}F$ which are stable under $\gamma$ (see \S\ref{sss:ASF Gr GLn}). Consider the action of $\Gm$ on $V\otimes_{k}F=U\otimes_{k}F\oplus U^{*}\otimes_{k}F$ such that $Ut^{i}$ has weight $2i-1$ and $U^{*}t^{i}$ has weight $2i$. This induces a $\Gm$-action on $\Gr_{G}$ and on $\sX_{\gamma}$. We first consider the fixed points $\Gr_{G}^{\Gm}$. A lattice $\L\in\Gr_{G}$ is fixed under $\Gm$ if and only if it is the $t$-adic completion of
\begin{equation*}
\Span\{(A_{i}\oplus B_{i})t^{i}; i\in\ZZ \}
\end{equation*}
where $\cdots \subset A_{-1}\subset A_{0}\subset A_{1}\subset \cdots \subset U$ is a filtration of $U$ such that $A_{N}=0$ for $N\ll 0$ and $A_{N}=U$ for $N\gg 0$, and $\cdots \subset B_{-1}\subset B_{0}\subset B_{1}\subset \cdots \subset U^{*}$ is a similar filtration of $U^{*}$, such that
\begin{equation}\label{AB dual}
B_{i-1}=A^{\bot}_{-i}, \forall i\in\ZZ
\end{equation}
under the duality pairing between $U$ and $U^{*}$. The last condition reflects the fact that $\L$ is self-dual  under the symplectic form.

A lattice $\L\in\sX_{\gamma}^{\Gm}$ is then determined by two filtrations $A_{\bu}$ of $U$ and $B_{\bu}$ of $U^{*}$, dual in the sense \eqref{AB dual}, with the extra condition that
\begin{equation*}
\eta A_{i}\subset B_{i}; \quad \xi B_{i}\subset A_{i+1}, \quad\forall i\in\ZZ.
\end{equation*}
We summarize the data into the following diagram
\begin{equation*}
\xymatrix{ A_{-1}\ar@{^{(}->}[rr]\ar[dr]^{\eta} && A_{0}\ar@{^{(}->}[rr]\ar[dr]^{\eta} && A_{1}\ar[dr]^{\eta}& \cdots\\
\cdots & A^{\bot}_{0}\ar@{^{(}->}[rr]\ar[ur]^{\xi} && A^{\bot}_{-1}\ar@{^{(}->}[rr]\ar[ur]^{\xi} && A_{-2}^{\bot}}
\end{equation*}

For example, when $\dim V=4$, hence $\dim U=2$, there are two possibilities. The first possibility is $A_{-1}=0\subset A_{0}=U$, which corresponds to the standard lattice $V\otimes_{k}\calO_{F}$. The second possibility is $A_{-1}=0\subset A_{0}\subset A_{1}=U$, and $A_{0}\subset U$ is a line satisfying $\xi A^{\bot}_{0}\subset A_{0}$, i.e., $A^{\bot}_{0}\subset U^{*}$ is an isotropic line under the quadratic form $X$. There are two such lines $A^{\bot}_{0}$, giving two other $\Gm$-fixed points. Therefore, $\sX_{\gamma}^{\Gm}$ consists of 3 points. We have $\sX_{\gamma}=C_{1}\cup C_{2}$ where $C_{i}\cong\PP^{1}$, the two components intersect at the standard lattice, and each $C_{i}$ contains one of the remaining $\Gm$-fixed points.

\subsubsection{The Bernstein-Kazhdan example}\label{sss:BK example} In \cite[Appendix]{KL}, Bernstein and Kazhdan gave the first example of an irreducible component of an affine Springer fiber which was not a rational variety. We keep the same notation as in \S\ref{ex:Sp}. Let $\Fl_{\bP}$ be the partial affine flag variety of $G=\Sp(V)$ classifying pairs of lattices $\L'\subset\L$ such that $\L^{\vee}=\L$ and $\L'^{\vee}=t^{-1}\L$.
Let $\gamma$ be as in \S\ref{ex:Sp}. Then the same $\Gm$ acts on $\sX_{\bP,\gamma}$, and the fixed points $\L'\subset \L$ can be described by two pairs of filtrations $(A_{\bu}, B_{\bu})$ and $(A'_{\bu}, B'_{\bu})$, where $(A_{\bu}, B_{\bu})$ is the kind of filtration of $U$ and $U^{*}$ as described in \S\ref{ex:Sp}, and $(A'_{\bu}, B'_{\bu})$ is similar except that \eqref{AB dual} is replaced by $B'_{i}=A_{-i}'^{\bot}$. Moreover, the inclusion $\L'\subset \L$ is equivalent to $A'_{i}\subset A_{i}$ and $B'_{i}\subset B_{i}$, for all $i$.

Consider for example $\dim V=6$ and we fix the dimension of the filtrations:
\begin{equation*}
A_{-1}=0, \quad \dim A_{0}=2, \quad A_{1}=U;\quad  \dim A'_{0}=1, \quad A'_{1}=U.
\end{equation*}
Such filtrations $(A_{\bu}, B_{\bu}; A'_{\bu}, B'_{\bu})$ are determined by the complete flag $0\subset A'_{0}\subset A_{0}\subset U$ satisfying $\xi A^{\bot}_{0}\subset A_{0}$ and that $\eta A'_{0}\subset A'^{\bot}_{0}$. In other words,  $A^{\bot}_{0}$ is an isotropic line in $U^{*}$ under the quadratic form $X$, and $A'_{0}$ is an isotropic line in $U$ under the quadratic form $Y$. The pair $(A^{\bot}_{0}, A'_{0})$ determines a point in $Q(X)\times Q(Y)\subset \PP(U)\times\PP(U^{*})$, the product of conics defined by $X$ and $Y$. The incidence relation $A'_{0}\subset A_{0}$ defines a curve of bidegree $(2,2)$ in $Q(X)\times Q(Y)\cong\PP^{1}\times \PP^{1}$, which is then a curve of genus one. Therefore, a connected component of $\sX^{\Gm}_{\bP,\gamma}$ is a curve of genus  one. Consider the points in $\sX_{\bP,\gamma}$ that contract to this curve, and take its closure $Z$. One can show that $\dim Z=3=\dim\sX_{\bP,\gamma}$. Hence $\sX_{\bP,\gamma}$ contains an irreducible component $Z$ which is irrational. We refer to the appendix of \cite{KL} for more details.

\subsubsection{``Subregular'' affine Springer fibers}
When $\sX_{\gamma}$ or $\sY_{\gamma}$ is one-dimensional, we may call them subregular affine Springer fibers, by analog with subregular Springer fibers discussed in \S\ref{sss:subreg}. If $\dim\sY_{\gamma}=1$, it is a union of $\PP^{1}$'s, hence we can define its dual graph. In \cite[Prop 7.7]{KL}, the dual graphs of the subregular affine Springer fibers in $\Fl$ are classified, and they are almost always the extended Dynkin diagrams of simply-laced groups, except that they can also be infinite chains in type $A$ (see Example \ref{P1 chain}).

\subsection{Geometric Properties of affine Springer fibers}

\subsubsection{Non-reducedness} The ind-scheme $\tX_{\gamma}$ is {\em never} reduced if $G$ is nontrivial and $\gamma$ is regular semisimple in $\frg(F)$. For example, in the case considered in \S\ref{sss:reduction rs}, $\tX_{\gamma}$ is isomorphic to $\Gr_{T}=LT/L^{+}T$. We have seen in \S\ref{sss:centralizer} that for non-reduced rings $R$, elements in $LT(R)=R\lr{t}^{\times}$ can have nilpotent leading coefficients. Therefore $\Gr_{T}(R)$ is not just $\xcoch(T)$, which is $\Gr^{\red}_{T}$. This shows that $\Gr_{T}$ is non-reduced, hence $\tX_{\gamma}$ is non-reduced.

The next theorem is the fundamental finiteness statement about $\sX_{\gamma}$.

\begin{theorem}[Essentially Kazhdan and Lusztig {\cite[Prop 2.1]{KL}}]\label{th:ft quot} Let $\gamma\in \frg(F)$ be a regular semisimple element. Then
\begin{enumerate}
\item There exists a closed subscheme $Z\subset \sX_{\gamma}$ which is projective over $k$, such that $\sX_{\gamma}=\cup_{\ell\in\L_{\gamma}}\ell\cdot Z$.
\item The ind-scheme $\sX_{\gamma}$ is a scheme locally of finite type over $k$. 
\item The action of $\L_{\gamma}$ on $\sX_{\gamma}$ is free, and the quotient $\L_{\gamma}\bs\sX_{\gamma}$ (as an fppf sheaf on $k$-algebras) is representable by a proper algebraic space over $k$.
\end{enumerate}
\end{theorem}

We sketch a proof of this theorem below in three steps.

\subsubsection{First reduction} We show that part (1) of the theorem implies (2) and (3). Let $Z$ be a projective subscheme as in (1). To show (2), we would like to show that any $x\in \sX_{\gamma}(k)$ has an open neighborhood which is a scheme of finite type. By the $\L_{\gamma}$-action we may assume $x\in Z(k)$. Since $Z$ is of finite type, the set $\Sigma:=\{\ell\in \L_{\gamma}|Z\cap\ell\cdot Z\neq\varnothing\}$ is finite. Let $U=\sX_{\gamma}-\cup_{\ell\notin\Sigma}\ell\cdot Z$, then $U$ is an open neighborhood of $Z$, hence an open  neighborhood of $x$. Moreover, $U$ is contained in the finite union $\cup_{\ell\in\Sigma}\ell\cdot Z$, hence contained in some Schubert variety $\Gr_{G,\leq\l}$. Hence $U$ is an open subset of the projective scheme $\sX_{\gamma}\cap\Gr_{G,\leq\l}$, therefore $U$ is itself a scheme of finite type. To show (3), note that the fppf sheaf quotient $\L_{\gamma}\bs\sX_{\gamma}$ is a separated algebraic space because it is the quotient of $\sX_{\gamma}$ by the \'etale equivalence relation $\L_{\gamma}\times\sX_{\gamma}\subset \sX_{\gamma}\times\sX_{\gamma}$ (given by the action and projection maps). By (1), there is a surjection $Z\to \L_{\gamma}\bs\sX_{\gamma}$ from a projective scheme $Z$, which implies that $\L_{\gamma}\bs\sX_{\gamma}$ is proper.

\subsubsection{Proof of (1) when $\gamma$ lies in a split torus}\label{sss:split case} We first consider the case where $\gamma$ lies in a split torus. By $G(F)$-conjugation, we may assume $\gamma\in\frt(F)$. In this case $\L_{\gamma}\cong\xcoch(T)$. Fix a Borel subgroup $B$ containing $T$ and let $N$ be the unipotent radical of $B$. The Iwasawa decomposition of $LG$ gives
\begin{equation*}
\Gr_{G}=LN\cdot \L_{\gamma}L^{+}G/L^{+}G=\sqcup_{\l\in\xcoch(T)}LN\cdot \l(t)L^{+}G/L^{+}G.
\end{equation*}
Let $X:=\sX_{\gamma}\cap (LN\cdot L^{+}G/L^{+}G)\subset \sX_{\gamma}$. It is enough to show that $X$ lies in some affine Schubert variety $\Gr_{\leq \l}$, for then its closure $Z:=\ov{X}$ in $\Gr_{\leq\l}$ satisfies the condition in (1). For later use, we note that the translations $\ell\cdot X$ for $\ell\in \L_{\gamma}$ are disjoint and cover $\sX_{\gamma}$. 

Fixing an ordering of the positive roots of $T$ with respect to $B$, we may write an element  $u\in N(F)$ uniquely as
\begin{equation}\label{root gp}
u=\prod_{\alpha>0}x_{\alpha}(c_{\alpha})
\end{equation}
where $c_{\alpha}\in F$ and $x_{\alpha}: \Ga\to N$ is the root group corresponding to $\alpha$. To show that $X$ is contained in an affine Schubert variety, it suffices to give a lower bound for the valuations of $c_{\alpha}$ appearing in \eqref{root gp} for any $u\in N(F)$ such that $[u]\in X(k)$. We may expand $\Ad(u^{-1})\gamma\in\frb(F)$ in terms of the root decomposition
\begin{equation*}
\Ad(u^{-1})\gamma=\gamma+\sum_{\alpha>0}P_{\alpha}(\gamma; c)e_{\alpha}
\end{equation*} 
where $e_{\alpha}\in\frg_{\alpha}$ is a fixed basis for each root space and $P_{\alpha}(\gamma;c)$ is an $F$-valued polynomial function in $\{c_{\alpha}\}_{\alpha>0}$ and linear in $\gamma$. Induction on the height of $\alpha$ shows that $P_{\alpha}(\gamma;c)$ takes the following form
\begin{equation}\label{conj u gamma}
P_{\alpha}(\gamma;c)=\jiao{\alpha,\gamma}c_{\alpha}+\sum_{\beta<\alpha}\jiao{\beta,\gamma}P^{\beta}_{\alpha}(c)
\end{equation}
where $P^{\beta}_{\alpha}(c)$ is a polynomial involving only $\{c_{\alpha'}\}_{\alpha'<\alpha}$, and homogeneous of degree $\alpha$ (we define $\deg c_{\alpha'}:=\alpha'\in\xch(T)$). 

Let $n=\max_{\alpha>0}\{\val_{F}\jiao{\alpha,\gamma}\}$. This is finite because $\gamma$  is regular semisimple. If $[u]\in X(k)$, i.e., $\Ad(u^{-1})\gamma\in\frg(\calO_{F})$, then induction on the height of $\alpha$ shows that
\begin{equation*}
\val_{F}(c_{\alpha})\geq -(2\textup{ht}(\alpha)-1)n,
\end{equation*} 
which gives the desired lower bound and shows that $X$ lies in an affine Schubert variety.

\subsubsection{Proof of (1) in the general case} In the general case, we give a simplified argument compared to the original one in \cite{KL}, following the same idea. We make a base change to $F'=k\lr{t^{1/m}}$ over which $\gamma$ can be conjugated into a split torus. Let $\Gr'_{G}$ be the affine Grassmannian of $G$ defined using the field $F'$ in place of $F$ (so that $\Gr'_{G}(k)=G(F')/G(\calO'_{F})$), and let $\sX'_{\gamma}\subset\Gr'_{G}$ be the corresponding affine Springer fiber. Then both $\Gr'_{G}$ and $\sX'_{\gamma}$ carry an action of $\Gamma:=\Gal(F'/F)\cong\mu_{m}$ induced from its action on $F'$, and we have a closed embedding $\sX_{\gamma}\incl (\sX'_{\gamma})^{\Gamma}$. 

Let $\L'_{\gamma}=\xcoch(G_{\gamma}\otimes F')\incl G_{\gamma}(F')$ be the lattice constructed using the field $F'$. There is an action of $\Gamma$ on $\xcoch(G_{\gamma}\otimes F')$ with fixed points $\xcoch(G_{\gamma})\cong \L_{\gamma}$. This action induces an action of $\Gamma$ on $\L'_{\gamma}$, but it does not respect the embedding $\xcoch(G_{\gamma}\otimes F')\incl G_{\gamma}(F')$. The fixed points $(\L'_{\gamma})^{\Gamma}=\xcoch(G_{\gamma}\otimes F')^{\Gamma}$ may not lie in $G_{\gamma}(F)$, however it always contains $\L_{\gamma}$ with finite index. 

From the proof in the split case in \S\ref{sss:split case} we have a finite type locally closed subscheme $X'\subset\sX'_{\gamma}$ coming from the Iwasawa decomposition, such that $\sX'_{\gamma}$ can be decomposed as the {\em disjoint union} $\sqcup_{\ell\in\L'_{\gamma}}\ell\cdot X'$  (not as schemes but as constructible sets). We identify $\ell\in \L'_{\gamma}$ with an element in $\xcoch(G_{\gamma}\otimes F')$ then $\sigma(\ell)$ makes sense for $\sigma\in\Gamma$. One checks that $\sigma(\ell\cdot X')=\sigma(\ell)\cdot X'$ for $\ell\in\L'_{\gamma}$ even though $\sigma$ does not respect the embedding $\L'_{\gamma}\to G_{\gamma}(F')$. Therefore $(\sX'_{\gamma})^{\Gamma}\subset (\L'_{\gamma})^{\Gamma}\cdot X'$. Choosing representatives $C$ for the finite coset space $(\L'_{\gamma})^{\Gamma}/\L_{\gamma}$, we see that $(\sX'_{\gamma})^{\Gamma}$ is contained in $\L_{\gamma}\cdot (C\cdot X')$. Since $C\cdot X'$ is of finite type, $Z'=\ov{C\cdot X'}\cap (\sX'_{\gamma})^{\Gamma}$ is a projective subscheme of $\sX'_{\gamma}$ whose $\L_{\gamma}$-translations cover $(\sX'_{\gamma})^{\Gamma}$. Finally the projective subscheme $Z=Z'\cap \sX_{\gamma}$ of $\sX_{\gamma}$ satisfies the requirement of (1).  This finishes the proof of Theorem \ref{th:ft quot}. \qed

\subsubsection{Reduction to Levi} The proof of Theorem \ref{th:ft quot} in the split case in \S\ref{sss:split case} gives more information. In the Iwasawa decomposition, let $S_{\l}\subset \Gr_{G}$ be the $LN$-orbit of $\l(t)$, for $\l\in\xcoch(T)$. This is called a {\em semi-infinite orbit}, because it has infinite dimension and also has infinite codimension in $\Gr_{G}$.  Let $C_{\l}:=\sX_{\gamma}\cap S_{\l}$, then $C_{\l}=\l(t)\cdot X$ in the notation of \S\ref{sss:split case}. The formula \eqref{conj u gamma} implies that $C_{\l}\neq\varnothing$ (or equivalently $X\neq\varnothing$, or $\sX_{\gamma}\neq\varnothing$) if and only if $\jiao{\alpha,\gamma}\in\calO_{F}$ for all roots $\alpha$, and if so, $C_{\l}$ is isomorphic to an almost affine space (namely an iterated $\AA^{1}$-bundle) of dimension
\begin{equation*}
\dim C_{\l}=\sum_{\alpha>0}\val_{F}(\jiao{\alpha,\gamma})=\frac{1}{2}\val_{F}\Delta(\gamma)
\end{equation*}
Here $\Delta(\gamma)$ is the determinant of the adjoint action of $\gamma$ on $\frg(F)/\frt(F)$. Therefore $\sX_{\gamma}$ can be decomposed into almost affine spaces of the same dimension indexed by $\xcoch(T)$. However, this decomposition is not a stratification: the closure $\ov{C}_{\l}$ of $C_{\l}$ will intersect other $C_{\l'}$ but certainly not a union of such $C_{\l'}$'s. 

The decomposition $\sX_{\gamma}=\sqcup C_{\l}$ in the split case has a generalization. Suppose $P$ is a parabolic subgroup of $G$ with unipotent radical $N_{P}$ and a Levi subgroup $M_{P}$.  Let $\mathfrak{m}_{P}=\Lie M_{P}$ and suppose $\gamma\in\mathfrak{m}_{P}(F)$ is regular semisimple as an element in  $\frg(F)$.  Using the generalized Iwasawa decomposition $G(F)=N_{P}(F)M_{P}(F)G(\calO_{F})$, there is a well-defined map $\Gr_{G}(k)\to \Gr_{M_{P}}(k)$ sending $nm G(\calO_{F})$ to $m M_{P}(\calO_{F})$, for $n\in N_{P}(F)$ and $m\in M_{P}(F)$. However this map does not give a map of ind-schemes. Nevertheless the fibers of this map have natural structure of infinite dimensional affine spaces. Restricting this map to $\sX_{\gamma}$ we get
\begin{equation*}
\tau: \sX_{\gamma}(k)\to \sX^{M_{P}}_{\gamma}(k)
\end{equation*}
where $\sX^{M_{P}}_{\gamma}\subset \Gr_{M_{P}}$ is the affine Springer fiber for $\gamma$ and the group $M_{P}$. The fibers of $\tau$, if non-empty, are almost affine spaces of dimension $\frac{1}{2}\val_{F}\Delta^{G}_{M_{P}}(\gamma)$, where
\begin{equation}\label{Delta G L}
\Delta^{G}_{M_{P}}(\gamma):=\det(\ad(\gamma)|\frg(F)/\mathfrak{m}_{P}(F)).
\end{equation}
Assume $T\subset M_{P}$, then the connected components of $\Gr_{M_{P}}$ are indexed by $\xcoch(T)/R^{\vee}_{M_{P}}$ where $R^{\vee}_{M_{P}}$ is the coroot lattice of $M_{P}$. If we decompose $\sX^{M_{P}}_{\gamma}$ into connected components $\sX^{M_{P}}_{\gamma}(\l)$ for $\xcoch(T)/R^{\vee}_{M_{P}}$ and taking their preimages under $\tau$, we get a decomposition $\sX_{\gamma}=\sqcup\sX_{\gamma,\l}$ into locally closed sub-ind-schemes indexed also by $\l\in\xcoch(T)/R^{\vee}_{M_{P}}$. One can show that $\sX_{\gamma,\l}$, if non-empty, is an almost affine space bundle over $\sX^{M_{P}}_{\gamma}(\l)$ with fiber dimension $\frac{1}{2}\val_{F}\Delta^{G}_{M_{P}}(\gamma)$.

\subsubsection{Connectivity and equidimensionality}\label{sss:equidim} When $G$ is simply-connected,  $\Fl$ is connected, and in this case the affine Springer fiber $\sY_{\gamma}$ is also connected. See \cite[\S4, Lemma 2]{KL}. As a consequence, when $G$ is simply-connected,  $\sX_{\bP,\gamma}$ is connected for all parahoric $\bP$ because the natural projection $\sY_{\gamma}\to \sX_{\bP,\gamma}$ is surjective. 

In \cite{KL}, it is also shown that $\sY_{\gamma}$ is equidimensional. The argument there is similar to Spaltenstein's the proof of the connectivity and equidimensionality for Springer fibers in \cite{Spa77}.

\subsubsection{The dimension formula}
By Theorem \ref{th:ft quot}, the dimension of $\sX_{\gamma}$ is well-defined, and is the dimension of $L_{\gamma}\bs\sX_{\gamma}$ as an algebraic space. To state a formula for $\dim\sX_{\gamma}$, we need some more notation.

Consider the adjoint action $\ad(\gamma):\frg(F)\to \frg(F)$. The kernel of this map is $\frg_{\gamma}(F)$, and the induced endomorphism $\ov{\ad(\gamma)}$ on $\frg(F)/\frg_{\gamma}(F)$ is invertible. Let $\Delta(\gamma)\in F^{\times}$ be the determinant of $\ov{\ad(\gamma)}$. This is consistent with our earlier definition of $\Delta(\gamma)$ in the case $\gamma$ lies in a split torus $\frt(F)$.

On the other hand, recall $\xcoch(G_{\gamma})$ is the group of $F$-rational cocharacters of $G_{\gamma}$, which is also the rank of the maximal $F$-split subtorus of $G_{\gamma}$. Let
\begin{equation*}
c(\gamma)=r-\rk_{\ZZ}\xcoch(G_{\gamma}).
\end{equation*}
Then $c(\gamma)$ is also the rank of the maximal $F$-anisotropic subtorus $G_{\gamma}$.

\begin{theorem}[Bezrukavnikov {\cite{Bez}}, conjectured by Kazhdan-Lusztig \cite{KL}]\label{th:aff dim} Let $\gamma\in \frg(F)$ be a regular semisimple element. Then we have 
\begin{equation}\label{aff Spr dim}
\dim\sX_{\gamma}=\frac{1}{2}(\val_{F}\Delta(\gamma)-c(\gamma)).
\end{equation}
\end{theorem}

\subsubsection{Sketch of proof}\label{dim formula discussion} A key role in the proof is played by the notion of {\em regular points} of $\sX_{\gamma}$. We have an evaluation map $\ev: \sX_{\gamma}\to [\frg/G]$ sending $[g]\in\sX_{\gamma}\subset \Gr_{G}$ to the reduction of $\Ad(g^{-1})\gamma$ modulo $t$, which is well-defined up to the adjoint action by $G$.  We say $[g]\in \sX_{\gamma}$ is a {\em regular point} if $\ev([g])$ lies in the open substack $[\frg^{\reg}/G]$ of $[\frg/G]$. Let $\sX^{\reg}_{\gamma}\subset\sX_{\gamma}$ be the open sub-ind-scheme of $\sX_{\gamma}$ consisting of regular points. It can be shown that $\sX^{\reg}_{\gamma}$ is non-empty. Denote the preimage of $\sX^{\reg}_{\gamma}$ in $\sY_{\gamma}$ by $\sY^{\reg}_{\gamma}$, then the projection map $\sY^{\reg}_{\gamma}\to\sX^{\reg}_{\gamma}$ is an isomorphism. Since $\sY^{\reg}_{\gamma}$ is equidimensional as mentioned in \S\ref{sss:equidim}, we see that $\dim\sX^{\reg}_{\gamma}=\dim\sY^{\reg}_{\gamma}=\dim\sY_{\gamma}$. Of course we have $\dim\sX_{\gamma}\leq\dim \sY_{\gamma}$, therefore we must have $\dim \sX^{\reg}_{\gamma}=\dim \sX_{\gamma}=\dim \sY_{\gamma}$. It remains to calculate the dimension of $\sX^{\reg}_{\gamma}$.

Recall we have defined the local Picard group $P_{a}$ in \S\ref{sss:sym ASF}. The action of $P_{a}$ on $\sX_{\gamma}$ preserves $\sX^{\reg}_{\gamma}$, and in fact $\sX^{\reg}_{\gamma}$ is a torsor under $P_{a}$.  Therefore it suffices to compute the dimension of $P_{a}=LG_{\gamma}/L^{+}J_{a}$. 

Consider the projection $\chi_{\frt}: \frt\to\frc=\frt\sslash W$. Pullback along $a=\chi(\gamma):\Spec\calO_{F}\to\frc$ we get a finite morphism $\chi_{a}: \chi^{-1}_{\frt}(a)=\Spec A\to \Spec\calO_{F}$, for a finite flat $\calO_{F}$-algebra $A$. There is a close relationship between the group scheme $J_{a}$ and $J'_{a}=(\Res^{A}_{\calO_{F}}(T\otimes_{k}A))^{W}$, where the Weyl group $W$ acts diagonally on both $T$ and $A$.  In fact $J_{a}$ and $J'_{a}$ are equal up to connected components in their special fibers. In particular, $\dim P_{a}=\dim P'_{a}$, where $P'_{a}=LG_{\gamma}/L^{+}J'_{a}$. It is not hard to see that
\begin{equation*}
\dim P'_{a}=\dim (\frt\otimes \tilA)^{W}-\dim (\frt\otimes A)^{W}=\dim (\frt\otimes (\tilA/A))^{W}
\end{equation*}
where $\tilA$ is the normalization of $A$. From this one deduces the dimension formula \eqref{aff Spr dim}. \qed

\begin{remark} Using the relation between affine Springer fibers and Hitchin fibers, Ng\^o \cite[Cor 4.16.2]{NgoFL} showed that $\sX^{\reg}_{\gamma}$ is in fact dense in $\sX_{\gamma}$. In particular, each irreducible component of $\sX_{\gamma}$ is a rational variety. However, this rationality property is false for affine Springer fibers in more general affine partial flag varieties, as we saw in Bernstein-Kazhdan's example in \S\ref{sss:BK example}. 
\end{remark}

\subsubsection{Purity} It was conjectured by Goresky, Kottwitz and MacPherson \cite{GKM} that the cohomology of affine Springer fibers should be pure (in the sense of Frobenius weights if $k=\overline{\FF}_{p}$, or in the sense of Hodge structures if $k=\CC$). The purity of affine Springer fibers would allow the authors of \cite{GKM} to prove the Fundamental Lemma for unramified elements using localization techniques in equivariant cohomology. This purity conjecture is still open in general. In \cite{GKMPurity}, a class of affine Springer fibers called {\em equivalued} were shown to be pure.

\subsubsection{Invariance under perturbation}  Suppose $a,a'\in\frc(\calO_{F})\cap\frc^{\rs}(F)$. We say $a\equiv a'\mod t^{N}$ if $a$ and $a'$ have the same image under the map $\frc(\calO_{F})\to\frc(\calO_{F}/t^{N})$. In \cite[Prop 3.5.1]{NgoFL}, it is shown that for fixed $a\in\frc(\calO_{F})\cap\frc^{\rs}(F)$, there exists some $N$ (depending on $a$) such that whenever $a'\equiv a\mod t^{N}$, we have isomorphisms
\begin{equation*}
\tX_{a}\cong\tX_{a'} \textup{ and } P_{a}\cong P_{a'}
\end{equation*}
in a way compatible with the actions. Therefore, we may say that $\sX_{a}$ varies locally constantly with $a$ under the $t$-adic topology on $\frc(\calO_{F})$. 

For example,  consider the case $G=\GL_{n}$ and let $a\in\frc(\calO_{F})$ correspond to a characteristic polynomial $P(x)=x^{n}+a_{1}x^{n-1}+\cdots+a_{n}$ whose roots are in $\calO_{F}$ and are distinct modulo $t$. Then for any $a'\equiv a\mod t$, the characteristic polynomial of $a'$ also has distinct roots modulo $t$. In this case, $\tX_{a}$ and $\tX_{a'}$ are both torsors under $P_{a}\cong P_{a'}\cong \Gr_{T}$.

\subsection{Affine Springer representations}
In this subsection we introduce an analog of Springer's $W$-action on $\cohog{*}{\sB_{e}}$ in the affine situation.  

\subsubsection{The affine Weyl group} We view $W$ as a group of automorphisms of the cocharacter lattice $\xcoch(T)$, where $T$ is a fixed maximal torus of $G$. The {\em extended affine Weyl group} $\tilW$ is the semidirect product
\begin{equation*}
\tilW=\xcoch(T)\rtimes W.
\end{equation*}
When $G$ is simply-connected, so that $\xcoch(T)$ is spanned by coroots, $\tilW$ is a Coxeter group with simple reflections $\{s_{0},s_{1}, \cdots,s_{r}\}$ in bijection with the nodes of the extended Dynkin diagram of $G$.  In general, $\tilW$ is a semidirect product of the {\em affine Weyl group} $\Wa=(\ZZ R^{\vee})\rtimes W$ (where $\ZZ R^{\vee}$ is the coroot lattice), which is a Coxeter group, and an abelian group $\Om\cong \xcoch(T)/\ZZ R^{\vee}$. The group $\tilW$ naturally acts on the affine space $\xcoch(T)_{\RR}$ by affine transformations, where $\xcoch(T)$ acts by translations.

\begin{theorem}[Lusztig {\cite{L96}}, Sage {\cite{Sage}}]\label{th:aff action} There is a canonical action of $\tilW$  on $\homog{*}{\sY_{\gamma}}$. 
\end{theorem}

Since $\sY_{\gamma}$ is not of finite type, the $\ell$-adic homology $\homog{*}{\sY_{\gamma}}$ is understood as the inductive limit $\varinjlim_{n}\homog{*}{\sY_{\gamma,n}}$, whenever we present $\sY_{\gamma}$ as a union of projective subschemes $\sY_{\gamma,n}$.

\subsubsection{Sketch of the construction of the $\tilW$-action} We consider only the case $G$ is simply-connected so that $\tilW=\Wa$ is generated by affine simple reflections $s_{0}, \cdots, s_{r}$. For each parahoric  subgroup $\bP\subset LG$ we have a corresponding affine Springer fiber $\tX_{\bP,\gamma}$. For $\bP$ containing a fixed Iwahori subgroup $\bI$, we have a projection $\pi_{\bP,\gamma}:\tY_{\gamma}\to \tX_{\bP,\gamma}$. 

Let $L_{\bP}$ be  the Levi quotient of $\bP$ and $\frl_{\bP}=\Lie L_{\bP}$. We have an evaluation map $\ev_{\bP,\gamma}: \tX_{\bP,\gamma}\to [\frl_{\bP}/L_{\bP}]$ defined as follows. For $[g]\in\Fl_{\bP}$ such that $\Ad(g^{-1})\gamma\in\Lie \bP$, we send the coset $[g]=g\bP$ to the image of $\Ad(g^{-1})\gamma$ under the projection $\Lie \bP\to\frl_{\bP}$. This is well-defined up to the adjoint action of $L_{\bP}$. We have a Cartesian diagram
\begin{equation*}
\xymatrix{\tY_{\gamma}\ar[d]^{\pi_{\bP,\gamma}}\ar[r]^{\ev_{\bI,\gamma}} & [\wt{\frl}_{\bP}/L_{\bP}]\ar[d]^{\pi_{\frl_{\bP}}}\\
\tX_{\bP,\gamma}\ar[r]^{\ev_{\bP,\gamma}} & [\frl_{\bP}/L_{\bP}]}
\end{equation*}
where $\pi_{\frl_{\bP}}$ is the Grothendieck alteration for the reductive group $\frl_{\bP}$. By the Springer theory for $L_{\bP}$, we have a $W(L_{\bP})$-action on the direct image complex $\bR\pi_{\frl_{\bP,*}}\DD$ (where $\DD$ stands for the dualizing complex for $[\wt{\frl}_{\bP}/L_{\bP}]$). By proper base change, we get an action of $W(L_{\bP})$ on $\bR\pi_{\bP,\gamma,*}\DD_{\tY_{\gamma}}$, and hence on $\homog{*}{\sY_{\gamma}}$. Taking a standard parahoric $\bP$ corresponding to the $i$-th node in the extended Dynkin diagram, then $W(L_{\bP})=\jiao{s_{i}}$, and we get an involution $s_{i}$ acting on $\homog{*}{\sY_{\gamma}}$. To check the braid relation between  $s_{i}$ and $s_{j}$ for neighboring nodes $i$ and $j$, we may choose a standard parahoric $\bP$ such that $W(L_{\bP})=\jiao{s_{i},s_{j}}$ and the braid relation holds because $W(L_{\bP})$ acts on $\homog{*}{\sY_{\gamma}}$. This shows that $\tilW$ acts on $\homog{*}{\sY_{\gamma}}$. \qed

Despite the simplicity of the construction of the $\tilW$-action on the homology of affine Springer fibers, the calculation of these actions are quite difficult. One new feature here is that the action of $\tilW$ on $\homog{*}{\sY_{\gamma}}$ {\em may not be semisimple}, as we shall see in the next example.

\subsubsection{An example in $\SL_{2}$}\label{sss:tilW SL2 ram} Consider the case $G=\SL_{2}$ and the element $\gamma=\mat{0}{t^{2}}{t}{0}$. This is a subregular case. The affine Springer fiber $\sY_{\gamma}$ has two irreducible components $C_{0}$ and $C_{1}$ both isomorphic to $\PP^{1}$. Here $C_{0}$ parametrizes chains of lattices $\{\L_{i}\}$ where $\L_{0}$ is the standard lattice $\calO^{2}_{F}$ and $\L_{-1}\subset \L_{0}$ is varying. The other component $C_{1}$ parametrizes chains of lattices $\{\L_{i}\}$ where $\L_{-1}=t\calO_{F}\oplus\calO_{F}$ and $\L_{0}$ is varying. The fundamental classes $[C_{0}], [C_{1}]$ give a natural basis for the top homology group $\homog{2}{\sY_{\gamma}}$. One can show that the action of the affine Weyl group $\tilW=\jiao{s_{0},s_{1}}$ on $\homog{2}{\sY_{\gamma}}$ takes the following form under the basis $[C_{0}]$ and $[C_{1}]$:
\begin{equation*}
s_{0}=\mat{-1}{2}{0}{1}; \quad s_{1}=\mat{1}{0}{2}{-1}.
\end{equation*} 
We see that $\homog{2}{\sY_{\gamma}}$ is a {\em nontrivial extension} of the sign representation of $\tilW$ by the trivial representation spanned by $[C_{0}]+[C_{1}]$. One can also canonically identify the $\tilW$-module $\homog{2}{\sY_{\gamma},\ZZ}$ with the affine coroot lattice of the loop group $LG$.

\subsection{Comments and generalizations}

\subsubsection{Relation with orbital integrals} As we will see in \S\ref{s:orb}, the cohomology and point-counting of affine Springer fibers are closely related to orbital integrals on $p$-adic groups $G(\FF_{q}\lr{t})$. 

\subsubsection{Extended symmetry} The $\tilW$-action on $\homog{*}{\sY_{\gamma}}$ can be extended to an action of the wreath product $\Sym(\xch(T))\rtimes \tilW$. For homogeneous affine Springer fibers $\sY_{\gamma}$ (those admitting a torus action coming from loop rotation), the equivariant cohomology group $\upH^{*}_{\Gm}(\sY_{\gamma})$ admits an action of the {\em graded double affine Hecke algebra},  which is a deformation of $\Sym(\xch(T))\rtimes \tilW$. For details we refer to \cite{OY}. Vasserot and Varagnolo \cite{Vass} \cite{VV} constructed an action of the {\em double affine Hecke algebra} on the $K$-groups of affine Springer fibers.

\subsubsection{The group version} Taking $\gamma\in G(F)$ instead of in $\frg(F)$, one can similarly define the group version of  affine Springer fibers, which we still denote by $\tX_{\gamma}$ with reduced structure $\sX_{\gamma}$. For a $k$-algebras $R$, we have
\begin{equation}\label{gp AFS}
\tX_{\gamma}(R)=\{[g]\in \Gr_{G}(R)|g^{-1}\gamma g\in L^{+}G(R)\}.
\end{equation}
However, in the group version, the definition above admits an interesting generalization. Recall the $L^{+}G$-double cosets in $LG$ are indexed by dominant cocharacters $\l\in\xcoch(T)^{+}$. For $\l\in\xcoch(T)^{+}$ we denote the corresponding double coset by $(LG)_{\l}$, which is the preimage of the Schubert stratum $\Gr_{G,\l}$ under the projection $LG\to \Gr_{G}$. Similarly we may define $(LG)_{\leq \l}$ to be the preimage of the closure $\Gr_{G,\leq \l}$ of $\Gr_{G,\l}$. One can replace the condition $g^{-1}\gamma g\in L^{+}G(R)$ in \eqref{gp AFS} by $g^{-1}\gamma g\in (LG)_{\l}(R)$ or $g^{-1}\gamma g\in (LG)_{\leq\l}(R)$, and take reduced structures to obtain reduced generalized affine Springer fibers $\sX_{\l,\gamma}$ and $\sX_{\leq\l, \gamma}$. We have an open embedding $\sX_{\l,\gamma}\incl\sX_{\leq\l,\gamma}$, whose complement is the union of $\sX_{\mu,\gamma}$ for dominant coweights $\mu\leq\l$. The motivation for introducing $\sX_{\leq\l,\gamma}$ is to give geometric interpretation of orbital integrals of spherical Hecke functions on $G(F)$.

A.Bouthier has established the fundamental geometric properties of $\sX_{\l,\gamma}$, parallel to Theorem \ref{th:ft quot} and Theorem \ref{th:aff dim}. 

\begin{theorem}[Bouthier \cite{Bouthier}]  Let $\gamma\in G(F)$ be regular semisimple, and let $\l\in\xcoch(T)^{+}$.
\begin{enumerate}
\item The generalized affine Springer fiber $\sX_{\l,\gamma}$ is non-empty if and only if $[\nu_{\gamma}]\leq\l$, where $[\nu_{\gamma}]\in\xcoch(T)^{+}_{\QQ}$ is the Newton point of $\gamma$, see \cite[\S2]{KV}.
\item The ind-scheme $\sX_{\l,\gamma}$ is locally of finite type.
\item We have
\begin{equation*}
\dim\sX_{\l,\gamma}=\jiao{\rho,\l}+\frac{1}{2}(\val_{F}\Delta(\gamma)-c(\gamma))
\end{equation*}
where $\rho\in\xch(T)_{\QQ}$ is half the sum of positive roots, and $\Delta(\gamma)$ and $c(\gamma)$ are defined similarly as in the Lie algebra situation.
\end{enumerate}
\end{theorem}
The proof of this theorem uses the theory of Vinberg semigroups, which is a kind of compactification of $G$.

\subsubsection{} In \cite{KV}, Kottwitz and Viehmann defined two generalizations of affine Springer fibers for elements  $\gamma$ in the Lie algebra $\frg(F)$.

\subsubsection{} As an analog of Hessenberg varieties, one can also consider the following situation. Let $(\rho,V)$ be a linear representation of a reductive group $G$ over $k$. Let $\L\subset V\otimes F$ be an $\calO_{F}$-lattice stable under $G(\calO_{F})$. For $v\in V\otimes F$ we may define a sub-ind-scheme $\tX_{v}$ of $\Gr_{G}$ 
\begin{equation*}
\tX_{\L, v}(R)=\{[g]\in\Gr_{G}(R)|\rho(g^{-1})v\in\L\widehat{\otimes}_{k}R\}.
\end{equation*}
Let $\sX_{\L,v}$ be the reduced structure of $\tX_{\L,v}$. The cohomology of these ind-schemes are related to orbital integrals that appear in relative trace formulae.

\subsection{Exercises}

\subsubsection{}\label{exII:SO parahoric} Let $G=\SO(V,q)$ for some vector space $V$ over $k$ equipped with a quadratic form $q$. Give an interpretation of the parahoric subgroups and affine partial flag varieties of $LG$ in terms of self-dual lattice chains in $V\otimes_{k}F$, in the same style as in \S\ref{sss:Sp parahoric}.

\subsubsection{}\label{exII:SL2} Verify the descriptions of the affine Springer fibers for $G=\SL_{2}$ given in \S\ref{P1 chain} and \S\ref{SL2 ram dim 1}.

\subsubsection{} Let $G=\SL_{2}$ and $\gamma=\mat{0}{t^{n}}{1}{0}\in\frg(F)$. Describe $\sY_{\gamma}$.

\subsubsection{} Let $G=\SL_{2}$ and $\gamma=\mat{0}{t^{n}}{1}{0}\in\frg(F)$. Construct a nontrivial $\Gm$-action on $\sX_{\gamma}$ involving loop rotations (i.e., the action scales $t$) and determine its fixed points.

\subsubsection{} Let $G=\SL_{2}$ and $\gamma=\mat{t^{n}}{0}{0}{-t^{n}}\in\frg(F)$. Let $T\subset G$ be the diagonal torus, then $G_{\gamma}=T\otimes_{k}F$. What is the regular locus $\sX^{\reg}_{\gamma}$ (see \S\ref{dim formula discussion})? Study the $L^+T$-orbits on $\sX_{\gamma}$. 

\subsubsection{} In the setup of \S\ref{sss:L gamma}, show that the action of $\L_{\gamma}$ on $\Gr_{G}$ is free, which implies that its action on $\sX_{\gamma}$ is free. Show also that the permutation action of $\L_{\gamma}$ on the set of irreducible components of $\sX_{\gamma}$ is free.

\subsubsection{} For $G=\GL_{n}$, let $L\subset G$ be the Levi subgroup consisting of block diagonal matrices with sizes of blocks $n_{1},\cdots, n_{s}$, $\sum_{i}n_{i}=n$. Let $\gamma=(\gamma_{1},\cdots,\gamma_{s})\in\frl(F)$ be regular semisimple as an element in $\frg(F)$. What is the invariant $\Delta^{G}_{L}(\gamma)$ (see \eqref{Delta G L}) in terms of familiar invariants of the characteristic polynomials of the $\gamma_{i}$?

\subsubsection{}\label{exII:SL3} Let $G=\SL_{3}$ and $\gamma=\diag(x_{1}t,x_{2}t,x_{3}t)\in\frg(F)$, with $x_{i}\in k$ pairwise distinct and $x_{1}+x_{2}+x_{3}=0$. Describe the affine Springer fibers $\sX_{\gamma}$ and $\sY_{\gamma}$.

Note: this is a good exercise if you have a whole day to kill.

\subsubsection{} For $G=\Sp_{6}$ and $\gamma=\mat{0}{tX}{Y}{0}$ as in Example \ref{ex:Sp}, describe the $\Gm$-fixed points on $\sX_{\gamma}$ and $\sY_{\gamma}$.

\subsubsection{} Verify the calculations in \S\ref{sss:tilW SL2 ram}.

\subsubsection{} Let $G=\SL_{2}$ and let $\gamma=\mat{t}{0}{0}{-t}$. Describe the affine Springer fiber $\sY_{\gamma}$. What is the action of $\tilW=\jiao{s_{0},s_{1}}$ on $\homog{2}{\sY_{\gamma}}$ in terms of the basis given by the irreducible components of $\sY_{\gamma}$?


\section{Lecture III: Orbital integrals}\label{s:orb}
The significance of affine Springer fibers in representation theory is demonstrated by their close relationship with orbital integrals. Orbital integrals are certain integrals that appear in the harmonic analysis of $p$-adic groups. Just as conjugacy classes of a finite group are fundamental to understanding its representations, orbital integrals are fundamental to understanding representations of $p$-adic groups. In certain cases, orbital integrals can be interpreted as counting points on affine Springer fibers.

\subsection{Integration on a $p$-adic group}\label{ss:int}

\subsubsection{The setup} Let $F$ be a local non-archimedean field, i.e., $F$ is either a finite extension of $\QQ_{p}$ or a finite extension of $\FF_{p}\lr{t}$. Then $F$ has a discrete valuation $\val: F^{\times}\to \ZZ$ which we normalize to be surjective. Let $\calO_{F}$ be the valuation ring of $F$ and $k$ be the residue field. Therefore, unlike in the previous sections,  $k$ is a finite field. We assume that $\chark$ is large with respect to the groups in question.

\subsubsection{Haar measure and integration} Let $G$ be an algebraic group over $F$. The topological group $G(F)$ is locally compact and totally disconnected. It has a right invariant Haar measure $\mu_{G}$ which is unique up to a scalar. For a measurable subset $S\subset G(F)$, we denote its volume under $\mu_{G}$ by $\vol(S,\mu_{G})$. Fixing a compact open subgroup $K_{0}\subset G(F)$, we may normalize the Haar measure $\mu_{G}$ so that $K_{0}$ has volume $1$. For example, if we choose an integral model $\calG$ of $G$ over $\calO_{F}$, we may take $K_{0}=\calG(\calO_{F})$.  

With the Haar measure $\mu_{G}$ one can integrate smooth (i.e., locally constant) compactly supported functions on $G(F)$ with complex values. We denote this function space by $\sS(G(F))$ (where $\sS$ stands for Schwarz). For $f\in\sS(G(F))$, the integral
\begin{equation*}
\int_{G(F)}f\mu_{G}
\end{equation*} 
can be calculated as follows. One can find a subgroup $K\subset K_{0}$ of finite index such that $f$ is right $K$-invariant, i.e., $f(gx)=f(g)$ for all $g\in G(F)$ and $x\in K$ (see Exercise \ref{exIII:K inv}). Then the integral above becomes a weighted counting in the coset $G(F)/K$:
\begin{equation*}
\int_{G(F)}f\mu_{G}=\vol(K,\mu_{G})\sum_{[g]\in G(F)/K}f(g)=\frac{1}{[K_{0}:K]}\sum_{[g]\in G(F)/K}f(g).
\end{equation*}
It is easy to check that the right side above  is independent of the choice of $K$ as long as $f$ is right $K$-invariant and $K$ has finite index in $K_{0}$.

\subsubsection{Variant} Let $H\subset G$ be a subgroup defined over $F$ together with a Haar measure $\mu_{H}$ on it. Consider a function $f\in\sS(H(F)\bs G(F))$, i.e., $f$ is a left $H(F)$-invariant, locally constant function on $G(F)$ whose support is compact modulo $H(F)$, we may define the integral
\begin{equation}\label{int quot}
\int_{H(F)\bs G(F)}f\frac{\mu_{G}}{\mu_{H}}.
\end{equation} 
This integral is calculated in the following way. Again we choose a finite index subgroup $K\subset G(F)$ such that $f$ is right $K$-invariant. Then the integral \eqref{int quot} can be written as a weighted sum over double cosets $H(F)\bs G(F)/K$:
\begin{equation*}
\int_{H(F)\bs G(F)}f\frac{\mu_{G}}{\mu_{H}}=\frac{1}{[K_{0}:K]}\sum_{[g]\in H(F)\bs G(F)/K}\frac{f(g)}{\vol(H(F)\cap gKg^{-1}, \mu_{H})}.
\end{equation*}

\subsection{Orbital integrals} 

\subsubsection{Definition of orbital integrals} We continue with the setup of \S\ref{ss:int}. We denote the Lie algebra of $G$ by $\frg(F)$ to emphasize that it is a vector space over $F$. Let $\varphi\in\sS(\frg(F))$ and $\gamma\in\frg(F)$. Consider the map $G(F)\to \frg(F)$ given by $g\mapsto \Ad(g^{-1})\gamma$. Then the composition $f: g\mapsto \varphi(\Ad(g^{-1})\gamma)$ is a smooth function on $G(F)$.  Then $f$ is locally constant, left invariant under the centralizer $G_{\gamma}(F)$ of $\gamma$ in $G(F)$, and has compact support modulo $G_{\gamma}(F)$. 

Fix Haar measures $\mu_{G}$ on $G(F)$ and $\mu_{G_{\gamma}}$ on $G_{\gamma}(F)$. The following integral is then a special case of \eqref{int quot} (except that we write the integration variable $g\in G_{\gamma}(F)\bs G(F)$ explicit below while not in \eqref{int quot})
\begin{equation*}
O_{\gamma}(\varphi):=\int_{G_{\gamma}(F)\bs G(F)}\varphi(\Ad(g^{-1})\gamma)\frac{\mu_{G}}{\mu_{G_{\gamma}}}.
\end{equation*}
Such integrals are called {\em orbital integrals} on the Lie algebra $\frg(F)$. We may similarly define orbital integrals on the group $G(F)$ by replacing $\gamma$ with an element in $G(F)$ and $\varphi$ with an element in $\sS(G(F))$.

\subsubsection{Specific situation} For the rest of the section we will restrict to the following situation. Let $G$ be a split reductive group over $F$. We may fix an integral model of $G$ by base changing the corresponding Chevalley group scheme from $\ZZ$ to $\calO_{F}$. In the following we will regard $G$ as a reductive group scheme over $\calO_{F}$. We normalize the Haar measure $\mu_{G}$ on $G(F)$ by requiring that $K_{0}=G(\calO_{F})$ have volume $1$.

The Lie algebra $\frg(F)$ contains a canonical lattice $\frg(\calO_{F})$ coming from the integral model over $\calO_{F}$. We will be most interested in the orbital integral of the characteristic function $\varphi=\one_{\frg(\calO_{F})}$ of the lattice $\frg(\calO_{F})$. 

\subsubsection{The centralizer of $\gamma$} Suppose $\gamma$ is regular semisimple so that its centralizer $G_{\gamma}$ is a torus over $F$. Let $F^{\ur}=F\wh{\otimes}_{k}\kbar$, which is a complete discrete valuation field whose residue field is algebraically closed. We continue to let $t$ denote a uniformizer of $F$, which is also a uniformizer in $F^{\ur}$. Using $t$, the construction in \S\ref{sss:L gamma} gives an embedding $\xcoch(G_{\gamma}\otimes_{F}F^{\ur})\incl  G_{\gamma}(F^{\ur})$ whose image we still denote by $\L_{\gamma}$. This embedding being $\Gk$-equivariant, $\L_{\gamma}$ carries an action of $\Gk$. We think of $\L_{\gamma}$ as an \'etale group scheme over $k$, then the notation $\L_{\gamma}(k)$ makes sense, and it is just the $\Gk$-invariants in $\L_{\gamma}$ if we regard the latter as a plain group. Then $\L_{\gamma}(k)\subset G_{\gamma}(F)$ is a discrete and cocompact subgroup.

\subsubsection{Centralizers in $\GL_{n}$}\label{sss:GLn orb} Let $G=\GL_{n}$ and let $\gamma$ be a regular semisimple element in $\frg(F)$ which is not necessarily diagonalizable over $F$. Assume either $\textup{char}(F)=0$ or $\chark>n$. As in \S\ref{sss:GLn centralizer}, the characteristic polynomial $P(x)=x^{n}+a_{1}x^{n-1}+\cdots+a_{n}\in F[x]$ of $\gamma$ is separable, hence the $F$-algebra $F[x]/(P(x))$ is isomorphic to a product of fields $F_{1}\times \cdots \times F_{m}$. We have
\begin{equation*}
G_{\gamma}(F)\cong F^{\times}_{1}\times\cdots\times F^{\times}_{m}.
\end{equation*}
In this case, the lattice $\L_{\gamma}(k)\subset G_{\gamma}(F)=F^{\times}_{1}\times\cdots\times F^{\times}_{m}$ consists of elements of the form $(t^{d_{1}},\cdots, t^{d_{m}})$ for $(d_{1},\cdots, d_{m})\in\ZZ^{m}$. The quotient $\L_{\gamma}(k)\bs G_{\gamma}(F)$ is isomorphic to $\prod_{i=1}^{m}F^{\times}_{i}/t^{\ZZ}$. Each $F^{\times}_{i}/t^{\ZZ}$ fits into an exact sequence
\begin{equation}\label{Fi seq}
0\to \calO^{\times}_{F_{i}}\to F^{\times}_{i}/t^{\ZZ}\to \ZZ/e_{i}\ZZ\to 0
\end{equation}
where $e_{i}$ is the ramification degree of the extension $F_{i}/F$, therefore the quotient $\L_{\gamma}(k)\bs G_{\gamma}(F)$ is compact.

\subsubsection{Orbital integrals in terms of counting}\label{sss:X gamma} Consider the following subset of $G(F)/G(\calO_{F})$
\begin{equation*}
X_{\gamma}:=\{[g]\in G(F)/G(\calO_{F})|\Ad(g^{-1})\gamma\in\frg(\calO_{F})\}.
\end{equation*}
This is a set-theoretic version of the affine Springer fiber. 

The group $G_{\gamma}(F)$ acts on $X_{\gamma}$ by the rule $G_{\gamma}(F)\ni h:[g]\mapsto[hg]$. For any free abelian group $L\subset G_{\gamma}(F)$, its action on $G(F)/G(\calO_{F})$ by left translation is free (because the stabilizers are necessarily finite), hence it acts freely on $X_{\gamma}$.

More generally, for any discrete cocompact subgroup $L\subset G_{\gamma}(F)$, the quotient groupoid $L\bs X_{\gamma}$ is {\em finitary}, i.e.,  it has finitely many isomorphism classes and the automorphism group of each object is finite. For a finitary groupoid $Y$, we define the cardinality of $Y$ to be 
\begin{equation}\label{card groupoid}
\#Y:=\sum_{y\in \textup{Ob}(Y)/\cong}\frac{1}{\#\Aut(y)}
\end{equation}

The next lemma follows directly from the definitions, whose proof is left to the reader as Exercise \ref{exIII:Ogamma}.
\begin{lemma}\label{l:Ogamma} Let $\gamma$ be a regular semisimple element in $\frg(F)$. Let $L\subset G_{\gamma}(F)$ be any discrete cocompact subgroup. We have
\begin{equation*}
O_{\gamma}(\one_{\frg(\calO_{F})})=\frac{1}{\vol(G_{\gamma}(F)/L,\mu_{G_{\gamma}})}\#\left(L\bs X_{\gamma}\right)
\end{equation*}
with the cardinality on the right side interpreted as in \eqref{card groupoid}.
\end{lemma}

\subsubsection{The case $G=\GL_{n}$ and fractional ideals}\label{sss:GLn orb formula} We continue with the situation in \S\ref{sss:GLn orb}. Under the identification of $G(F)/G(\calO_{F})$ with the set of $\calO_{F}$-lattices in $F^{n}$ (see \eqref{Gr lattices}), we have
\begin{equation*}
X_{\gamma}\cong\{\textup{lattices }\L\subset F^{n}|\gamma \L\subset \L\}.
\end{equation*} 
The bijection sends $[g]\in X_{\gamma}$ to the lattice $\L=g\calO^{n}_{F}$. 

We give another interpretation of $X_{\gamma}$. Let $P(x)=x^{n}+a_{1}x^{n-1}+\cdots+a_{n}$ be the characteristic polynomial of $\gamma$. Let
\begin{equation*}
A=\calO_{F}[x]/(P(x))
\end{equation*}
be the commutative $\calO_{F}$-subalgebra of $\frg(F)$ generated by $\gamma$. The ring of total fractions $\Frac(A)$ is a finite \'etale $F$-algebra of degree $n$, and $A$ is an order in it. The canonical action of $A$ on $F^{n}$ realizes $F^{n}$ as a free $\Frac(A)$-module of rank $1$.  Recall a {\em fractional $A$-ideal} is a finitely generated $A$-submodule $M\subset \Frac(A)$. If we choose an element $e\in F^{n}$ as a basis for the $\Frac(A)$-module structure, we get a bijection
\begin{equation}\label{frac ideals}
\{\mbox{fractional $A$-ideals}\}\bij X_{\gamma}
\end{equation}
which sends $M\subset \Frac(A)$ to the $\calO_{F}$-lattice $M\cdot e\subset F^{n}$.

Using the algebra $A$, we have a canonical isomorphism
\begin{equation*}
G_{\gamma}(F)\cong \Frac(A)^{\times}.
\end{equation*}
This isomorphism intertwines the action of $G_{\gamma}(F)$ on $X_{\gamma}$ by left translation and the action of $\Frac(A)^{\times}$ on the set of fractional $A$-ideals by multiplication. 

When $A$ happens to be a product of Dedekind domains (i.e., $A$ is the maximal order in $\Frac(A)$), all fractional $A$-ideals are principal, which is the same as saying that the action of $G_{\gamma}(F)\cong \Frac(A)^{\times}$ on $X_{\gamma}$ is transitive. In general, principal fractional ideals form a homogeneous space $\Frac(A)^{\times}/A^{\times}$ under $G_{\gamma}(F)$; the difficulty in counting $X_{\gamma}$ in general is caused by the singularity of the ring $A$.

We normalize the Haar measure $\mu_{G_{\gamma}}$ on $G_{\gamma}(F)\cong\prod_{i}F^{\times}_{i}$ so that $\prod_{i}\calO^{\times}_{F_{i}}$ gets volume $1$. Choose a uniformizer $t_{i}\in F_{i}$, we may form the lattice $L_{0}=t^{\ZZ}_{1}\times \cdots\times t^{\ZZ}_{m}\subset G_{\gamma}(F)$. Now $G_{\gamma}(F)/L_{0}\cong \prod_{i}\calO_{F_{i}}^{\times}$ has volume $1$. Therefore Lemma \ref{l:Ogamma} gives
\begin{equation}\label{GLn L0}
O_{\gamma}(\one_{\frg(\calO_{F})})=\#(L_{0}\bs X_{\gamma}).
\end{equation}
Using \eqref{frac ideals}, we may interpret \eqref{GLn L0} as saying that $O_{\gamma}(\one_{\frg(\calO_{F})})$ is the number of fractional $A$-ideals up to multiplication by the powers of the $t_{i}$'s.

\subsection{Relation with affine Springer fibers}\label{ss:orb ASF}
From this subsection we restrict to the case $F$ is a local function field, i.e., $F=k\lr{t}$ for a finite field $k=\FF_{q}$. Let $\gamma\in\frg(F)$ be regular semisimple. The definitions of the affine Grassmannian and the affine Springer fiber $\sX_{\gamma}$ we gave in \S\ref{s:ASF} make sense when the base field $k$ is a finite field, so we have a sub-ind-scheme $\sX_{\gamma}$ of $\Gr_{G}$, both defined over $k$. 

The following lemma is clear from the definitions.
\begin{lemma}\label{l:X gamma k} The set of $k$-rational points $\sX_{\gamma}(k)$ is the same as the set $X_{\gamma}$ defined in \S\ref{sss:X gamma}, both as subsets of $\Gr_{G}(k)=G(F)/G(\calO_{F})$.
\end{lemma}

\subsubsection{} If we base change from $k$ to $\kbar$,  by Theorem \ref{th:ft quot} we know that $\L_{\gamma}\bs \sX_{\gamma,\kbar}$ is a proper algebraic space over $\kbar$. The proof there actually shows that this algebraic space is defined over $k$, which we denote by $\L_{\gamma}\bs \sX_{\gamma}$. We emphasize here that $\L_{\gamma}$ is viewed as an \'etale group scheme over $k$ whose $\kbar$-points is the plain group used to be denoted $\L_{\gamma}$.

In view of Lemma \ref{l:Ogamma} and Lemma \ref{l:X gamma k}, it is natural to expect that the orbital integral $O_{\gamma}(\one_{\frg(\calO_{F})})$ can be expressed as the number of $k$-points on the quotient $\L_{\gamma}\bs \sX_{\gamma}$. Such a relationship takes its cleanest form when $G=\GL_{n}$.

\subsubsection{The case of $\GL_{n}$}\label{sss:GLn orb coho}
In the situation of \S\ref{sss:GLn orb formula}, upon choosing uniformizers $t_{i}\in F_{i}$, we defined the lattice $L_{0}=t^{\ZZ}_{1}\times \cdots\times t^{\ZZ}_{m}$. Base change from $k$ to $\kbar$, we may similarly define a lattice $\wt\L_{\gamma}$ as in \S\ref{sss:GLn centralizer}, using the same choice of uniformizers $t_{i}\in F_{i}$ (note that $F_{i}\wh\otimes_{k}\kbar$ may split into a product of fields, but $t_{i}\otimes 1$ will project to a uniformizer in each factor). The $\Gk$-action on $\wt\L_{\gamma}$ gives it the structure of an \'etale group scheme over $k$, just as $\L_{\gamma}$. We have $\wt\L_{\gamma}(k)=L_{0}$. There is an analog of Theorem \ref{th:ft quot} if we replace $\L_{\gamma}$ with $\wt\L_{\gamma}$. In particular, $\wt\L_{\gamma}\bs\sX_{\gamma}$ is a proper algebraic space over $k$ admitting a surjective map from a projective scheme.

\begin{theorem}\label{th:GLn orb} Let $G=\GL_{n}$. Let $\gamma$ be a regular semisimple element in $\frg(F)$. We fix the Haar measure on $G_{\gamma}(F)$ such that its maximal compact subgroup gets volume $1$. Then we have
\begin{equation}\label{GLn orb}
O_{\gamma}(\one_{\frg(\calO_{F})})=\#(\wt\L_{\gamma}\bs\sX_{\gamma})(k)=\sum_{i}(-1)^{i}\Tr\left(\Frob_{k}, \cohog{i}{\wt\L_{\gamma}\bs\sX_{\gamma,\kbar},\Ql}\right).
\end{equation}
\end{theorem}

The second equality in \eqref{GLn orb} follows from the Grothendieck-Lefschetz trace formula. Comparing \eqref{GLn orb} with \eqref{GLn L0}, we only need to argue that $L_{0}\bs X_{\gamma}=\wt\L_{\gamma}(k)\bs \sX_{\gamma}(k)$ is the same as $(\wt\L_{\gamma}\bs\sX_{\gamma})(k)$. Let us first make some general remarks about $k$-points on quotient stacks.

\subsubsection{$k$-points of a quotient}\label{sss:k pts quot} We consider  a quotient stack $\sY=[A\bs X]$ where $X$ is a scheme over $k$ and $A$ is an algebraic group over $k$ acting on $X$. Then, by definition, $\sY(k)$ is the groupoid of pairs $(S, f)$ where $S\to \Spec k$
is a left $A$-torsor, and $f:S\to X$ is an $A$-equivariant morphism. The isomorphism class of the $A$-torsor $S$ is classified by the Galois cohomology $\cohog{1}{k, A}:=\cohog{1}{\Gk,A(\kbar)}$. For each class $\xi\in \cohog{1}{k,A}$, let $S_{\xi}$ be an $A$-torsor over $k$ with class $\xi$. We may define a twisted form of $X$ over $k$ by $X_{\xi}:=A\bs (S_{\xi}\times X)$. We also have the inner form $A_{\xi}:=\Aut_{A}(S_{\xi})$ of $A$ acting on the $k$-scheme $X_{\xi}$. It is easy to see that $A$-equivariant morphisms $f:S_{\xi}\to X$ are in bijection with $X_{\xi}(k)$. Therefore we get a decomposition of groupoids
\begin{equation}\label{k pts quot}
\sY(k)=[A\bs X](k)\cong\bigsqcup_{\xi\in\cohog{1}{k,A}}A_{\xi}(k)\bs X_{\xi}(k).
\end{equation}

\subsubsection{Proof of Theorem \ref{th:GLn orb}} We have reduced to showing that $(\wt\L_{\gamma}\bs\sX_{\gamma})(k)= \wt\L_{\gamma}(k)\bs \sX_{\gamma}(k)$. By \eqref{k pts quot}, it suffices to show that $\cohog{1}{k,\wt\L_{\gamma}}=0$. We use the notation from \S\ref{sss:GLn orb coho}. From the definition of $\wt\L_{\gamma}$ we see that, as a $\Gk$-module, it is of the form
\begin{equation*}
\wt\L_{\gamma}\cong \bigoplus_{i=1}^{m}\Ind^{\Gk}_{\Gal(\kbar/k_{i})}\ZZ
\end{equation*}
where $k_{i}$ is the residue field of $F_{i}$. Therefore
\begin{equation*}
\cohog{1}{k,\wt\L_{\gamma}}\cong\bigoplus_{i=1}^{m}\cohog{1}{k_{i},\ZZ}=0.
\end{equation*}
\qed

\subsection{Stable orbital integrals} The setup is the same as \S\ref{ss:orb ASF}. For general $G$, the generalization of the formula \eqref{GLn orb} is not straightforward. Namely, the orbital integral $O_{\gamma}(\one_{\frg(\calO_{F})})$ by itself does not have a cohomological interpretation. The problem is that we may not be able to find an analog of $\wt\L_{\gamma}$ with vanishing first Galois cohomology so that $(\wt\L_{\gamma}\bs\sX_{\gamma})(k)$ is simply $\wt\L_{\gamma}(k)\bs \sX_{\gamma}(k)$. In view of formula \eqref{k pts quot},  a natural fix to this problem is to consider the twisted forms of $\sX_{\gamma}$ altogether. This suggests taking not just the orbital integral $O_{\gamma}(\one_{\frg(\calO_{F})})$ but a sum of several orbital integrals $O_{\gamma'}(\one_{\frg(\calO_{F})})$.

\subsubsection{Stable conjugacy} Fix a regular semisimple element $\gamma\in\frg(F)$. An element $\gamma'\in\frg(F)$ is called {\em stably conjugate} to $\gamma$ if it is in the same $G(\ov{F})$-orbit of $\gamma$. Equivalently, $\gamma'$ is stably conjugate to $\gamma$ if $\chi(\gamma')=\chi(\gamma)\in\frc(F)$. For $\gamma'$ stably conjugate to $\gamma$, one can attach a Galois cohomology class $\inv(\gamma,\gamma')\in\cohog{1}{F,G_{\gamma}}$ which becomes trivial in $\cohog{1}{F,G}$. The assignment $\gamma'\mapsto\inv(\gamma,\gamma')$ gives a bijection of pointed sets
\begin{equation}\label{stable class}
\{\gamma'\in\frg(F) \mbox{ stably conjugate to } \gamma\}/G(F) \cong \ker(\cohog{1}{F,G_{\gamma}}\to \cohog{1}{F,G}).
\end{equation}

\subsubsection{The case of $\GL_{n}$} When $G=\GL_{n}$, we use the notation from \S\ref{sss:GLn centralizer}. We have
\begin{equation*}
\cohog{1}{F,G_{\gamma}}\cong \prod_{i=1}^{m}\cohog{1}{F_{i},\Gm}=\jiao{1}.
\end{equation*}
Therefore, by \eqref{stable class}, all elements stably conjugate to $\gamma$ are in fact $G(F)$-conjugate to $\gamma$. Of course this statement can also be proved directly using companion matrices. 

\subsubsection{The case of $\SL_{2}$}\label{sss:st conj SL2} We consider the case $G=\SL_{2}$. Let  $a\in k^{\times}-(k^{\times})^{2}$. Let
\begin{equation*}
\gamma=\mat{0}{at}{t}{0}, \quad \gamma'=\mat{0}{at^{2}}{1}{0}
\end{equation*}
be two regular semisimple elements in $\sl_{2}(F)$. Since they have the same determinant $-at^{2}$, they are stably conjugate to each other. However, they are not conjugate to each other under $\SL_{2}(F)$. One can show that the stable conjugacy class of $\gamma$ consists of exactly two $\SL_{2}(F)$-orbits represented by $\gamma$ and $\gamma'$, see Exercise \ref{exIII:stable conj}.

\subsubsection{Definition of stable orbital integrals} Let $\varphi\in\sS(\frg(F))$. We define the {\em stable orbital integral} of $\varphi$ with respect to $\gamma$ to be
\begin{equation}\label{SO}
SO_{\gamma}(\varphi)=\sum_{\gamma'}O_{\gamma'}(\varphi)
\end{equation}
where $\gamma'$ runs over $G(F)$-orbits of elements $\gamma'\in \frg(F)$ that are stably conjugate to $\gamma$. For $\gamma'$ stably conjugate to $\gamma$, we have a {\em canonical isomorphism} $G_{\gamma}\cong G_{\gamma'}$ as $F$-groups. Therefore, once we fix a Haar measure on $G_{\gamma}(F)$, we get a canonical Haar measure on the other centralizers $G_{\gamma'}(F)$. It is with this choice that we define $O_{\gamma'}(\varphi)$ in  \eqref{SO}.

\subsubsection{Stable part of the cohomology}\label{sss:st H} The quotient group scheme $LG_{\gamma}/\L_{\gamma}$ acts on $\L_{\gamma}\bs\sX_{\gamma}$.  The component group $\pi_{0}(LG_{\gamma}/\L_{\gamma})$ is an \'etale group scheme over $k$ whose $\kbar$-points acts on $\cohog{*}{L_{\gamma}\bs\sX_{\gamma,\kbar},\Ql}$. We define the {\em stable part} of this cohomology group to be the invariants under this action
\begin{equation*}
\cohog{*}{\L_{\gamma}\bs\sX_{\gamma,\kbar}}_{\st}=\cohog{*}{\L_{\gamma}\bs\sX_{\gamma,\kbar}}^{\pi_{0}(LG_{\gamma}/\L_{\gamma})(\kbar)}.
\end{equation*}
It turns out that if we replace $\L_{\gamma}$ with any other $\Gk$-stable free abelian subgroup $\L\subset G_{\gamma}(F^{\ur})$ commensurable with $\L_{\gamma}$, the similarly defined stable part cohomology is canonically isomorphic to the above one.

\begin{theorem}[Special case of Goresky-Kottwitz-MacPherson {\cite[Th 15.8]{GKM}} and Ng\^{o} {\cite[Cor 8.2.10]{NgoFL}}]\label{th:orb coho} We have
\begin{equation}\label{d}
SO_{\gamma}(\one_{\frg(\calO_{F})})=\frac{1}{\vol(K_{\gamma},\mu_{G_{\gamma}})}\sum_{i}(-1)^{i}\Tr\left(\Frob_{k}, \cohog{i}{\L_{\gamma}\bs\sX_{\gamma,\kbar},\Ql}_{\st}\right).
\end{equation}
Here $K_{\gamma}\subset G_{\gamma}(F)$ is the parahoric subgroup of the torus $G_{\gamma}$. 
\end{theorem}

This cohomological interpretation of the stable orbital integral is the starting point of the proof of the Fundamental Lemma (see \cite{GKM} and \cite{NgoFL}).

\subsubsection{}\label{sss:K gamma} Let us briefly comment on the definition of the parahoric subgroup $K_{\gamma}$. Let $\tilK_{\gamma}\subset G_{\gamma}(F)$ be the maximal compact subgroup. Then there is a canonical smooth group scheme $\calG_{\gamma}$ over $\calO_{F}$ whose $F$-fiber is $G_{\gamma}$ and whose $\calO_{F}$ points is $\tilK_{\gamma}$. This group scheme $\calG_{\gamma}$ is the finite type N\'eron model for the torus $G_{\gamma}$. Now let $\calG^{\circ}_{\gamma}\subset\calG_{\gamma}$ be the open subgroup scheme obtained by removing the non-neutral component of the special fiber of $\calG_{\gamma}$. Then, by definition, $K_{\gamma}=\calG^{\circ}_{\gamma}(\calO_{F})$. The positive loop group $L^{+}\calG^{\circ}_{\gamma}$ is, up to nilpotents, the neutral component of $LG_{\gamma}$.

\subsubsection{Sketch of proof of Theorem \ref{th:orb coho}}\label{sss:pf SO} We sketch an argument which is closer to that of \cite{GKM} than that of \cite{NgoFL}. One can show that there exists an \'etale $k$-subgroup $\wt\L\subset LG_{\gamma}$ commensurable with $\L_{\gamma}$, which maps {\em onto} the \'etale group scheme $\pi_{0}(LG_{\gamma})$ over $k$ (and the kernel is necessarily finite). The group $\wt\L$ is an analog of $\wt\L_{\gamma}$ defined in \S\ref{sss:GLn orb coho}. We form the quotient stack $[\wt\L\bs \sX_{\gamma}]$. Then we have an isomorphism
\begin{equation}\label{st coho quot}
\cohog{*}{\L_{\gamma}\bs\sX_{\gamma,\kbar},\Ql}_{\st}\cong \cohog{*}{[\wt\L\bs\sX_{\gamma,\kbar}],\Ql}.
\end{equation}
In fact, by the discussion in the end of \S\ref{sss:st H}, the left side above does not change if we replace $\L_{\gamma}$ by a commensurable lattice, so by shrinking $\L_{\gamma}$ we may assume $\L_{\gamma}\subset \wt\L$. On the other hand, the right side above can be computed by the Leray spectral sequence associated with the map $\L_{\gamma}\bs\sX_{\gamma,\kbar}\to[\wt\L\bs\sX_{\gamma,\kbar}]$ which is a torsor under the finite discrete group $(\wt\L/\L_{\gamma})(\kbar)$, therefore $\cohog{*}{[\wt\L\bs\sX_{\gamma,\kbar}],\Ql}=\cohog{*}{\L_{\gamma}\bs\sX_{\gamma,\kbar},\Ql}^{(\wt\L/\L_{\gamma})(\kbar)}$. Since $(\wt\L/\L_{\gamma})(\kbar)$ surjects onto $\pi_{0}(LG_{\gamma}/\L_{\gamma})(\kbar)$ by the choice of $\wt\L$, we see that
\begin{equation*}
\cohog{*}{\L_{\gamma}\bs\sX_{\gamma,\kbar},\Ql}^{\pi_{0}(LG_{\gamma}/\L_{\gamma})(\kbar)}=\cohog{*}{\L_{\gamma}\bs\sX_{\gamma,\kbar},\Ql}^{(\wt\L/\L_{\gamma})(\kbar)},
\end{equation*}
from which \eqref{st coho quot} follows.

By the discussion in \S\ref{sss:k pts quot} and \eqref{k pts quot}, we have
\begin{equation}\label{L bs X decomp}
[\wt\L\bs\sX_{\gamma,\kbar}](k)\cong\bigsqcup_{\xi\in\cohog{1}{k,\wt\L}} \wt\L(k)\bs\sX_{\gamma,\xi}(k).
\end{equation}
The inclusion $\wt\L\subset LG_{\gamma}$ gives
\begin{equation*}
\theta: \cohog{1}{k,\wt\L}\surj\cohog{1}{k,LG_{\gamma}}\cong\cohog{1}{F, G_{\gamma}}.
\end{equation*}
The first surjection follows from Lang's theorem because the quotient $LG_{\gamma}/\wt\L$ is connected; the second follows from another theorem of Lang which says that $\cohog{1}{F^{\ur}, G_{\gamma}}$ vanishes
\footnote{See \cite[Ch.X, \S7, p.170, Application and Example (b)]{Serre}. Let $K$ be a complete discrete valuation field with perfect residue field, and $K^{\ur}$ its maximal unramified extension. Then Lang's theorem asserts that $K^{\ur}$ is a $C_{1}$-field. Therefore $\cohog{1}{K^{\ur}, A}=0$ for any torus $A$ over $K$.}. For each $\xi\in\cohog{1}{k,\wt\L}$ such that $\sX_{\gamma,\xi}(k)\neq\varnothing$, one can show that the image of $\theta(\xi)$ in $\cohog{1}{F,G}$ is trivial. Therefore, by \eqref{stable class}, to each $\xi\in\cohog{1}{k,\wt\L}$ we can attach an element $\gamma'$ stably conjugate to $\gamma$, unique up to $G(F)$-conjugacy, such that $\inv(\gamma,\gamma')=\theta(\xi)$. One can show that $\sX_{\gamma,\xi}(k)$ is in bijection with the set $X_{\gamma'}$. Therefore,  \eqref{L bs X decomp} implies
\begin{equation}\label{pre SO}
\#[\wt\L\bs\sX_{\gamma}](k)=\#\ker(\theta)\sum_{\gamma'}\#(\wt\L(k)\bs X_{\gamma'})
\end{equation}
where the sum is over the $G(F)$-orbits of those $\gamma'$ stably conjugate to $\gamma$. Applying Lemma \ref{l:Ogamma} to the discrete cocompact subgroup $\wt\L(k)$, we have
\begin{equation*}
O_{\gamma'}(\one_{\frg(\calO_{F})})=\frac{1}{\vol(G_{\gamma}(F)/\wt\L(k),\mu_{G_{\gamma}})}\#(\wt\L(k)\bs X_{\gamma'}).
\end{equation*}
Plugging this into the right side of \eqref{pre SO}, we get
\begin{equation}\label{next SO}
\SO_{\gamma}(\one_{\frg(\calO_{F})})=\frac{1}{\#\ker(\theta)\vol(G_{\gamma}(F)/\wt\L(k),\mu_{G_{\gamma}})}\#[\wt\L\bs\sX_{\gamma}](k).
\end{equation} 
By the Grothendieck-Lefschetz trace formula,   $\#[\wt\L\bs\sX_{\gamma,\kbar}](k)$ is equal to the alternating Frobenius trace on $\cohog{*}{[\wt\L\bs\sX_{\gamma,\kbar}],\Ql}$, which, by \eqref{st coho quot}, can be identified with $\cohog{*}{\L_{\gamma}\bs\sX_{\gamma,\kbar},\Ql}_{\st}$. Therefore the theorem follows from the identity \eqref{next SO} together with the volume identity
\begin{equation}\label{vol id}
\#\ker(\theta)\vol(G_{\gamma}(F)/\wt\L(k),\mu_{G_{\gamma}})=\vol(K_{\gamma},\mu_{G_{\gamma}}).
\end{equation}
To show this, let $C=L^{+}\calG^{\circ}_{\gamma}\cap \wt\L$ (where $\calG^{\circ}_{\gamma}$ is the connected N\'eron model of $G_{\gamma}$ whose $\calO_{F}$ points is $K_{\gamma}$, see \S\ref{sss:K gamma}). This is a finite \'etale group over $k$. We have a short exact  sequence of group ind-schemes over $k$
\begin{equation*}
1\to C\to L^{+}\calG^{\circ}_{\gamma}\times\wt\L\to (LG_{\gamma})^{\red}\to 1
\end{equation*}
The associated six term exact sequence for $\Gk$-cohomology gives
\begin{equation*}
1\to C(k)\to K_{\gamma}\times \wt\L(k)\to G_{\gamma}(F)\to \cohog{1}{k,C}\to \ker(\theta)\to1
\end{equation*}
from which we get \eqref{vol id}, using that $\#C(k)=\#\cohog{1}{k,C}$. \qed

\subsection{Examples in $\SL_{2}$}
By Theorem \ref{th:orb coho}, in order to calculate orbital integrals, we need to know not just the geometry of the affine Springer fiber $\sX_{\gamma}$, but also the action of Frobenius on its cohomology. Having already seen many examples of affine Springer fibers over an algebraically closed field in \S\ref{ss:aff examples}, our emphasis here will be on the Frobenius action. 

In this subsection we let $G=\SL_{2}$ and assume $\chark>2$. We will compute several orbital integrals in this case and verify Theorem \ref{th:orb coho} in these cases by explicit calculations.

\subsubsection{Unramified case: $\gamma$}\label{unram sl2} Let $\gamma=\mat{0}{at}{t}{0}\in\sl_{2}(F)$ be a regular semisimple element with $a\in k^{\times}-(k^{\times})^{2}$. 

Let $E\subset\mathfrak{gl}_{2}(F)$ be the centralizer of $\gamma$ in $\mathfrak{gl}_{2}(F)$. Then $E=\{u+v\gamma|u,v\in F\}$ is an unramified quadratic extension of $F$ obtained by adjoining $\sqrt{a}$. Therefore we have $E=k_{E}\lr{t}$, with $k_{E}= k(\sqrt{a})$. We have $G_{\gamma}(F)=(E^{\times})^{\Nm=1}=\ker(\Nm: E^{\times}\to F^{\times})$, which is compact. We fix a Haar measure on $G_{\gamma}(F)$ with total volume 1. 

Let $\sX_{\gamma}$ be the affine Springer fiber of $\gamma$, which is a scheme over $k$. Lemma \ref{l:Ogamma} implies that
\begin{equation*}
O_{\gamma}(\one_{\frg(\calO_{F})})=\#X_{\gamma}=\#\sX_{\gamma}(k).
\end{equation*}
In \S\ref{P1 chain} we have shown that $\sX_{\gamma,\kbar}$ is an infinite union of rational curves $C_{n}\cong\PP^{1}$ indexed by the integers $n\in\ZZ$. Since $\gamma$ is diagonalizable over $E$,  each component $C_{n}$ is in fact defined over $k_{E}=k(\sqrt{a})$. The lattice $\L_{\gamma}\subset G_{\gamma}(F^{\ur})$ is contained in $G_{\gamma}(E)\cong E^{\times}$, and is generated by the uniformizer $t\in E$. We label the components $C_{n}$ so that $t\in \L_{\gamma}$ sends $C_{n}$ to $C_{n+1}$. Let $x_{n+1/2}:=C_{n}\cap C_{n+1}$, which is a $k_{E}$-point of $\sX_{\gamma}$.  

The action of  the nontrivial involution $\sigma\in\Gal(k_{E}/k)$ on $G_{\gamma}(E)$ is by inversion, hence it also acts on $\L_{\gamma}$ by inversion.  The standard lattice $\calO^{2}_{F}$ lies in both $C_{0}$ and $C_{1}$, hence it is the point $x_{1/2}$. Therefore the point $x_{1/2}$ is fixed by $\sigma$ since it is defined over $k$. Since the action of $\sigma$ on $\sX_{\gamma,k_{E}}$ is compatible with its action on $\L_{\gamma}$ (by inversion), the only possibility is that
\begin{equation*}
\sigma(C_{n})=C_{1-n}, \quad \sigma(x_{n+1/2})=x_{-n+1/2}, \quad\forall n\in\ZZ.
\end{equation*}
The action of $\sigma$ can be represented by the picture
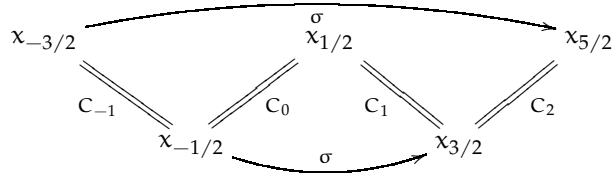
\begin{equation*}
\xymatrix{ x_{-3/2} \ar@/^{1pc}/[rrrr]_{\sigma}
\ar@{=} '[dr]_{C_{-1}} '[rr]_{C_{0}}  '[drrr]_{C_{1}}  [rrrr]_{C_{2}}
&  & x_{1/2}  & & x_{5/2} \\
& x_{-1/2}\ar@/_1pc/[rr]^{\sigma} &  & x_{3/2}}
\end{equation*}
Here each double line represents a $\PP^{1}$. From this we see that $X_{\gamma}=\sX_{\gamma}(k)$ consists of only one point $x_{1/2}$, namely the standard lattice $\calO^{2}_{F}$. This implies that
\begin{equation}\label{Oga}
O_{\gamma}(\one_{\frg(\calO_{F})})=1.
\end{equation}

\subsubsection{Unramified case: $\gamma'$} Now consider the element $\gamma'=\mat{0}{at^{2}}{1}{0}\in\sl_{2}(F)$. In Exercise \ref{exIII:stable conj} we see that $\gamma'$ is stably conjugate to $\gamma$ but not conjugate to $\gamma$ under $\SL_{2}(F)$. However, $\gamma'$ is conjugate to $\gamma$ under $\SL_{2}(E)$. Therefore, the affine Springer fiber $\sX_{\gamma'}$ still looks the same as $\sX_{\gamma}$ over $k_{E}$, but the action of $\sigma\in\Gal(k_{E}/k)$ is different.

Consider the component of $\sX_{\gamma',k_{E}}$ whose $k_{E}$-points consist of $\gamma'$-stable lattices $\L\subset E^{2}$ such that $t\calO_{E}\oplus\calO_{E}\subset \L\subset \calO_{E}\oplus t^{-1}\calO_{E}$. This component is cut out by conditions defined over $k$, so it is stable under $\sigma$, and we call this component $C'_{0}$. We label the other components of $\sX_{\gamma',k_{E}}$ by $C'_{n}$ ($n\in\ZZ$) so that the generator $t\in\L_{\gamma'}$ sends $C'_{n}$ to $C'_{n+1}$. Let $x'_{n+1/2}=C'_{n}\cap C'_{n+1}\in \sX_{\gamma'}(k_{E})$. 
Since the action of $\sigma$ on $\sX_{\gamma,k_{E}}$ is compatible with its action on $\L_{\gamma'}$ (by inversion), the only possibility is that
\begin{equation*}
\sigma(C'_{n})=C'_{-n}, \quad \sigma(x'_{n+1/2})=x'_{-n-1/2},\quad\forall n\in\ZZ.
\end{equation*}
The action of $\sigma$ can be represented by the picture
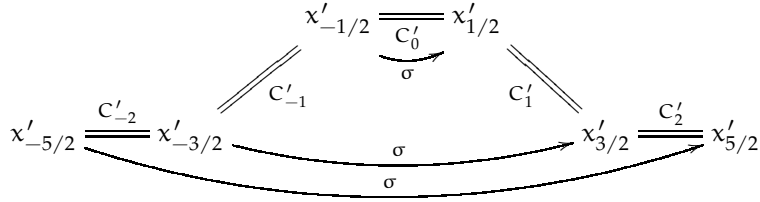
\begin{equation*}
\xymatrix{ & & x'_{-1/2} \ar@/_{1pc}/@<-1ex>[r]_{\sigma}& x'_{1/2}\\
x'_{-5/2}  \ar@/_2pc/[rrrrr]^{\sigma} 
\ar@{=} '[r]^{C'_{-2}} '[urr]_{C'_{-1}}  '[urrr]_{C'_{0}}  '[rrrr]_{C'_{1}} [rrrrr]^{C'_{2}}
&  x'_{-3/2} \ar@/_{1pc}/[rrr]^{\sigma}& & & x'_{3/2}  & x'_{5/2}}
\end{equation*}
Therefore no point $x'_{n+1/2}$ is defined over $k$. The component $C'_{0}$ is the only one that is defined over $k$, and it has to be isomorphic to $\PP^{1}$ over $k$ because it is so over $k_{E}$ (there are no nontrivial Brauer-Severi varieties over a finite field). We see that $\sX_{\gamma'}(k)=C'_{0}(k)\cong\PP^{1}(k)$ has $q+1$ elements. Therefore
\begin{equation}\label{Ogb}
O_{\gamma'}(\one_{\frg(\calO_{F})})=q+1.
\end{equation}
Adding up \eqref{Oga} and \eqref{Ogb} we get
\begin{equation*}
SO_{\gamma}(\one_{\frg(\calO_{F})})=q+2.
\end{equation*}
 
\subsubsection{Unramified case: cohomology} The quotient $\L_{\gamma}\bs\sX_{\gamma}$ is a nodal rational curve obtained from $\PP^{1}$ by glueing two $k$-points into a nodal point. 

Now let us consider the quotient $\L_{\gamma'}\bs\sX_{\gamma'}$. Over $k_{E}$ this is also a nodal rational curve consisting of a unique node $y$ which is image of all $x_{j}$. While $y$ is a $k$-point of the quotient $\L_{\gamma'}\bs\sX_{\gamma'}$, none of its preimages $x_{j}$ are defined over $k$. Therefore, the $q+1$ points in $C'_{0}(k)$ still map injectively to the quotient, in addition to the point $y$. We conclude that $(\L_{\gamma'}\bs\sX_{\gamma'})(k)$ consists of $q+2$ points.  Since $LG_{\gamma'}/\L_{\gamma'}$ is connected, the stable part of the cohomology of $\L_{\gamma'}\bs\sX_{\gamma'}$ is the whole $\cohog{*}{\L_{\gamma'}\bs\sX_{\gamma'}}$, and the alternating sum of Frobenius trace on it is the cardinality of $(\L_{\gamma'}\bs\sX_{\gamma'})(k)$. In this special case we have verified the formula \eqref{d}. We remark that the action of $\sigma$ on the 1-dimensional space $\cohog{1}{\L_{\gamma'}\bs\sX_{\gamma',\kbar},\Ql}$ is by $-1$, and the Grothendieck-Lefschetz trace formula for Frobenius reads $\#(\L_{\gamma'}\bs\sX_{\gamma'})(k)= 1+1+q$ instead of $1-1+q=q$, the latter being the number of $k$-points on a nodal rational curve obtained by identifying two $k$-points on $\PP^{1}$.

\subsubsection{Ramified case: orbital integrals}\label{sss:ram orb SL2} Consider the elements $\gamma=\mat{0}{t^{2}}{t}{0}$ and $\gamma'=\mat{0}{at^{2}}{a^{-1}t}{0}$ where $a\in k^{\times}-(k^{\times})^{2}$. Again $\gamma$ and $\gamma'$ are stably conjugate but not conjugate in $\SL_{2}(F)$. 

Let $E$ be the centralizer of $\gamma$ in $\gl_{2}(F)$. Then $E$ is a ramified quadratic extension of $F$, and $G_{\gamma}(F)=(E^{\times})^{\Nm=1}$. Similarly, let $E'$ be the centralizer of $\gamma'$ in $\gl_{2}(F)$. Then $E'$ is another ramified quadratic extension of $F$, and $G_{\gamma'}(F)=(E'^{\times})^{\Nm=1}$. We choose Haar measures on compact groups $G_{\gamma}(F)$ and $G_{\gamma'}(F)$ with total volume $1$.  

In both cases, $\sX_{\gamma}$ and $\sX_{\gamma'}$ are isomorphic to $\PP^{1}$ as varieties over $k$. In fact we have shown in \S\ref{SL2 ram dim 1} that these varieties are isomorphic to $\PP^{1}$ over $\kbar$, hence they must be isomorphic to $\PP^{1}$ over $k$ as well.  Both $\sX_{\gamma}(k)$ and $\sX_{\gamma'}(k)$  consist of lattices $t\calO_{F}\oplus\calO_{F}\subset \L\subset \calO_{F}\oplus t^{-1}\calO_{F}$, therefore $\sX_{\gamma}=\sX_{\gamma'}$ as subvarieties of $\Gr_{G}$. By Lemma \ref{l:Ogamma}, we have
\begin{equation*}
O_{\gamma}(\one_{\frg(\calO_{F})})=O_{\gamma'}(\one_{\frg(\calO_{F})})=\#\sX_{\gamma}(k)=q+1.
\end{equation*}
Therefore
\begin{equation}\label{ram so}
SO_{\gamma}(\one_{\frg(\calO_{F})})=2(q+1).
\end{equation}

\subsubsection{Ramified case: cohomology} In the setup of \S\ref{sss:ram orb SL2}, $\L_{\gamma}=0$. The component group of $LG_{\gamma}$ is $\ZZ/2\ZZ$, but its action on $\sX_{\gamma}$ is trivial. Therefore the stable part of the cohomology is the whole $\cohog{*}{\sX_{\gamma},\Ql}$, on which the alternating sum of the Frobenius gives the cardinality of $\sX_{\gamma}(k)\cong\PP^{1}(k)$.
However, the parahoric subgroup of $G_{\gamma}(F)=(E^{\times})^{\Nm=1}$ has index $2$ in it ($K_{\gamma}$ consists of those $e\in (E^{\times})^{\Nm=1}$ whose reduction in $k$ is $1$). Therefore, the right side of formula \eqref{d} gets a factor $2=\vol(K_{\gamma},d_{\gamma}g)^{-1}$ in front of $\#\sX_{\gamma}(k)$. This is consistent with \eqref{ram so}, and we have checked the formula \eqref{d} in our special case.

\subsection{Remarks on the Fundamental Lemma} Let us go back to the situation in  \S\ref{unram sl2}. What happens if we take the difference of $O_{\gamma'}(\one_{\frg(\calO_{F})})$ and $O_{\gamma}(\one_{\frg(\calO_{F})})$ instead of their sum? Is there a geometric interpretation of this difference analogous to Theorem \ref{th:orb coho}?

\subsubsection{The $\kappa$-orbital integral} The linear combination $O_{\gamma}(\one_{\frg(\calO_{F})})-O_{\gamma'}(\one_{\frg(\calO_{F})})$ is an example of $\kappa$-orbital integrals. More generally, let $\kappa$ be a character of $\cohog{1}{F,G_{\gamma}}$, then we define the $\kappa$-orbital integral of $\varphi\in \sS(\frg(F))$ to be
\begin{equation*}
O^{\kappa}_{\gamma}(\varphi)=\sum_{\gamma'}\kappa(\inv(\gamma,\gamma'))O_{\gamma'}(\varphi)
\end{equation*}
where the sum is over the $G(F)$-orbits in the stable conjugacy class of $\gamma$.

\subsubsection{Statement of the Fundamental Lemma} The Langlands-Shelstad conjecture, also known as the {\em Fundamental Lemma}, states that the $\kappa$-orbital integral of $\one_{\frg(\calO_{F})}$ for $\gamma\in\frg(F)$ is equal to the stable orbital integral of an element $\gamma_{H}\in\frh(F)$ for a smaller group $H$, up to a simple factor. In formula, the Fundamental Lemma is the identity
\begin{equation*}
O^{\kappa}_{\gamma}(\one_{\frg(\calO_{F})})=\Delta(\gamma,\gamma_{H})\cdot SO_{\gamma_{H}}(\one_{\frh(\calO_{F})}).
\end{equation*}
The smaller group $H$ depends on both $G$ and $\kappa$, and is called the {\em endoscopic group} of $(G,\kappa)$. The number $\Delta(\gamma,\gamma_{H})$ is called the {\em transfer factor}, which turns out to be an integer power of $q$ (depending on $\gamma$ and $\gamma_{H}$) if $\gamma$ is chosen appropriately from its stable conjugacy class, and the measures on $G_{\gamma}$ and $H_{\gamma_{H}}$ are chosen properly.

\subsubsection{A simple case} In the situation of \S\ref{unram sl2}, take the nontrivial character $\kappa$ on $\cohog{1}{F,G_{\gamma}}\cong\ZZ/2\ZZ$, the corresponding endoscopic group $H$ is isomorphic to the torus $G_{\gamma}$; but in general it is not always isomorphic to a subgroup of $G$. The Fundamental Lemma in this case is the identity
\begin{equation*}
O_{\gamma'}(\one_{\frg(\calO_{F})})-O_{\gamma}(\one_{\frg(\calO_{F})})=q=q\cdot SO_{\gamma_{H}}(\one_{\frh(\calO_{F})})
\end{equation*}
where $\gamma_{H}=\gamma$ if we identify $H$ with $G_{\gamma}$. 

On the other hand, in the ramified situation \S\ref{sss:ram orb SL2}, the $\kappa$-orbital integral of $\gamma$ for the nontrivial $\kappa$ vanishes. This maybe explained without calculating the orbital integrals explicitly, for in general, $\kappa$ must factor through a further quotient of $\cohog{1}{F,G_{\gamma}}$ for $O^{\kappa}_{\gamma}(\one_{\frg(\calO_{F})})$ to be possibly nonzero. For the precise statement, see \cite[Prop 8.2.7]{NgoFL}.

\subsubsection{Comments on the proof} The Fundamental Lemma for general $G$ and function field $F$ was established by B-C. Ng\^o \cite{NgoFL}. There is a generalization of Theorem \ref{th:orb coho} to $\kappa$-orbital integrals, in which we replace the stable part of the cohomology of $\L_{\gamma}\bs\sX_{\gamma}$ by the $\kappa$-isotypic part. Using this generalization, Goresky, Kottwitz and MacPherson \cite{GKM} reformulated the Fundamental Lemma as a relation between cohomology groups of the affine Springer fiber $\sX_{\gamma}$ and its endoscopic cousin $\sX^{H}_{\gamma_{H}}$. They were also able to prove the Fundamental Lemma in some special but highly nontrivial cases. Ng\^o's proof builds on this cohomological reformulation, but also uses a new ingredient, namely Hitchin fibers, which can be viewed as a ``global'' analog of affine Springer fibers. This will be the topic of the next lecture.

\subsection{Exercises}

In these exercises, $k=\FF_{q}$ denotes a finite field with $\chark\neq2$, and $F=k\lr{t}$.

\subsubsection{}\label{exIII:K inv} Let $H$ be an algebraic group over $F$ and let $f\in\sS(H(F))$. Show that there exists a compact open subgroup $K\subset H(F)$ such that $f$ is both left and right invariant under $K$.

\subsubsection{}\label{exIII:Ogamma} Prove Lemma \ref{l:Ogamma}.

\subsubsection{}\label{exIII:stable conj} Let $a\in k^{\times}-(k^{\times})^{2}$ and $G=\SL_{2}$. Consider the matrices $\gamma,\gamma'$ as in \S\ref{sss:st conj SL2}. Show that  $\gamma$ and $\gamma'$ are stably conjugate but not conjugate under $\SL_{2}(F)$.

\subsubsection{} Let $G=\GL_{n}$ over $F$ and let $\gamma=\diag(\gamma_{1},\cdots,\gamma_{m})$ be a block diagonal matrix in $\frg(F)^{\rs}$, where $\gamma_{i}\in\mathfrak{gl}_{n_{i}}(F)^{\rs}$. Consider the asymptotic behavior of the orbital integral $O_{\gamma}(\one_{\frg(\calO_{F})})$  as $q=\#k$ tends to $\infty$. Find the smallest integer $d$ such that
\begin{equation*}
O_{\gamma}(\one_{\frg(\calO_{F})})=O(q^{d})
\end{equation*}
as $q\to\infty$. Note that the $O(\cdot)$ on the right side is analysts' $O$ while on the left side it means orbital integral.  Can you interpret $d$ in terms of the characteristic polynomials of the $\gamma_{i}$'s? 

For an explicit estimate of $O_{\gamma}(\one_{\frg(\calO_{F})})$, see \cite{Y-Orb}.

\subsubsection{} Let $G=\SL_{3}$ and $\gamma=\left(\begin{array}{ccc} 0 & 0 & t^{4}\\ 1 & 0 & 0 \\ 0 & 1& 0 \end{array}\right)\in\frg(F)$. Compute $O_{\gamma}(\one_{\frg(\calO_{F})})$. 

{\em Hint}: use the cell decomposition introduced in \S\ref{sss:LS}.

\subsubsection{} Let $G=\SL_{2}$. Let $f$ be the characteristic function of elements $X\in\frg(\calO_{F})$ such that the reduction $\ov{X}$ in $\frg(k)$ is regular nilpotent. Let $\gamma\in\frg(F)$ be a regular semisimple element.  
\begin{enumerate}
\item Show that $O_{\gamma}(f)=0$ unless $\det(\gamma)\in t\calO_{F}$.
\item When $\det(\gamma)\in t\calO_{F}$, show that
\begin{equation*}
O_{\gamma}(f)=O_{\gamma}(\one_{\frg(\calO_{F})})-O_{t^{-1}\gamma}(\one_{\frg(\calO_{F})}).
\end{equation*}
\end{enumerate}

\subsubsection{} Let $G=\GL_{2}$ and $\gamma=\mat{0}{t^{n}}{1}{0}$ for $n\geq 1$ odd. Let $G(F)_{d}$ be the set of $g\in G(F)$ with $\val_{F}(\det g)=d$. Fix the Haar measure on $G(F)$ such that $G(\calO_{F})$ has volume $1$. Show that, for any integer $d\geq0$,
\begin{equation*}
\int_{G(F)_{d}}\one_{\frg(\calO_{F})}(g^{-1}\gamma g)\one_{\calO^{2}_{F}}((0,1)g)dg
\end{equation*}
is the same as the number of closed subschemes $Z$ of the plane curve $y^{2}-t^{n}=0$ satisfying: (1) the underlying topological space of $Z$ is the point $(y,t)=(0,0)$; (2) $\dim_{k}\calO_{Z}=d$.

Note: this exercise relates orbital integrals to Hilbert schemes of points on curves. This relationship has been used in \cite{Y-Orb} to provide an estimate for orbital integrals for $\GL_{n}$. See also \S\ref{sss:MHilb} and \S\ref{sss:Hilb} for an global analog.


\section{Lecture IV: Hitchin fibration} 

During the second half of 1980s, Hitchin introduced the famous integrable system, the moduli space of Higgs bundles, in his study of gauge theory. Around the same time, Kazhdan and Lusztig introduced affine Springer fibers as natural analogs of Springer fibers. For more than 15 years these two objects stayed unrelated until B-C.Ng\^o saw a connection between the two. Ng\^o's fundamental insight can be summarized as saying that Hitchin fibers are global analogs of affine Springer fibers, while affine Springer fibers are local versions of Hitchin fibers. Here ``global'' refers to objects involving a global function field, or an algebraic curve, rather than just a local function field, or a formal disk.  This observation, along with ingenious technical work, allowed Ng\^o to prove the Fundamental Lemma for orbital integrals conjectured by Langlands and Shelstad. We will review Hitchin's integrable system in a slightly more general setting, and make precise its connection to affine Springer fibers.

\subsection{The Hitchin moduli stack}

\subsubsection{The setting} We are back to the setting in \S\ref{ss:notation}. In addition, we fix an algebraic curve $X$ over $k$ (assumed algebraically closed) which is smooth, projective and connected.

\subsubsection{The moduli stack of bundles}
There is a moduli stack $\Bun_{n}$ classifying vector bundles of rank $n$ over $X$.  For any $k$-algebra $R$, $\Bun_{n}(R)$ is the groupoid of rank $n$ vector bundles (locally free coherent sheaves) on $X_{R}:=X\times_{\Spec k}\Spec R$. The stack $\Bun_{n}$ is algebraic, see \cite[Th 4.6.2.1]{LMB}. Moreover, it is smooth and locally of finite type over $k$.

\subsubsection{$G$-torsors}
Recall a {\em (right) $G$-torsor} over $X$ is a scheme $\calE$ over $X$ with a fiberwise action of $G$, such that locally for the \'etale topology of $X$, $\calE$ becomes the $G\times X$ and the $G$-action becomes the right translation action of $G$ on the first factor.

For general reductive $G$, We have the moduli stack $\Bun_{G}$ of $G$-torsors over $X$. For a $k$-algebra $R$, the $R$-points of $\Bun_{G}$ is the groupoid of $G$-torsors over $X_{R}$. Then $\Bun_{G}$ is also a smooth algebraic stack locally of finite type over $k$.

\subsubsection{Associated bundles}\label{sss:ass bun} Let $(V,\rho)$ be a $k$-representation of $G$. Let $\calE$ be a $G$-torsor over $X$. Then there is a vector bundle $\rho(\calE)$ of $X$ whose total space is
\begin{equation*}
\textup{Tot}(\rho(\calE))=G\bs (\calE\times V)
\end{equation*}
where $G$ acts on $(e,v)\in\calE\times V$ by $g\cdot(e,v)=(eg^{-1}, \rho(g)v)$. The vector bundle $\rho(\calE)$ is said to be associated to $\calE$ and $\rho$.

When $G=\GL_{n}$, there is an equivalence of groupoids
\begin{equation}\label{GLn vec}
\{ \mbox{vector bundles $\calV$ of rank $n$ over $X$}\}\cong\{\mbox{$\GL_{n}$-torsors $\calE$ over $X$}\}. 
\end{equation}
The direction $\calV\mapsto\calE$ sends a vector bundle $\calV$ to the $\GL_{n}$-torsor of framings of $\calV$, namely take $\calE=\un{\Isom}_{X}(\calO^{n}_{X},\calV)$, with the natural action of $\GL_{n}$ on the trivial bundle $\calO^{n}_{X}$. The other direction $\calE\mapsto\calV$ sends a $\GL_{n}$-torsor $\calE$ over $X$ to the vector bundle $\St(\calE)$ associated to $\calE$ and the standard representation $\St$ of $\GL_{n}$. The equivalence \eqref{GLn vec} gives a canonical isomorphism of stacks $\Bun_{n}\cong\Bun_{\GL_{n}}$.

\subsubsection{}\label{sss:Bun Sp} For other classical groups $G$, $G$-torsors have more explicit descriptions in terms of vector bundles. For example, when $G=\SL_{n}$, a $G$-torsor over $X$ amounts to the same thing as a pair $(\calV,\iota)$ where $\calV$ is a vector bundle over $X$ of rank $n$ and $\iota$ is an isomorphism of line bundles $\wedge^{n}\calV\cong\calO_{X}$.

When $G=\Sp_{2n}$, the groupoid of $G$-torsors $\calE$ on $X$ is equivalent to the groupoid of pairs $(\calV,\omega)$ where $\calV$ is a vector bundle of rank $2n$ over $X$ and $\omega:\wedge^{2}\calV\to\calO_{X}$ is an $\calO_{X}$-linear map of coherent sheaves that gives a symplectic form on geometric fibers. The map in one direction sends a $\Sp_{2n}$-torsor $\calE$ to the pair $(\calV,\omega)$, where $\calV=\St(\calE)$ is the vector bundle associated to $\calE$ and the standard representation $\St$ of $\Sp_{2n}$, and the symplectic form $\omega$ on $\calV$ comes from the canonical map of $\Sp_{2n}$-representations $\wedge^{2}(\St)\to \one$ (where $\one$ is the trivial representation).

\subsubsection{Higgs bundles} Fix a line bundle $\calL$ over $X$. An {\em $\calL$-twisted Higgs bundle of rank $n$} over $X$  is a pair $(\calV,\varphi)$ where $\calV$ is a vector bundle over $X$ of rank $n$, and $\varphi:\calV\to\calV\otimes\calL$ is an $\calO_{X}$-linear map. The endomorphism $\varphi$ is called a {\em Higgs field} on $\calV$.

There is a moduli stack $\calM_{n,\calL}$ classifying $\calL$-twisted Higgs bundles of rank $n$ over $X$. The morphism $\calM_{n,\calL}\to \Bun_{n}$ forgetting the Higgs field is representable. Therefore $\calM_{n,\calL}$ is also an algebraic stack over $k$.

\subsubsection{$G$-Higgs bundles} An {\em $\calL$-twisted $G$-Higgs bundle} is a pair $(\calE,\varphi)$ where $\calE$ is a $G$-torsor over $X$ and $\varphi$ is a global section of the vector bundle $\Ad(\calE)\otimes\calL$ over $X$. Here, $\Ad(\calE)$ is the vector bundle associated to $\calE$ and the adjoint representation $(\frg,\Ad)$ of $G$, in the sense of \S\ref{sss:ass bun}. We call $\varphi$ an {\em $\calL$-twisted Higgs field} on $\calE$.

When $G=\GL_{n}$, the notion of $\calL$-twisted $G$-Higgs bundle is equivalent to that of an $\calL$-twisted Higgs bundle of rank $n$. In fact, to each $\calL$-twisted $G$-Higgs bundle $(\calE,\varphi)$, we get a Higgs bundle $(\calV=\St(\calE),\phi)$, where $\phi:\calV\to\calV\otimes\calL$ viewed as a global section of $\underline{\End}(\calV)\otimes\calL$ corresponds to $\varphi$ under the canonical isomorphism $\underline{\End}(\calV)\cong\Ad(\calE)$.

We also have the moduli stack $\calM_{G,\calL}$ of $\calL$-twisted Higgs $G$-torsors over $X$. The $R$-points of $\calM_{G,\calL}$ is the groupoid of $\calL_{R}$-twisted $G$-Higgs bundles on $X_{R}$, where $\calL_{R}$ denotes the pullback of $\calL$ to $X_{R}$. The forgetful morphism $\calM_{G,\calL}\to \Bun_{G}$ is representable, hence $\calM_{G,\calL}$ is an algebraic stack over $k$. When $G=\GL_{n}$, we have a canonical isomorphism of stacks $\calM_{n,\calL}\cong\calM_{\GL_{n},\calL}$. 

\subsubsection{Examples} When $G=\SL_{n}$, an $\calL$-twisted $G$-Higgs bundle over $X$ amounts to the same thing as a triple $(\calV,\iota,\varphi)$ where $\calV$ is a vector bundle over $X$ of rank $n$, $\iota:\wedge^{n}\calV\isom\calO_{X}$ and $\varphi: \calV\to\calV\otimes\calL$ satisfies $\Tr(\varphi)=0$.

When $G=\Sp_{2n}$, an $\calL$-twisted $G$-Higgs bundle over $X$ amounts to the same thing as a triple $(\calV,\omega,\varphi)$ where $\calV$ is a vector bundle over $X$ of rank $2n$, $\omega:\wedge^{2}\calV\to\calO_{X}$ is nondegenerate fiberwise,  and $\varphi: \calV\to\calV\otimes\calL$ such that for all local sections $u$ and $v$ of $\calV$,  $\omega(\varphi(u),v)+\omega(u,\varphi(v))=0$ as a local section of $\calL$.

\subsubsection{Hitchin moduli stack as a cotangent bundle}
In Hitchin's original paper \cite{Hi}, he considered the case where $\calL=\omega_{X}$ is the sheaf of $1$-forms on $X$.   This case is particularly important because $\calM_{G,\omega_{X}}$ is closely related to the cotangent bundle of $\Bun_{G}$. For a point $\calE\in\Bun_{G}(R)$ which is a $G$-torsor over $X_{R}$, the cotangent complex of $\Bun_{G}$ at $\calE$ is given by the derived global sections of the complex $\Ad(\calE)^{\vee}\otimes\omega_{X}$ over $X_{R}$. Using a Killing form on $\frg$ we may identify $\Ad(\calE)^{\vee}$ with $\Ad(\calE)$, therefore the cotangent complex of $\Bun_{G}$ at $\calE$ is $\bR\Gamma(X_{R}, \Ad(\calE)\otimes\omega_{X})$, which lives in degrees $0$ and $1$. In particular, when $\calE\in \Bun_{G}(k)$ has finite automorphism group (e.g., $\calE$ is stable), the {\em Zariski cotangent space} at $\calE$ is $\cohog{0}{X,\Ad(\calE)\otimes\omega_{X}}$, i.e., a cotangent vector of $\Bun_{G}$ at $\calE$ is the same thing as a $\omega_{X}$-twisted Higgs field on $\calE$. Therefore, $T^{*}\Bun_{G}$ (properly defined) and $\calM_{G,\omega_{X}}$ share an open substack $T^{*}\Bun^{s}_{G}$, where $\Bun^{s}_{G}$ is the open substack of stable $G$-bundles.

\subsection{Hitchin fibration}

\subsubsection{Hitchin fibration for $\GL_{n}$}
For an $\calL$-twisted Higgs bundle $(\calV,\varphi)$ on $X$, locally on $X$ we may view $\varphi$ as a matrix with entries which are local sections of $\calL$, and we may take the characteristic polynomial of this matrix. The coefficients of this polynomial are well-defined global sections of $\calL^{\otimes i}$, $1\leq i\leq n$. More intrinsically, $\varphi$ induces a map $\wedge^{i}\varphi:\wedge^{i}\calV\to\wedge^{i}\calV\otimes\calL^{\otimes i}$, and we may take 
\begin{equation*}
a_{i}(\varphi):=\Tr(\wedge^{i}\varphi)\in\cohog{0}{X,\calL^{\otimes i}}.
\end{equation*}
This way we have defined a morphism
\begin{equation*}
f: \calM_{n,\calL}\to\calA_{n,\calL}:=\prod_{i=1}^{n}\cohog{0}{X,\calL^{i}}
\end{equation*}
sending $(\calV,\varphi)$ to $(a_{1}(\varphi),\cdots, a_{n}(\varphi))$. We view $\calA_{n,\calL}$ as an affine space over $k$. The morphism $f$ is called the {\em Hitchin fibration} in the case $G=\GL_{n}$.

\subsubsection{Hitchin fibration in general}\label{sss:Hit fib general} For general connected reductive $G$ as in \S\ref{ss:notation}, the coefficients of the characteristic polynomial in the case of $\GL_{n}$ are replaced with the fundamental $G$-invariant polynomials on $\frg$. Recall that $\frc=\frg\sslash G=\Spec\Sym(\frg^{*})^{G}$. Chevalley's theorem says that $\Sym(\frg^{*})^{G}\cong\Sym(\frt^{*})^{W}$, and the latter is a polynomial ring in $r$ variables. We fix homogeneous generators $f_{1},\cdots, f_{r}$ of $\Sym(\frg^{*})^{G}$ as a $k$-algebra, whose degrees  $d_{1}\leq \cdots \leq d_{r}$ are canonically defined although $f_{i}$ are not canonical.  When $G$ is almost simple, the numbers $e_{i}=d_{i}-1$ are the {\em exponents} of $G$. Viewing $f_{i}$ as a symmetric multilinear function $\frg^{\otimes d_{i}}\to k$ invariant under $G$, for any $G$-torsor $\calE$ over $X$, $f_{i}$ induces a map of the associated bundles
\begin{equation*}
f_{i}: \Ad(\calE)^{\otimes d_{i}}\to\calO_{X}.
\end{equation*}
This further induces
\begin{equation*}
f^{\calL}_{i}:(\Ad(\calE)\otimes\calL)^{\otimes d_{i}}\to\calL^{\otimes d_{i}}.
\end{equation*}
If $\varphi$ is an $\calL$-twisted Higgs field on $\calE$,  we may evaluate $f^{\calL}_{i}$ on the section $\varphi^{\otimes d_{i}}$ of $(\Ad(\calE)\otimes\calL)^{\otimes d_{i}}$ to get
\begin{equation*}
a_{i}(\varphi):=f^{\calL}_{i}(\varphi^{\otimes d_{i}})\in\cohog{0}{X,\calL^{\otimes d_{i}}},  \quad i=1,2,\cdots, r.
\end{equation*}
The assignment $(\calE,\varphi)\mapsto(a_{i}(\varphi))_{1\leq i\leq r}$ defines the {\em Hitchin fibration} for $\calM_{G,\calL}$
\begin{equation*}
f=f_{G,\calL}: \calM_{G,\calL}\to \calA_{G,\calL}:=\prod_{i=1}^{r}\cohog{0}{X,\calL^{\otimes d_{i}}}.
\end{equation*}
The target $\calA_{G,\calL}$ is again viewed as an affine space over $k$, and is called the {\em  Hitchin base}.

A more intrinsic way to define the Hitchin base $\calA_{G,\calL}$ is the following. The affine scheme $\frc=\Spec \Sym(\frg^{*})^{G}$ is  equipped with a $\Gm$-action inducing the grading on its coordinate ring. Let $\Tot^{\times}(\calL)\to X$ be the complement of the zero section in the total space of $\calL$. Consider the $\calL$-twisted version of $\frc$ over $X$:
\begin{equation*}
\frc_{X,\calL}:=(\frc\times \Tot^{\times}(\calL))/\Gm
\end{equation*}
where $\l\in\Gm$ acts by $\l: (c,\wt{x})\mapsto (\l c, \l^{-1}\wt{x})$ on the two coordinates. This is an affine space bundle over $X$ whose fibers are isomorphic to $\frc$. Then $\calA_{G,\calL}$ can be canonically identified with the moduli space of sections of the map $\frc_{X,\calL}\to X$. In particular, every point $a\in\calA_{G,\calL}$ gives a map $a:X\to [\frc/\Gm]$.

\subsubsection{The generically regular semisimple locus} Trivializing $\calL$ at the generic point $\eta$ of $X$ and restricting $a_{i}$ to $\eta$, we have a polynomial $P_{a}(y)=y^{n}-a_{1}y^{n-1}+a_{2}y^{n-2}+\cdots +(-1)^{n}a_{n}\in F[y]$, where $F=k(X)$ is the function field of $X$. When $P_{a}(y)$ is a separable polynomial in $F[y]$, we call such an $a$ {\em generically regular semisimple}. The generic regular semisimplicity of $a$ is equivalent to the nonvanishing of the discriminant $\Delta(P_{a})\in\cohog{0}{X,\calL^{n(n-1)}}$ and therefore it defines an open subscheme $\Ah_{n,\calL}\subset \calA_{n,\calL}$. 

For general $G$, viewing $a\in\calA_{G,\calL}$ as a map $a:X\to [\frc/\Gm]$ (see the end of \S\ref{sss:Hit fib general}), we call $a$ {\em generically regular semisimple} if $a$ sends the generic point of $X$ into the open substack $[\frc^{\rs}/\Gm]$. This defines an open subscheme $\Ah_{G,\calL}\subset\calA_{G,\calL}$ generalizing the construction of $\Ah_{n,\calL}$ above.

\subsubsection{Geometric properties}\label{sss:geom Hit} When $\deg \calL>2g-2$, the stack $\calM|_{\Ah_{G,\calL}}$  is smooth, see \cite[Th 4.14.1]{NgoFL}. In this situation, the morphism $f_{G,\calL}$ is flat over $\Ah_{G,\calL}$, see \cite[Cor 4.16.4]{NgoFL}. When $G$ is semisimple, there is a further open dense subset $\calA^{\textup{ani}}_{G,\calL}\subset \Ah_{G,\calL}$ over which $\calM$ is a Deligne-Mumford stack and the map $f_{G,\calL}$ is proper, see \cite[Prop 6.1.3]{NgoFL}. Comparing to the infinite-dimensionality involved in the geometry of affine Springer fibers, the Hitchin fibration has much nicer geometric properties, and yet it is closely related to the affine Springer fibers, as we shall see in \S\ref{ss:rel Hit ASF}.

\subsubsection{Generalization}\label{sss:general Hitchin} Let $H$ be a reductive group over $k$ that fits into an exact sequence of reductive groups
\begin{equation*}
1\to H_{1}\to H\to A\to 1.
\end{equation*}
Let $(V,\rho)$ be a representation of $H$. Then we may consider pairs $(\calE,\psi)$ where $\calE$ is an $H$-torsor over $X$ and $\psi$ is a section of the associated bundle $\rho(\calE)$. Alternatively, such a pair is the same as a morphism $X\to [V/H]$. One can prove that there is an algebraic stack $\calM_{H, \rho}$ classifying such pairs. Every $H$-torsor induces an $A$-torsor, hence we have a morphism $\alpha: \calM_{H, \rho}\to \Bun_{A}$. Fix an $A$-torsor $\calE_{A}$ over $X$, we denote the preimage $\alpha^{-1}(\calE_{A})$ by $\calM_{H,\rho,\calE_{A}}$.

To recover the Hitchin moduli stack, we consider the case $H=G\times\Gm$ with $H_{1}=G$ and $A=\Gm$. Let $V=\frg$ with the action $\rho$ of $H$ defined as follows: $G$ acts by the adjoint representation and $\Gm$ acts by scaling on $V$. An $H$-torsor is a pair consisting of a $G$-torsor and a line bundle $\calL$ on $X$. Fixing the line bundle $\calL$ (which is equivalent to fixing an $A=\Gm$-torsor), we get an isomorphism
\begin{equation*}
\calM_{H,\rho, \calL}\cong\calM_{G, \calL}.
\end{equation*}

For general $(1\to H_{1}\to H\to A\to 1,V,\rho,\calE_{A})$ as above, we may define the analog of the Hitchin base as follows. Let $\frc_{H,\rho}=V\sslash H_{1}$. This is the analog of $\frc$, and it carries an action of $A$. Then we form the twisted version of $\frc_{H,\rho}$ over $X$
\begin{equation*}
\frc_{H,\rho,\calE_{A}}:=\calE_{A}\twtimes{A}\frc_{H,\rho}
\end{equation*}
Then we define $\calA_{H,\rho,\calE_{A}}$ to be the moduli space of sections to the map $\frc_{H,\rho,\calE_{A}}\to X$. The morphism $[V/H]\to [\frc_{H,\rho}/A]$ then induces the analog of the Hitchin fibration
\begin{equation*}
f_{H,\rho,\calE_{A}}: \calM_{H,\rho,\calE_{A}}\to \calA_{H,\rho,\calE_{A}}
\end{equation*}

\subsubsection{Example}\label{sss:MHilb} Consider $H=\GL(U)\times\Gm\times\Gm$, and $V=\End(U)\oplus U^{*}$. The action $\rho(g,s_{1}, s_{2})$ on $V$ is given by $(A,u^{*})\mapsto (s_{1}gAg^{-1}, s_{2}gu^{*})$. The moduli stack $\calM_{H,\rho}$ then maps to $\Pic(X)\times\Pic(X)$ by remembering only the two $\Gm$-torsors. Fixing $(\calL_{1},\calL_{2})\in \Pic(X)\times\Pic(X)$, its  preimage  in $\calM_{H,\rho}$ classifies triples $(\calU,\varphi,\beta)$ where $(\calU,\varphi)$ is a Higgs bundle over $X$ of rank $n=\dim V$, $\beta$ is an $\calO_{X}$-linear map $\calU\to \calL_{2}$. The Hitchin base in this case is the same as the classical Hitchin base $\calA_{n,\calL_{1}}$. Later in \S\ref{sss:Hilb} we will relate this moduli space to Hilbert schemes of curves.

\subsection{Hitchin fibers} An important observation made by Hitchin is that the fibers of the Hitchin fibration $f$ can be described in ``abelian'' terms, namely by line bundles on certain finite coverings of $X$. We elaborate on this observation for $G=\GL_{n}$ and $G=\Sp_{2n}$.

\subsubsection{The case of $\GL_{n}$ and the spectral curve}\label{sss:spectral curve} For $a=(a_{1},\cdots, a_{n})\in\calA_{n,\calL}$, one can define a curve $Y_{a}$ equipped with a degree $n$ morphism $\pi_{a}: Y_{a}\to X$. The construction is as follows. The total space of $\calL$ can be written as a relative spectrum over $X$
\begin{equation*}
\Sigma:=\Tot(\calL)=\Spec(\calO\oplus\calL^{\otimes -1}\oplus\calL^{\otimes -2}\oplus\cdots)=\Spec\Sym(\calO\oplus\calL^{\otimes -1})
\end{equation*}
Let $\pi: \Sigma\to X$ be the projection. Consider the map of coherent sheaves on $X$
\begin{equation*}
\iota_{a}: \calL^{\otimes -n}\to \pi_{*}\calO_{\Sigma}=\calO\oplus\calL^{\otimes -1}\oplus\calL^{\otimes -2}\oplus\cdots
\end{equation*}
given in coordinates by $((-1)^{n}a_{n}, (-1)^{n-1}a_{n-1},\cdots, -a_{1},1,0,\cdots)$. By adjunction $\iota_{a}$ induces a map $\iota'_{a}: \pi^{*}\calL^{\otimes -n}\to \calO_{\Sigma}$, whose image we denote by $\calI_{a}$.  Then $\calI_{a}$ is an ideal sheaf on $\Sigma$. Define the {\em spectral curve} $Y_{a}$ to be the closed subscheme of $\Sigma$ defined by $\calI_{a}$
\begin{equation*}
Y_{a}=\Spec(\calO\oplus\calL^{\otimes-1}\oplus\calL^{\otimes-2}\oplus\cdots)/\calI_{a}
\end{equation*}
If we trivialize $\calL$ on some open subset $U\subset X$, and view $a_{i}$ as functions on $U$, then $Y_{a}|_{U}$ is the subscheme of $U\times\AA^{1}$ defined by one equation $y^{n}-a_{1}y^{n-1}+a_{2}y^{n-2}+\cdots +(-1)^{n}a_{n}=0$ (where $y$ is the coordinate on $\AA^{1}$). The projection $\pi_{a}: Y_{a}\to X$ is finite flat of degree $n$. The curve $Y_{a}$ is called the {\em spectral curve} of $a$ since the fibers of $\pi_{a}$ are the roots of the characteristic polynomial $y^{n}-a_{1}y^{n-1}+a_{2}y^{n-2}+\cdots +(-1)^{n}a_{n}$.

When $a\in\Ah_{n,\calL}$,  the curve $Y_{a}$ is reduced and therefore smooth on a Zariski dense open subset, there is a moduli stack $\cPic(Y_{a})$ classifying torsion-free coherent $\calO_{Y_{a}}$-modules that are generically of rank $1$, see \cite{AK}. The usual Picard stack $\Pic(Y_{a})$ classifying line bundles on $Y_{a}$ is an open substack of $\cPic(Y_{a})$, and it acts on $\cPic(Y_{a})$ by tensoring.

\begin{prop}\label{p:cPic}  Suppose $a\in\Ah_{n,\calL}(k)$. Let $\calM_{a}$ be the fiber of $f:\calM_{n,\calL}\to \calA_{n,\calL}$ over $a$. Then there is a canonical isomorphism of stacks
\begin{equation*}
\cPic(Y_{a})\cong\calM_{a}.
\end{equation*}
\end{prop}

\subsubsection{Sketch of proof}\label{sss:pf cPic}
We give the morphism $\cPic(Y_{a})\to \calM_{a}$. For any coherent sheaf $\calF$ on $Y_{a}$, the direct image $\pi_{a,*}\calF$ is a coherent sheaf on $X$ equipped with a map $\varphi_{\calF}: \pi_{a,*}\calF\otimes\calL^{-1}\to\pi_{a,*}\calF$ because $\calF$ is an $\calO_{Y_{a}}$-module and $\calO_{Y_{a}}$ contains $\calL^{-1}$ as the second direct summand. When $\calF$ is torsion-free and generically rank $1$, $\calV:=\pi_{a,*}\calF$ is torsion-free over $X$ (hence a vector bundle) of rank $n$, and the map $\varphi_{\calF}$ induces a Higgs field $\varphi:\calV\to\calV\otimes\calL$. The assignment $\calF\mapsto (\calV,\varphi)$ defines the morphism $\cPic(Y_{a})\to \calM_{a}$, which can be shown to be an isomorphism. \qed

\subsubsection{The case $G=\Sp_{2n}$} In this case $a\in\calA_{G,\calL}$ is a tuple $(a_{1},\cdots, a_{n})$ with $a_{i}\in\cohog{0}{X,\calL^{\otimes2i}}$. For $a\in\calA_{G,\calL}$, we can similarly define a spectral curve $Y_{a}$ as the closed subscheme of the total space of $\calL$ cut out by the ideal locally generated by
\begin{equation*}
P_{a}(y)=y^{2n}+a_{1}y^{2n-2}+\cdots+a_{n}.
\end{equation*}
Note that $Y_{a}$ carries an involution  $\sigma(y)=-y$ under which the projection $\pi_{a}: Y_{a}\to X$ is invariant. Now suppose $P_{a}(y)$ is separable when restricted to the generic point of $X$, so that $Y_{a}$ is reduced. The involution $\sigma$ induces an involution $\sigma$ on the compactified Picard $\cPic(Y_{a})$. Let $(-)^{\vee}$ be the relative Serre duality functor on coherent sheaves on $Y_{a}$, i.e., $\calF^{\vee}=\un{\Hom}_{\calO_{X}}(\calF,\calO_{X})$ viewed as an $\calO_{Y_{a}}$-module in a natural way. Let $\cPrym(Y_{a};\sigma)$ be the moduli stack of pairs $(\calF,\iota)$ where $\calF\in\cPic(Y_{a})$ and $\iota$ is an isomorphism $\iota: \sigma^{*}\calF\cong\calF^{\vee}$ satisfying that $(\sigma^{*}\iota)^{\vee}=-\iota$. We called $\cPrym(Y_{a};\sigma)$ {\em the compactified Prym stack} of $Y_{a}$ with respect to the involution $\sigma$. Similar to the case of $\GL_{n}$, we have the following description of $\calM_{a}$.

\begin{prop} Suppose $a\in\calA_{G,\calL}$ such that $Y_{a}$ is reduced. Then there is a canonical isomorphism of stacks
\begin{equation*}
\cPrym(Y_{a})\cong\calM_{a}.
\end{equation*}
\end{prop}

\subsubsection{Sketch of proof}
For $a\in\calA_{G,\calL}$, points in $\calM_{a}$ are triples $(\calV,\omega, \varphi)$ where $\calV$ is a vector bundle of rank $2n$ on $X$, $\omega: \wedge^{2}\calV\to\calO_{X}$ is a symplectic form on $\calV$ as in \S\ref{sss:Bun Sp}, and $\varphi:\calV\to\calV\otimes\calL$ satisfies $\om(\varphi u, v)+\om(u, \varphi v)=0$ for local sections $u,v$ of $\calV$, and the characteristic polynomial of $\varphi$ is $P_{a}(y)$. By Proposition \ref{p:cPic}, the Higgs bundle $(\calV,\varphi)$ gives a point $\calF\in\cPic(Y_{a})$. The symplectic form can be viewed as an isomorphism $\jmath: (\calV,\varphi)\isom (\calV^{\vee}, -\varphi^{\vee})$ such that $\jmath^{\vee}=-\jmath$. Note that the Higgs bundle $(\calV^{\vee}, -\varphi^{\vee})$ corresponds to $\sigma^{*}\calF^{\vee}$ under the isomorphism in Proposition \ref{p:cPic}. The isomorphism $\jmath$ then turns into $\iota:  \sigma^{*}\calF\isom\calF^{\vee}$ satisfying $(\sigma^{*}\iota)^{\vee}=-\iota$.  \qed

\subsubsection{Example \ref{sss:MHilb} continued}\label{sss:Hilb}
Fix two line bundles $\calL_{1}$ and $\calL_{2}$ on $X$, and let $\calH_{\calL_{1},\calL_{2}}$ be the fiber of the moduli stack $\calM_{H,\rho}$ in \S\ref{sss:MHilb} over $(\calL_{1},\calL_{2})$. Then $\calH_{\calL_{1},\calL_{2}}$ classifies $(\calV,\varphi,\beta)$ where $\calV$ is a vector bundle of rank $n$ over $X$, $\varphi:\calU\to\calU\otimes\calL_{1}$ is a Higgs field and $\beta: \calU\to \calL_{2}$. We have a Hitchin-type map $\calH_{\calL_{1},\calL_{2}}\to\calA_{n,\calL_{1}}$ sending $(\calV,\varphi,\beta)$ to the $(a_{i}(\varphi))_{1\leq i\leq n}$. Let $a\in\Ah_{n,\calL_{1}}$ and let $\calH_{\calL_{1},\calL_{2},a}$ be the fiber of $\calH_{\calL_{1},\calL_{2}}$ over $a$. Consider the same spectral curve $Y_{a}$ as in \S\ref{sss:spectral curve}. Using Proposition \ref{p:cPic} we may identify $(\calU,\varphi)$ with a point $\calF\in\cPic(Y_{a})$. The map $\beta:\pi_{a,*}\calF=\calU\to\calL_{2}$ gives a map $b:\calF\to \pi^{!}_{a}\calL_{2}$ by adjunction. Here $\pi^{!}_{a}\calL_{2}\cong\pi^{*}_{a}\calL_{2}\otimes\omega_{Y_{a}/X}$, and the relative dualizing complex $\omega_{Y_{a}/X}=\omega_{Y_{a}}\otimes\pi_{a}^{*}\omega_{X}^{-1}$ is a line bundle on $Y_{a}$ because $Y_{a}$ is a planar curve hence Gorenstein. Since $\calF$ is torsion-free and generically a line bundle, $b$ is an injective map of coherent sheaves. Hence the data $(\calU,\varphi,\beta)$ turns into the data of a coherent subsheaf $\calF$ of the line bundle $\pi^{!}_{a}\calL_{2}$. Since $\pi^{!}_{a}\calL_{2}$ is a line bundle, the subsheaf $\calF$ is determined by the support of the quotient $(\pi^{!}_{a}\calL_{2})/\calF$, which is a zero-dimensional subscheme of $Y_{a}$.  We conclude that there is a canonical isomorphism
\begin{equation*}
\calH_{\calL_{1},\calL_{2},a}\cong\Hilb(Y_{a})
\end{equation*}
where $\Hilb(Y_{a})$ is the disjoint union of Hilbert schemes of zero-dimensional subschemes of $Y_{a}$ of various lengths.

To conclude, the Hilbert scheme of points of a family of planar curves can be realized as a Hitchin-type moduli space in the framework of \S\ref{sss:general Hitchin}.

\subsubsection{Symmetry on Hitchin fibers}\label{sss:sym Hit}
In the case $G=\GL_{n}$ and $a\in\Ah_{n,\calL}$, we have identified the Hitchin fiber $\calM_{a}$ with the compactified Picard stack $\cPic(Y_{a})$ of the spectral curve $Y_{a}$. The usual Picard stack $\Pic(Y_{a})$ acts on $\cPic(Y_{a})$ by tensor product. Therefore we have an action of $\Pic(Y_{a})$ on $\calM_{a}$. This action is simply transitive if $Y_{a}$ is smooth.

In the case $G=\Sp_{2n}$ and $a\in\Ah_{G,\calL}$, the Prym stack $\Prym(Y_{a};\sigma)$ is defined as the moduli stack of pairs $(\calL,\iota)$ where $\calL\in\Pic(Y_{a})$ and $\iota:\sigma^*\calL\cong\calL^{-1}$ satisfying $(\sigma^{*}\iota)^{\vee}=\iota$. Then $\Prym(Y_{a};\sigma)$ acts on $\calM_{a}\cong\cPrym(Y_{a};\sigma)$ by tensoring.

For general $G$ and $a\in\Ah_{G,\calL}$, one can similarly define a commutative group stack $\calP_{a}$ that acts on the corresponding Hitchin fiber $\calM_{a}$. In fact, via the map $a:X\to [\frc/\Gm]$, the regular centralizer group scheme $J$ in \S\ref{sss:sym ASF}, which descends to $[\frc/\Gm]$, pulls back to a smooth group scheme $J_{a}$ over $X$ which is generically a torus. The stack $\calP_{a}$ then classifies $J_{a}$-torsors over $X$.

\subsection{Relation with affine Springer fibers}\label{ss:rel Hit ASF}
In this subsection we will state a precise relationship between Hitchin fibers and affine Springer fibers, observed by Ng\^o. 

\subsubsection{} Fix a point $a\in\Ah_{G,\calL}(k)$. Viewing $a$ as a map $X\to [\frc/\Gm]$, let $U_{a}\subset X$ be the preimage of $[\frc^{\rs}/\Gm]$ under $a$. Since $a\in\Ah_{G,\calL}$, $U_{a}$ is non-empty hence the complement $X-U_{a}$ consists of finitely many $k$-points of $X$.  For each point $x\in X-U_{a}$, let $\calO_{x}$ be the completed local ring of $X$ at $x$, with fraction field $F_{x}$ and residue field $k(x)=k$.  We fix a trivialization of $\calL$ near $x$ and we may identify $a_{x}=a|_{\Spec\calO_{x}}$ as an element in $\frc(\calO_{x})\cap\frc^{\rs}(F_{x})$. Let $\ep(a_{x})\in\frg^{\rs}(F_{x})$ be the corresponding point in the Kostant section, then we have the affine Springer fiber $\sX_{a_{x}}:=\sX_{\ep(a_{x})}$ in the affine Grassmannian $\Gr_{G, x}=L_{x}G/L^{+}_{x}G$ (we put $x$ in the subscript to emphasize that the definition of the loop groups uses the field $F_{x}$, which is isomorphic to $k\lr{t}$). The loop group of the centralizer $LG_{\ep(a_{x})}$ acts on $\sX_{a_{x}}$, and the action  factors through the local Picard group $P_{a_{x}}$ as in \S\ref{sss:sym ASF}. On the other hand, we have the action of the global Picard stack $\calP_{a}$ on the Hitchin fiber $\calM_{a}$ mentioned in \S\ref{sss:sym Hit}. The product formula of Ng\^o roughly says that, modulo the actions of the local and global Picard stacks, $\calM_{a}$ is the product of the affine Springer fibers $\sX_{a_{x}}$ for all points $x\in X-U_{a}$.

\begin{theorem}[Product Formula, Ng\^o {\cite[Th 4.6]{NgoHit} and \cite[Prop 4.15.1]{NgoFL}}]\label{th:prod} For $a\in\Ah(k)$, there is a canonical morphism
\begin{equation*}
\calP_{a}\twtimes{\prod_{x\in X-U_{a}}P_{a_{x}}}\left(\prod_{x\in X-U_{a}}\sX_{a_{x}}\right)\to \calM_{a}
\end{equation*}
which is a homeomorphism of stacks. Here the notation $P\twtimes{H}Y$ (where $H$ acts on $P$ on the right and acts on $Y$ on the left) means the quotient of $P\times Y$ by the action of $H$ given by $h\cdot (p,y)=(ph^{-1}, hy)$, $p\in P, y\in Y$ and $h\in H$.
\end{theorem}

In the case $G=\GL_{n}$, the product formula can be reinterpreted in more familiar terms using the compactified Picard stack of spectral curves, which in fact makes sense for all reduced curves. Let $C$ be a reduced and projective curve over $k$. For each point $x\in C(k)$, one can define a local Picard group $P_{x}(C)$ whose $k$ points are $F^{\times}_{x}/\calO^{\times}_{x}$, where $\calO_{x}$ is the completed local ring of $C$ at $x$ and $F_{x}$ its ring of total fractions (which is a product of fields in general). There is also a local analog $\ov{P}_{x}(C)$ of $\cPic(C)$ whose $k$-points are the fractional $\calO_{x}$-ideals of $F_{x}$, compare \S\ref{sss:GLn orb formula}. Then the following variant of the product formula holds, whose proof is similar to that of Theorem \ref{th:prod}.

\begin{prop}\label{p:Pic prod} Let $C$ be a reduced and projective curve over $k$. Let $Z=C-C^{sm}$ be the singular locus of $C$. Then there is a canonical morphism 
\begin{equation*}
\Pic(C)\twtimes{\prod_{x\in Z}P_{x}(C)}\left(\prod_{x\in Z}\ov{P}_{x}(C)\right)\to\cPic(C)
\end{equation*}
which is a homeomorphism of stacks.
\end{prop}

The product formula provides a link between the geometry of Hitchin fibers and that of affine Springer fibers, the latter is closely related to orbital integrals as we have already seen. This link makes it possible to approach the Fundamental Lemma by studying the cohomology of Hitchin fibers. The advantage of using Hitchin fibers instead of affine Springer fibers is that  the Hitchin fibration has nicer geometric properties, as we have seen in \S\ref{sss:geom Hit}.

\subsection{A global version of the Springer action}

The product formula in Theorem \ref{th:prod} suggests that there should be a global analog of affine Springer theory where affine Springer fibers are replaced by Hitchin fibers. Such a theory was developed in a series of papers of the author starting with \cite{GS}.

\subsubsection{Iwahori level structure} We now define a Hitchin-type analog of the affine Springer fibers $\sY_{\gamma}$. We fix the curve $X$ and a line bundle $\calL$ on it as before. Let $\Mpar_{G,\calL}$ be the moduli stack classifying $(x,\calE,\varphi, \calE^{B}_{x})$ where $x\in X$, $(\calE,\varphi)$ is an $\calL$-twisted Higgs $G$-bundle, and $\calE^{B}_{x}$ is a reduction of the fiber $\calE_{x}$ at $x$ to a $B$-torsor $\calE^{B}_{x}$ (here $B\subset G$ is a Borel subgroup of $G$) compatible with $\varphi$. 

When $G=\GL_{n}$, $\Mpar_{G,\calL}$ classifies, in addition to an $\calL$-twisted Higgs bundle $(\calV,\varphi)$, a point $x\in X$ and a full flag of the fiber $\calV_{x}$. Such a full flag is the same data as a chain of coherent subsheaves $\calV_{i}\subset\calV$
\begin{equation*}
\calV(-x)=\calV_{0}\subset \calV_{1}\subset\cdots\subset\calV_{n-1}\subset\calV_{n}=\calV
\end{equation*}
such that $\calV_{i}/\calV_{i-1}$ has length $1$ supported at $x$. The compatibility condition between $\varphi$ and the full flag requires that $\varphi$ restrict to a map $\calV_{i}\to\calV_{i}\otimes\calL$ for each $0\leq i\leq n$. 

We have an analog of the Hitchin fibration
\begin{equation*}
\fpar_{G,\calL}:\Mpar_{G,\calL}\to \calA_{G,\calL}\times X
\end{equation*}
that records also the point $x\in X$ in the data, in addition to $a=f_{G,\calL}(\calE,\varphi)\in\calA_{G,\calL}$.

\begin{theorem}[See {\cite[Th 3.3.3]{GS}}]\label{th:GS} Suppose $\deg\calL>2g-1$, then there is a natural action of $\tilW=\xcoch(T)\rtimes W$ on the restriction of the complex $\bR\fpar_{!}\Ql$ to a certain open subset $(\Ah_{G,\calL}\times X)'$ of $\Ah_{G,\calL}\times X$. \footnote{When $\textup{char}(k)=0$, one can take $(\Ah_{G,\calL}\times X)'=\Ah_{G,\calL}\times X$; in general, the restriction to the open subset $(\Ah_{G,\calL}\times X)'$ does not limit applications to questions about affine Springer fibers.}
\end{theorem}

\subsubsection{Extended symmetry} Just as in the case of affine Springer fibers, the $\tilW$-action in Theorem \ref{th:GS} may be extended to an action of the graded version of the double affine Hecke algebra, see \cite[Th 6.1.6]{GS}. Also, there is a product formula relating the fiber $\Mpar_{a,x}$ over $(a,x)\in\Ah_{G,\calL}\times X$ and the product of the affine Springer fiber $\sY_{a_{x}}$ with affine Springer fibers $\sX_{a_{y}}$ for $y\neq x$.  The induced $\tilW$-action on the stalk $(\bR\fpar_{!}\Ql)_{a,x}\cong\cohog{*}{\Mpar_{a,x}}$ is compatible with the affine Springer action on $\cohog{*}{\sY_{a_{x}}}$ in Theorem \ref{th:aff action}. This connection between $\tilW$-actions on the cohomology of Hitchin fibers and affine Springer fibers can be used to prove results about affine Springer actions. See \cite{YSph} for such an application.

\subsection{Exercises}

\subsubsection{} Suppose $G=\SL_{n}$. Describe the Hitchin fibers over $\Ah_{G,\calL}$ in terms of spectral curves.

\subsubsection{} Suppose $G=\SL_{n}$. Compute the dimension of $\calA_{G,\calL}$ and of a Hitchin fiber $\calM_{a}$. When $\calL=\omega_{X}$, check that they have the same dimension. This is a numerical evidence that the Hitchin fibration in this case is a Lagrangian fibration.

\subsubsection{} Describe the Hitchin base for $G=\SO_{n}$.

\subsubsection{} For $G=\SO_{n}$, describe Hitchin fibers in terms of spectral curves.

\subsubsection{} Let $C$ be a rational curve over $k$ with a unique singularity $x_{0}$ which is unibranched (i.e., the preimage of $x_{0}$ in the normalization $\tilC\cong\PP^{1}$ is a single point).  Let $\overline{P}_{x_{0}}$ be the moduli space of fractional ideals for the completed local ring $\wh\calO_{C,x_{0}}$. Show that there is a canonical homeomorphism
\begin{equation*}
\overline{P}_{x_{0}}\to \cPic(C).
\end{equation*}
Explain why this is a special case of Proposition \ref{p:Pic prod}.

\end{document}